\documentclass[12pt,a4paper,draft]{article}
\usepackage{amsmath,mathrsfs,amssymb,amsthm,amscd}
\usepackage[]{fontenc}
\usepackage{xy} 
\usepackage{enumerate}
\xyoption{all}

\numberwithin{equation}{section}
\newcommand{\nnpar}[1]{\vskip 10pt \noindent {\scshape #1}}\nopagebreak%
\theoremstyle{plain}
\newtheorem{theorem}[equation]{Theorem}
\newtheorem{corollary}[equation]{Corollary}
\newtheorem{proposition}[equation]{Proposition}
\newtheorem{lemma}[equation]{Lemma}

\theoremstyle{definition}
\newtheorem{definition}[equation]{Definition}
\newtheorem{notation}[equation]{Notation}

\newtheorem{remark}[equation]{Remark}

\DeclareMathOperator{\Lie}{Lie}

\DeclareMathOperator{\tr}{tr}

\DeclareMathOperator{\cha}{\widehat{CH}}
\DeclareMathOperator{\ca}{\widehat{c}}
\DeclareMathOperator{\CH}{CH}

\DeclareMathOperator{\Ka}{\widehat{K}_{0}}
\DeclareMathOperator{\Ko}{K_{0}}
\DeclareMathOperator{\Kl}{K_{1}}

\DeclareMathOperator{\rata}{\widehat{Rat}}

\DeclareMathOperator{\diag}{diag}
\DeclareMathOperator{\dd}{d}
\DeclareMathOperator{\Img}{Im}
\DeclareMathOperator{\Int}{Int}

\DeclareMathOperator{\ch}{ch}
\DeclareMathOperator{\cl}{cl}

\DeclareMathOperator{\Ad}{Ad}
\DeclareMathOperator{\ad}{ad}

\DeclareMathOperator{\GL}{GL}

\DeclareMathOperator{\Hom}{Hom}

\DeclareMathOperator{\Id}{id}

\DeclareMathOperator{\Spec}{Spec}
\DeclareMathOperator{\Stab}{Stab}
\DeclareMathOperator{\supp}{supp}

\DeclareMathOperator{\amap}{a}

\DeclareMathOperator{\Res}{Res}

\DeclareMathOperator{\vmap}{v}

\DeclareMathOperator{\im}{Im\,}
\renewcommand{\Im}{{\im}}

\newcommand{\an}{\text{{\rm an}}}

\newcommand{\wlg}{\text{{\rm pre}}}

\newcommand{\growth}{\text{{\rm gth}}}

\newcommand{\GS}{\text{{\rm GS}}}

\newcommand{\QQ}{{\mathbb Q}}

\newcommand{\DD}{{D}}

\newcommand{\cc}{{\mathcal{C}}}

\newcommand{\lgi}{{\text{\rm l,ll}}}
\newcommand{\as}{{\text{\rm l,ll,a}}}
\newcommand{\pre}{{\textrm{pre}}}

\def\?{\ ???\ \immediate\write16{}%
\immediate\write16{Warning: There was still a question mark . . . }%
\immediate\write16{}}

\begin{document}
\setcounter{tocdepth}{2}
\setcounter{section}{-1}

\begin{titlepage}

\begin{center}

\vspace*{4cm}
{\Huge\bf Arithmetic characteristic}\\[1em]
{\Huge\bf  classes}\\[1em]
{\Huge\bf of automorphic vector bundles}

\vspace*{2cm}
{\LARGE J. I. Burgos Gil{\footnotemark
\footnotetext{Partially supported by Grant DGI BFM2000-0799-C02-01 and
  DGI BFM2003-02914}
}, J. Kramer, U. K\"uhn}

\end{center}

\end{titlepage}

\thispagestyle{empty}

\newpage
 
\phantom{aaa}

\begin{abstract}
  We develop a theory of arithmetic characteristic classes of (fully
  decomposed) automorphic vector bundles equipped with an invariant
  hermitian metric. These characteristic classes have values in an
  arithmetic Chow ring constructed by means of differential forms with
  certain log-log type singularities. We first study the cohomological
  properties of log-log differential forms, prove a Poincar\'e lemma
  for them and construct the corresponding arithmetic Chow groups.
  Then we introduce the notion of log-singular hermitian vector
  bundles, which is a variant of the good hermitian vector bundles
  introduced by Mumford, and we develop the theory of arithmetic
  characteristic classes. Finally we prove that the hermitian metrics
  of automorphic vector bundles considered by Mumford are not only
  good but also log-singular. The theory presented here provides the
  theoretical background which is required in the formulation of the
  conjectures of Maillot-Roessler in the semi-abelian case and which
  is needed to extend Kudla's program about arithmetic intersections
  on Shimura varieties to the non compact case.
\end{abstract}

\tableofcontents

\section{Introduction}
\label{sec:introduction}

\nnpar{The main goal.} The main purpose of this article is to extend
the arithmetic intersection theory and the theory of arithmetic
characteristic classes \`a la Gillet, Soul\'e to the category of
(fully decomposed) automorphic vector bundles equipped with the
natural equivariant hermitian metric on Shimura varieties of 
non-compact type. In order to achieve our main goal, an extension
of the formalism by Gillet, Soul\'e taking into account vector
bundles equipped with hermitian metrics allowing a certain type 
of singularities has to be provided. The main prerequisite for the 
present work is the article \cite{BurgosKramerKuehn:cacg}, where 
the foundations of cohomological arithmetic Chow groups are given. 
Before continuing to explain our main results and the outline of 
the paper below, let us fix some basic notations for the sequel.

Let $B$ denote a bounded, hermitian, symmetric domain. By definition, 
$B=G/K$, where $G$ is a semi-simple adjoint group and $K$ a maximal 
compact subgroup of $G$ with non-discrete center. Let $\Gamma$ be a 
neat arithmetic subgroup of $G$; it acts properly discontinuously
and fixed-point free on $B$. The quotient space $X=\Gamma\backslash 
B$ has the structure of a smooth, quasi-projective, complex variety.
The complexification $G_{\mathbb{C}}$ of $G$ yields the compact dual
$\check{B}$ of $B$ given by $\check{B}=G_{\mathbb{C}}/P_{+}\cdot K_
{\mathbb{C}}$, where $P_{+}$ is a suitable parabolic subgroup of $G$
provided by the Cartan decomposition of ${\rm Lie}(G)$. Every
$G_{\mathbb{C}}$-equivariant holomorphic vector bundle $\check{E}$ 
on $\check{B}$ defines a holomorphic vector bundle $E$ on $X$; $E$
is called an \emph{automorphic vector bundle}. An automorphic vector
bundle $E$ is called \emph{fully decomposed}, if $E=E_{\sigma}$ is 
associated to a representation $\sigma:P_{+}\cdot K_{\mathbb{C}} 
\longrightarrow{\rm GL}_{n}(\mathbb{C})$, which is trivial on $P_{+}$.  
Since $K$ is compact, every fully decomposed automorphic vector bundle 
$E$ admits a $G$-equivariant hermitian metric $h$.

Let us recall the following basic example. Let $\pi:\mathcal{B}_{g}^
{(N)}\longrightarrow\mathcal{A}_{g}^{(N)}$ denote the universal abelian 
variety over the moduli space of principally polarized abelian varieties 
of dimension $g$ with a level-$N$ structure ($N\ge 3$); let $e:\mathcal
{A}_{g}^{(N)}\longrightarrow\mathcal{B}_{g}^{(N)}$ be the zero section, 
and $\Omega=\Omega^{1}_{\mathcal{B}_{g}^{(N)}/\mathcal{A}_{g}^{(N)}}$
the relative cotangent bundle. The pull-back $e^{*}\Omega$ defines an
automorphic vector bundle on $\mathcal{A}_{g}^{(N)}$, which is equipped
with a natural hermitian metric $h$. Another example of an automorphic
vector bundle on $\mathcal{A}_{g}^{(N)}$ is the Hodge bundle $\omega$
obtained as the determinant line bundle $\omega=\det(e^{*}\Omega)$;
the corresponding hermitian automorphic line bundle $(\det(e^{*}\Omega),
\det(h))$ is denoted by $\overline{\omega}$.

\nnpar{Background results.} Let $(E,h)$ be a automorphic hermitian 
vector bundle on $X=\Gamma\backslash B$, and $\overline{X}$ a smooth 
toroidal compactification of $X$. In \cite{Mumford:Hptncc}, D.~Mumford 
has shown that the automorphic vector bundle $E$ admits a canonical 
extension $E_{1}$ to $\overline{X}$ characterized by a suitable extension
of the hermitian metric $h$ to $E_{1}$. However, the extension of 
$h$ to $E_{1}$ is no longer a smooth hermitian metric, but inherits 
singularities of a certain type. On the other hand, it is remarkable 
that this extended hermitian metric behaves in many aspects like a
smooth hermitian metric. In this respect, we will now discuss various
definitions which were made in the past in order to extract basic 
properties for these extended hermitian metrics.

In \cite{Mumford:Hptncc}, D.~Mumford introduced the concept of 
good forms and good hermitian metrics. The good forms are differential 
forms, which are smooth on the complement of a normal crossing divisor 
and have certain singularities along this normal crossing divisor; 
the singularities are modeled by the singularities of the Poincar\'e 
metric. The good forms have the property of being locally integrable 
with zero residue. Therefore, they define currents, and the map from 
the complex of good forms to the complex of currents is a morphism of
complexes. The good hermitian metrics are again smooth hermitian metrics
on the complement of a normal crossing divisor and have logarithmic 
singularities along the divisor in question. Moreover, the entries of 
the associated connection matrix are good forms. The Chern forms of 
good hermitian vector bundles, i.e., of vector bundles equipped with 
good hermitian metrics, are good forms, and the associated currents
represent the Chern classes in cohomology. Thus, in this sense, the 
good hermitian metrics behave like smooth hermitian metrics. In the 
same paper, D.~Mumford proves that automorphic hermitian vector bundles 
are good hermitian vector bundles. 

In \cite{Faltings:EaVZ}, G.~Faltings introduced the concept of a
hermitian metric on line bundles with logarithmic singularities along 
a closed subvariety. He showed that the heights associated to line 
bundles equipped with singular hermitian metrics of this type have 
the same finiteness properties as the heights associated to line 
bundles equipped with smooth hermitian metrics. The Hodge bundle 
$\overline{\omega}$ on $\mathcal{A}_{g}^{(N)}$ equipped with the 
Petersson metric provides a prominent example of such a hermitian 
line bundle; it plays a crucial role in Faltings' proof of the 
Mordell conjecture. Recall that the height of an abelian variety 
$A$ with respect to $\overline{\omega}$ is referred to as the Faltings 
height of $A$. It is a remarkable fact that if $A$ has complex 
multiplication of abelian type, then its Faltings height is essentially 
given by a special value of the logarithmic derivative of a Dirichlet
$L$-series. It is conjectured by P.~Colmez that in the general case 
the Faltings height is essentially given by a special value of the 
logarithmic derivative of an Artin $L$-series.

In \cite{Kuehn:gainc}, the third author introduced the concept 
of logarithmically singular hermitian line bundles on arithmetic
surfaces. He provided an extension of arithmetic intersection theory 
(on arithmetic surfaces) suited for such logarithmically singular 
hermitian line bundles. The prototype of such a line bundle is the
automorphic hermitian line bundle $\overline{\omega}$ on the modular
curve $\mathcal{A}_{1}^{(N)}$. U.~K\"uhn calculated its arithmetic 
self-intersection $\overline{\omega}^{2}$ number to 
\begin{displaymath}
\overline{\omega}^{2}=d_{N}\cdot\zeta_{\mathbb{Q}}(-1)\left(\frac
{\zeta'_{\mathbb{Q}}(-1)}{\zeta_{\mathbb{Q}}(-1)}+\frac{1}{2}\right); 
\end{displaymath}
here $\zeta_{\mathbb{Q}}(s)$ denotes the Riemann zeta function and
$d_{N}$ equals the degree of the classifying morphism of $\mathcal
{A}_{1}^{(N)}$ to the coarse moduli space $\mathcal{A}_{1}^{(1)}$.

In \cite{BurgosKramerKuehn:cacg}, an abstract formalism was developed,
which allows to associate to an arithmetic variety $\mathcal{X}$ 
arithmetic Chow groups $\cha^{*}(\mathcal{X},\mathcal{C})$ with 
respect to a cohomological complex $\mathcal{C}$ of a certain type. 
This formalism is an abstract version of the arithmetic Chow groups 
introduced in \cite{Burgos:CDB}. In \cite{BurgosKramerKuehn:cacg}, 
the arithmetic Chow ring $\cha^{*}(\mathcal{X},\mathcal{D}_{\pre})_
{\mathbb{Q}}$ was introduced, where the cohomological complex $\mathcal
{D}_{\pre}$ in question is built from pre-log and pre-log-log differential 
forms. This ring allows us to define arithmetic self-intersection 
numbers of automorphic hermitian line bundles on arithmetic varieties 
associated to $X=\Gamma\backslash B$. It is expected that these
arithmetic self-intersection numbers play an important role for 
possible extensions of the Gross-Zagier theorem to higher dimensions 
(cf. conjectures of S.~Kudla).

In \cite{BruinierBurgosKuehn:bpaihs}, J.~Bruinier, J.~Burgos, and 
U.~K\"uhn used the theory developed in \cite{BurgosKramerKuehn:cacg}  
in order to obtain an arithmetic generalization of the Hirzebruch-Zagier
theorem on the generating series for cycles on Hilbert modular varieties. 
Recalling that Hilbert modular varieties parameterize abelian surfaces
with multiplication by the ring of integers $\mathcal{O}_{K}$ of a real
quadratic field $K$, a major result in \cite{BruinierBurgosKuehn:bpaihs} 
is the following formula for the arithmetic the self-intersection number 
of the automorphic hermitian line bundle $\overline{\omega}$ on the 
moduli space of abelian surfaces with real multiplication by $\mathcal
{O}_{K}$ with a fixed level-$N$ structure
\begin{displaymath}
\overline{\omega}^{3}=-d_{N}\cdot\zeta_{K}(-1)\left(\frac{\zeta_{K}'
(-1)}{\zeta_{K}(-1)}+\frac{\zeta'_{\mathbb{Q}}(-1)}{\zeta_{\mathbb{Q}}
(-1)}+\frac{3}{2}+\frac{1}{2}\log(D_{K})\right);
\end{displaymath}
here $D_{K}$ is the discriminant of $\mathcal{O}_{K}$, $\zeta_{K}(s)$ 
is the Dedekind zeta function of $K$, and, as above, $d_{N}$ is the 
degree of the classifying morphism obtained by forgetting the level-$N$
structure.  

As another application of the formalism developed in \cite
{BurgosKramerKuehn:cacg}, we derived a height pairing with respect 
to singular hermitian line bundles for cycles in any codimension. 
Recently, G.~Freixas \cite{Freixas:_singul_heigh} has proved 
finiteness results for our height pairing, thus generalizing both 
Faltings' results mentioned above and the finiteness results of 
J.-B.~Bost, H.~Gillet and C.~Soul\'e \cite{BostGilletSoule:HpvpGf} 
in the smooth case.

The main achievement of the present paper is to give constructions
of arithmetic intersection theories, which are suited to deal with
all of the above vector bundles equipped with hermitian metrics
having singularities of a certain type such as the automorphic
hermitian vector bundles on Shimura varieties of non-compact type. 

For a perspective view of applications of the theory developed 
here we refer to the conjectures of V.~Maillot and D.~Roessler
\cite{MaillotRoessler:cdl}, K.~K\"ohler \cite{Koehler:hppat} and 
the program due to S.~Kudla \cite{Kudla:desgf}, \cite{Kudla:icm},
\cite{kudla:msri}.

\nnpar{Arithmetic characteristic classes.} 
We recall from \cite{Soule:lag} that the arithmetic $K$-group $\Ka
(\mathcal{X})$ of an arithmetic variety $\mathcal{X}$ \`a la Gillet,
Soul\'e is defined as the free group of pairs $(\overline{E},\eta)$ 
of a hermitian vector bundle $\overline{E}$ and a smooth differential 
form $\eta$ modulo the relation    
\begin{displaymath}
(\overline{S},\eta')+(\overline{Q},\eta'')=(\overline{E},\eta'+\eta''+
\widetilde{\ch}(\overline{\mathcal{E}})),
\end{displaymath}
for every short exact sequence of vector bundles (equipped with 
arbitrary smooth hermitian metrics)
\begin{displaymath}
\overline{\mathcal{E}}:\,0\longrightarrow\overline{S}\longrightarrow
\overline{E}\longrightarrow\overline{Q}\longrightarrow 0,
\end{displaymath} 
and for any smooth differential forms $\eta',\eta''$; here $\widetilde
{\ch}(\overline{\mathcal{E}})$ denotes the (secondary) Bott-Chern form 
of $\overline{\mathcal{E}}$. 

In \cite{Soule:lag}, H.~Gillet and C.~Soul\'e attached to the elements 
of $\Ka(\mathcal{X})$, represented by hermitian vector bundles $\overline
{E}=(E,h)$, arithmetic characteristic classes $\widehat{\phi}(\overline
{E})$, which lie in the ``classical'' arithmetic Chow ring $\cha^{*}
(\mathcal{X})_{\mathbb{Q}}$. A particular example of such an arithmetic 
characteristic class is the arithmetic Chern character $\widehat{\ch}
(\overline{E})$, whose definition also involves the Bott-Chern form
$\widetilde{\ch}(\overline{\mathcal{E}})$.

In order to be able to carry over the concept of arithmetic characteristic
classes to the category of vector bundles $E$ over an arithmetic variety
$\mathcal{X}$ equipped with a hermitian metric $h$ having singularities
of the type considered in this paper, we proceed as follows: Letting $h_
{0}$ denote an arbitrary smooth hermitian metric on $E$, we have the 
obvious short exact sequence of vector bundles
\begin{displaymath}
\overline{\mathcal{E}}:\,0\longrightarrow 0\longrightarrow(E,h)
\longrightarrow(E,h_{0})\longrightarrow 0,
\end{displaymath}
to which is attached the Bott-Chern form $\widetilde{\phi}(\overline
{\mathcal{E}})$ being no longer smooth, but having certain singularities.
Formally, we then set
\begin{displaymath}
\widehat{\phi}(E,h):=\widehat{\phi}(E,h_{0})+\amap\left(\widetilde
{\phi}(\overline{\mathcal{E}})\right),
\end{displaymath}
where $\amap$ is the morphism mapping differential forms into arithmetic
Chow groups. In order to give meaning to this definition, we need to know 
the singularities of $\widetilde{\phi}(\overline{\mathcal{E}})$; moreover,
we have to show the independence of the (arbitrarily chosen) smooth 
hermitian metric $h_{0}$.

Once we can control the singularities of $\widetilde{\phi}(\overline
{\mathcal{E}})$, the abstract formalism developed in \cite
{BurgosKramerKuehn:cacg} reduces our task to find a cohomological 
complex $\mathcal{C}$, which contains the elements $\widetilde{\phi}
(\overline{\mathcal{E}})$, and has all the properties we desire for
a reasonable arithmetic intersection theory. Once the complex 
$\mathcal{C}$ is constructed, we obtain an arithmetic $K$-theory 
with properties depending on the complex $\mathcal{C}$, of course.

The most naive way to construct an arithmetic intersection theory 
for automorphic hermitian vector bundles would be to only work with
good forms and good metrics. This procedure is doomed to failure 
for the following two reasons: First, the complex of good forms is
not a Dolbeault complex. However, this first problem can be easily
solved by imposing that it is also closed under the differential
operators $\partial$, $\bar{\partial}$, and $\partial\bar{\partial}$.
The second problem is that the complex of good forms is not big 
enough to contain the singular Bott-Chern forms which occur. For 
example, if $\mathcal{L}$ is a line bundle, $h_{0}$ a smooth hermitian
metric and $h$ a singular metric, which is good along a divisor $D$ 
(locally, in some open coordinate neighborhood, given by the equation 
$z=0$), the Bott-Chern form (associated to the first Chern class)
$\widetilde{{\rm c}}_{1}(\mathcal{L};h,h_{0})$ encoding the change 
of metrics grows like $\log\log(1/|z|)$, whereas the good functions 
are bounded. 

The solution of these problems led us to consider the $\mathcal{D}_
{\log}$-complexes $\mathcal{D}_{\pre}$ made by pre-log and pre-log-log 
forms and its subcomplex $\mathcal{D}_{l,ll}$ consisting of log and 
log-log forms. We emphasize that neither the complex of good forms 
nor the complex of pre-log-log forms are contained in each other.
We also note that if one is interested in arithmetic intersection 
numbers, the results obtained by both theories agree.

\nnpar{Discussion of results.} 
The $\mathcal{D}_{\log}$-complex $\mathcal{D}_{\pre}$ made out of
pre-log and pre-log-log forms could be seen as the complex that 
satisfies the minimal requirements needed to allow log-log singularities
along a fixed divisor as well as to have an arithmetic intersection 
theory with arithmetic intersection numbers in the proper case (see 
\cite{BurgosKramerKuehn:cacg}). As we will show in theorem \ref
{thm:secondary-good}, the Bott-Chern forms associated to the change 
of metrics between a smooth hermitian metric and a good metric belong 
to the complex of pre-log-log forms. Therefore, we can define arithmetic
characteristic classes of good hermitian vector bundles in the 
arithmetic Chow groups with pre-log-log forms. If our arithmetic 
variety is proper, then we can use this theory to calculate arithmetic 
Chern numbers of automorphic hermitian vector bundles of arbitrary 
rank. However, the main disadvantage of $\mathcal{D}_{\pre}$ is that 
we do not know the size of the associated cohomology groups.

The $\mathcal{D}_{\log}$-complex $\mathcal{D}_{l,ll}$ made out of
log and log-log forms is a subcomplex of $\mathcal{D}_{\pre}$. The
main difference is that \emph{all} the derivatives of the component 
functions of the log and log-log forms have to be bounded, which 
allows us to use an inductive argument to prove a Poincar\'e lemma, 
which implies that the associated Deligne complex computes the usual
Deligne-Beilinson cohomology (see theorem \ref{thm:fq}). For this
reason we have better understanding of the arithmetic Chow groups 
with log-log forms (see theorem \ref{thm:16}).  

Since a good form is in general not a log-log form, it is not true 
that the Chern forms of a good hermitian vector bundle are log-log
forms. Hence, we introduce the notion of log-singular hermitian 
metrics, which have, roughly speaking, the same relation to log-log 
forms as the good hermitian metrics to good forms. We then show that 
the Bott-Chern forms associated to the change of metrics between 
smooth hermitian metrics and log-singular hermitian metrics are 
log-log forms. As a consequence, we can define the Bott-Chern forms 
for short exact sequences of vector bundles equipped with log-singular 
hermitian metrics. These Bott-Chern forms have an axiomatic characterization 
similar to the Bott-Chern forms for short exact sequences of vector 
bundles equipped with smooth hermitian metrics. The Bott-Chern forms 
are the main ingredients in order to extend the theory of arithmetic
characteristic classes to log-singular hermitian vector bundles.  

The price we have to pay in order to use log-log forms is that 
it is more difficult to prove that a particular form is log-log: 
we have to bound all derivatives. Note however that most pre-log-log 
forms which appear are also log-log forms (see for instance section
\ref{sec:shim-vari}). On the other hand, we point out that the 
theory of log-singular hermitian vector bundles is not optimal for 
several other reasons. The most important one is that it is not 
closed under taking sub-objects, quotients and extensions. For
example, let
\begin{displaymath}
0\longrightarrow (E',h')\longrightarrow(E,h)\longrightarrow(E'',h'') 
\longrightarrow 0
\end{displaymath}
be a short exact sequence of hermitian vector bundles such that 
the metrics $h'$ and $h''$ are induced by $h$. Then, the assumption 
that $h$ is a log-singular hermitian metric does not imply that
the hermitian metrics $h'$ and $h''$ are log-singular, and vice 
versa. In particular, automorphic hermitian vector bundles that 
are not fully decomposed can always be written as successive extensions 
of fully decomposed automorphic hermitian vector bundles, whose metrics
are in general not log-singular. A related question is that the 
hermitian metric of a unipotent variation of polarized Hodge structures 
induced by the polarization is in general not log-singular. These
considerations suggest that one should further enlarge the notion 
of log-singular hermitian metrics. 

Since the hermitian vector bundles defined on a quasi-projective
variety may have arbitrary singularities at infinity, we also 
consider differential forms with arbitrary singularities along 
a normal crossing divisor. Using these kinds of differential 
forms we are able to recover the arithmetic Chow groups \`a la 
Gillet, Soul\'e for quasi-projective varieties. 

\nnpar{Outline of paper.} The set-up of the paper is as follows. 
In section \ref{sec:loglog}, we introduce several complexes of 
singular differential forms and discuss their relationship. Of 
particular importance are the complexes of log and log-log forms 
for which we prove a Poincar\'e lemma allowing us to characterize 
their cohomology by means of their Hodge filtration. In section
\ref{sec:arithmetic-chow-ring-log-log}, we introduce and study
arithmetic Chow groups with differential forms which are log-log 
along a fixed normal crossing divisor $D$. We also consider 
differential forms having arbitrary singularities at infinity; 
in particular, we prove that for $D$ being the empty set, the
arithmetic Chow groups defined by Gillet, Soul\'e are recovered. 
In section \ref{sec:char-class}, we discuss several classes of 
singular hermitian metrics; we prove that the Bott-Chern forms 
associated to the change of metrics between a smooth hermitian 
metric and a log-singular hermitian metric are log-log forms. 
We also show that the Bott-Chern forms associated to the change 
of metrics between a smooth hermitian metric and a good hermitian 
metric are pre-log-log. This allows us to define arithmetic 
characteristic classes of log-singular hermitian vector bundles. 
Finally, in section \ref{sec:shim-vari}, after having given a 
brief recollection of the basics of Shimura varieties, we prove 
that the fully decomposed automorphic vector bundles equipped 
with an equivariant hermitian metric are log-singular hermitian 
vector bundles. In this respect many examples are provided to 
which the theory developed in this paper can be applied.

\vspace{5mm}

\noindent
\emph{Acknowledgements:} In course of preparing this manuscript,
we had many stimulating discussions with many colleagues. We
would like to thank them all. In particular, we would like to
express our gratitude to J.-B.~Bost, J.~Bruinier, P.~Guillen, 
W.~Gubler, M.~Harris, S.~Kudla, V.~Maillot, D.~Roessler, C.~Soul\'e, 
J.~Wildeshaus. Furthermore, we would like to thank EAGER, the
Arithmetic Geometry Network, the Newton Institute (Cambridge), 
and the Institut Henri Poincar\'e (Paris) for partial support 
of our work.

\section{Log and log-log differential forms} 
\label{sec:loglog}

In this section, we will introduce several complexes of differential
forms with singularities along a normal crossing divisor $\DD$, and we will
discuss their basic properties.

The first one $\mathscr{E}^{\ast }_{X}\langle \DD\rangle$ is a
complex with logarithmic growth conditions in the spirit of
\cite{HarrisPhong:cdcli}. Unlike in \cite{HarrisPhong:cdcli}, where
the authors consider only differential forms of type $(0,q)$, here we
consider the whole Dolbeault complex and we show that it is an acyclic
resolution of the complex of holomorphic forms with logarithmic poles
along the normal crossing divisor $\DD$, i.e.,  this complex
computes the cohomology of the complement of $\DD$. Another difference
with \cite{HarrisPhong:cdcli} is that, in order to be able to prove
the Poincar\'e lemma for such forms, we need to impose growth
conditions to all derivatives of the functions. Note that a similar
condition has been already considered in \cite{HarrisZucker:BcSvIII}.

The second complex $\mathscr{E}^{\ast }_{X}\langle\langle
\DD\rangle\rangle$ contains differential forms with singularities of
log-log type along a normal crossing divisor $\DD$, and is related
with the complex of good forms in the sense of \cite{Mumford:Hptncc}.
As the complex of good forms, it contains the Chern forms of fully
decomposed automorphic hermitian vector bundles and is functorial with
respect to certain inverse images. Moreover all the differential forms
belonging to this complex are locally integrable with zero residue.
The new property of this complex is that it satisfies a Poincar\'e
lemma that implies that this complex is quasi-isomorphic to the
complex of smooth differential forms, i.e., this complex computes the
cohomology of the whole variety. The main interest of this complex, as
we shall see in subsequent sections, is that it contains also the
Bott-Chern forms associated to fully decomposed automorphic vector
bundles.  Note that neither the complex of good forms in the sense of
\cite{Mumford:Hptncc} nor the complex of log-log forms are contained
in each other.

The third complex $\mathscr{E}^{\ast }_{X}\langle \DD_1 \langle \DD_2
\rangle\rangle$ that we will introduce is a mixture of the previous
complexes.  It is formed by differential forms which are log along a
normal crossing divisor $\DD_{1}$ and log-log along another normal
crossing divisor $\DD_{2}$. This complex computes the cohomology of
the complement of $\DD_{1}$.

By technical reasons we will introduce several other complexes.

\subsection{Log forms}

\nnpar{General notations.}
Let $X$ be a complex manifold of dimension $d$. We will
denote by $\mathscr{E}^{\ast }_{X}$ the sheaf of complex smooth
differential forms over $X$.

Let $\DD$ be a normal crossing divisor on $X$. Let $V$ be an open
coordinate subset of $X$ with coordinates $z_{1},\dots ,z_{d}$; we put
$r_{i}=|z_{i}|$. We will say that $V$ is adapted to $\DD$ if
the divisor $\DD\cap V$ is given
by the equation $z_{1}\dots z_{k}=0$ and the coordinate
neighborhood $V$ is small enough; more precisely, we will assume that
all the coordinates satisfy $r_{i}\le 1/e^{e}$, which implies that $\log 
1/r_{i}>e$ and $\log (\log (1/r_{i}))>1$.

We will denote by by $\Delta
_{r}\subset \mathbb{C}$ the open disk of radius $r$ centered at $0$,
by $\overline {\Delta }_{r }$ the closed 
disk, and we will write 
$\Delta^{\ast}_{r}=\Delta
_{r}\setminus \{0\}$ and $\overline {\Delta } ^{\ast}_{r }=\overline
{\Delta }_{r}\setminus \{0\}$.

If $f$ and $g$ are two functions with non-negative real values, we
write $f\prec g$ if there exists a real constant $C>0$ such that
$f(x)\le C \cdot g(x)$ for all $x$ in the domain of definition under
consideration.

\nnpar{multi-indices.}  We collect here all the conventions we will
use about multi-indices.

\begin{notation} \label{def:11}
For any multi-index $\alpha =(\alpha _{1},\dots
,\alpha _{d})\in \mathbb{Z}_{\ge
  0}^{d}$ we write 
\begin{gather*}
  |\alpha |=\sum_{i=1}^{d}\alpha _{i},\qquad 
  z^{\alpha }= \prod_{i=1}^{d}z_{i}^{\alpha _{i}},\qquad 
\bar z^{\alpha }= \prod_{i=1}^{d}
 \bar z_{i}^{\alpha _{i}},\\
 r^{\alpha }=\prod_{i=1}^{d} r_{i}^{\alpha _{i}},\qquad
 (\log(1/r))^{\alpha }= \prod_{i=1}^{d} (\log(1/r_{i}))^{\alpha _{i}},\\
  \frac{\partial^{|\alpha |}}{\partial z^{\alpha }}f=
  \frac{\partial^{|\alpha |}}{\prod_{i=1}^{d}\partial z_{i}^{\alpha
  _{i}}}f,\qquad
  \frac{\partial^{|\alpha |}}{\partial \bar z^{\alpha }}f=
  \frac{\partial^{|\alpha |}}{\prod_{i=1}^{d}\partial \bar 
    z_{i}^{\alpha _{i}}}f.\ 
\end{gather*}

If $\alpha $ and $\beta $ are multi-indices we write
$\beta \ge \alpha $ if, for all $i=1,\dots ,d$, $\beta _{i}\ge \alpha
_{i}$. 
We denote by $\alpha +\beta $ the multi-index with components
$\alpha _{i}+\beta _{i}$. If $1\le i\le d$ we will denote by $\gamma
^{i}$ the multi-index with all the entries zero except the $i$-th
entry that takes the value $1$. More generally if $I$ is a subset of
$\{1,\dots ,d\}$ 
we will denote by $\gamma ^{I}$ the multi-index 
\begin{displaymath}
  \gamma ^{I}_{i}=
  \begin{cases}
    1,\  &i\in I,\\
    0,\ &i\not\in I.
  \end{cases}
\end{displaymath}

We will denote by $\underline n$ the constant multi-index  
\begin{displaymath}
  \underline n_{i}= n.
\end{displaymath}
In particular $\underline 0$ is the
multi-index $\underline 0=(0,\dots ,0)$.

If $\alpha $ is a multi-index and $k\ge 1$ is an integer, we will
denote by $\alpha ^{\le k}$ the multi-index
\begin{displaymath}
  \alpha ^{\le k}_{i}=
    \begin{cases}
    \alpha _{i},\  &i\le k,\\
    0,\ &i> k.
  \end{cases}
\end{displaymath}

For a multi-index $\alpha $, the order function associated to $\alpha
$,
\begin{displaymath}
  \Phi _{\alpha }:\{1,\dots ,|\alpha |\}\longrightarrow \{1,\dots ,d\}
\end{displaymath}
is given by
\begin{displaymath}
     \Phi _{\alpha }(i)=k, \text{ if }
     \sum_{j=1}^{k-1}\alpha _{j} < i \le \sum_{j=1}^{k}\alpha _{j}. 
\end{displaymath}
\end{notation}

\nnpar{Log forms.}
We introduce now a complex of differential forms with
logarithmic growth along a normal crossing divisor. This complex
can be used to compute the cohomology of a non compact algebraic
complex manifold, with its usual Hodge filtration. It
contains the $C^{\infty}$ logarithmic Dolbeault complex defined in
\cite{Burgos:CDc} but it is much bigger and, in particular,
it contains also the log-log differential forms defined later. In contrast
with the pre-log forms introduced in \cite{BurgosKramerKuehn:cacg}, in the
definition given here we impose growth conditions to the differential
forms and to all their derivatives. 

The problem of the weight filtration of the complex of log forms will
not be treated here. 

Let $X$ be a complex manifold of dimension $d$, $\DD$ a normal
crossing divisor, $U=X\setminus \DD$ and $\iota :U\longrightarrow X$ the
inclusion.

\begin{definition}  \label{def:log} Let $V$ be a coordinate
  neighborhood adapted to $\DD$. For every integer $K\ge 0$, we say
  that a smooth complex 
  function $f$ on $V\setminus \DD$  
  has \emph{logarithmic growth along $\DD$ of order $K$} if there
  exists an integer $N_{K}$ such that, 
  for every pair of multi-indices $\alpha, \beta  \in \mathbb{Z}_{\ge
  0}^{d}$, with $|\alpha +\beta |\le K$, it holds the inequality
  \begin{equation}
    \left |\frac{\partial^{|\alpha |}}{\partial z^{\alpha }}
      \frac{\partial^{|\beta |}}{\partial \bar z^{\beta }}
      f(z_{1},\dots ,z_{d})\right | \prec \frac{ \left
        |\prod_{i=1}^{k} \log(1/r_{i}) 
      \right |^{N_{K}}}{|z^{\alpha^{\le k}}\bar z^{\beta ^{\le k}}|}.
  \end{equation}
  We say that $f$ 
  has \emph{logarithmic growth along $\DD$ of infinite order} if it
  has logarithmic growth along $\DD$ of order $K$ for all $K\ge 0$. 
  The \emph{sheaf of differential forms on $X$ with
   logarithmic growth of infinite order along $\DD$}, denoted
   $\mathscr{E}^{\ast}_{X}\langle \DD\rangle $, is the sub-algebra of
   $\iota _{\ast 
   }\mathscr{E}^{\ast }_{U}$ generated, in each coordinate
 neighborhood adapted to $\DD$, by the functions with
 logarithmic growth of infinite order along $\DD$ and the differentials
 \begin{equation}\label{eq:23}
   \begin{alignedat}{2}
     & \frac{\dd z_{i}}{z_{i}},\  \frac{\dd\bar {z}_{i}}{\bar
       {z}_{i}},&
     \qquad \text{for } i&=1,\dots ,k,\\
     & \dd z_{i},\  \dd\bar {z}_{i},&
     \qquad \text{for } i&=k+1,\dots ,d.
   \end{alignedat}
 \end{equation}
 For shorthand, a differential form with logarithmic growth of infinite
 order along $\DD$ is called log along $\DD$ or, if 
 $\DD$ is understood, a log form.
\end{definition}

\nnpar{The Dolbeault algebra of log forms.}
The sheaf $\mathscr{E}^{\ast }_{X}\langle \DD\rangle $ inherits from
$\iota _{\ast}\mathscr{E}^{\ast}_{U}$ a real
structure and a bigrading. Moreover, it is clear that, if $\omega $ is
a log form then $\partial \omega $ and $\bar
\partial \omega $ are also log forms. Therefore $\mathscr{E}^{\ast
}_{X}\langle \DD\rangle $ is a sheaf of Dolbeault algebras.
We will use all the notations of
\cite{BurgosKramerKuehn:cacg} \S5 concerning Dolbeault algebras. For
the convenience of the reader we will recall these notations in
section \ref{sec:dolb-algebr-deligne}. In
particular, from the structure of Dolbeault algebra, there is a well
defined Hodge filtration denoted by $F$.

\nnpar{Pre-log forms.} Recall that, in \cite{BurgosKramerKuehn:cacg}
section 7.2 there is introduced the sheaf of pre-log forms denoted  
$\mathscr{E}^{\ast }_{X}\langle \DD\rangle_{\wlg} $. It is clear that
there is an inclusion of sheaves
\begin{displaymath}
 \mathscr{E}^{\ast }_{X}\langle \DD\rangle \subset 
\mathscr{E}^{\ast }_{X}\langle \DD\rangle _{\wlg}.
\end{displaymath}

\nnpar{The cohomology of the complex of log forms.} 
Let $\Omega^{\ast}_{X}(\log \DD)$ be the sheaf of holomorphic forms
with logarithmic poles along $\DD$ \cite{Deligne:THII}. 
Then the more general theorem  \ref{thm:fq} implies: 
\begin{theorem} \label{thm:logfq} The inclusion 
\begin{displaymath}
  \begin{CD}
    \Omega ^{\ast }_{X}(\log \DD) @>>>\mathscr{E}^{\ast }_{X}\langle
    \DD\rangle
  \end{CD}
\end{displaymath}
is a filtered quasi-isomorphism with respect to the Hodge filtration.
\end{theorem}
In other words, this complex is a resolution of the sheaf of
holomorphic forms with logarithmic poles along $\DD$,
$\Omega ^{\ast }_{X}(\log \DD)$. Thus if $X$ is a compact K\"ahler
manifold, the complex
of global sections $\Gamma (X,\mathscr{E}^{\ast }_{X}\langle \DD\rangle
)$ computes the cohomology of the open complex manifold $U=X\setminus
\DD$ with its Hodge 
filtration. 

Note that corollary \ref{thm:logfq} implies that there is an
isomorphism in the derived category $R\iota _{\ast} \underline
{\mathbb{C}}_{U}\longrightarrow \mathscr{E}^{\ast }_{X}\langle
\DD\rangle$. This isomorphism is compatible with the real structures.
Hence the complex $\mathscr{E}^{\ast }_{X}\langle \DD\rangle$ also
provides the real structure of the cohomology of $U$.

\nnpar{Inverse images.} 
The complex of log forms is
functorial with respect to inverse images.
More precisely, we have the following result.

\begin{proposition} \label{prop:4}
Let $f:X\longrightarrow Y$ be a morphism of complex
manifolds of dimension $d$ and $d'$. Let $\DD_{X}$, $\DD_{Y}$ be normal
crossing divisors on $X$ and
$Y$ respectively, satisfying $f^{-1}(\DD_{Y})\subseteq \DD_{X}$.  If
$\eta$ is a section of $\mathscr{E}^{\ast}_{Y}\langle
\DD_{Y}\rangle$, then $f^{\ast}\eta$ is a section of
$\mathscr{E}^{\ast}_{X}\langle \DD_{X}\rangle$.
\end{proposition}
\begin{proof} Let $p$ be a point of $X$.
  Let $V$ and $W$ be open coordinate neighborhood of $p$ and $f(p)$
  respectively, adapted to $\DD_{X}$
  and $\DD_{Y}$, and such that $f(V)\subset W$. Let $k$ and
  $k'$ be the number of components of $V\cap \DD_{X}$ and $W\cap
  \DD_{Y}$. Then the 
  condition $f^{-1}(\DD_{Y})\subseteq \DD_{X}$ implies that $f$ can be
  written as
  \begin{equation}\label{eq:11}
    f(x_{1},\dots ,x_{d})=(z_{1},\dots ,z_{d'}),
  \end{equation}
  with
  \begin{displaymath}
    z_{i}=
    \begin{cases}
      x_{1}^{a_{i,1}}\cdots x_{k}^{a_{i,k}}u_{i}, &\text{ if }i\le
      k',\\
      w_{i},& \text{ if }i>k', 
    \end{cases}
  \end{displaymath}
  where $u_{1}, \dots , u_{k'}$ are holomorphic functions that do not
  vanish in $V$, the $a_{i,j}$ are non negative integers and
  $w_{k'+1},\dots ,w_{d'}$ are holomorphic 
  functions. Shrinking $V$
  if necessary, we may assume that the 
  functions $u_{j}$ are holomorphic and do not vanish in a
  neighborhood of the adherence of $V$.

  For $1\le i \le k'$ we have
  \begin{displaymath}
    f^{\ast}\left(\frac{\dd z_{i}}{z_{i}}\right)=\sum_{j=1}^{k}
    a_{i,j} \frac{ \dd x_{j}}{x_{j}} + \frac{ \dd u_{i}}{u_{i}}.
  \end{displaymath}
  since the function $1/u_{i}$ is holomorphic in a neighborhood of the
  adherence of $V$, the function $1/u_{i}$ and all its derivatives are
  bounded. If follows
  that $f^{\ast}(\dd z_{i}/z_{i})$ is a log form  (along $\DD_{X}$). The
  same argument 
  shows that $f^{\ast}(\dd \bar z_{i}/\bar z_{i})$ is a log form.

  If a function $g$ on $W$ satisfies 
  \begin{displaymath}
    |g(z_{1},\dots ,z_{d'})|\prec \left |
      \prod_{i=1}^{k'} \log (1/|z_{i}|)\right | ^{N},
  \end{displaymath}
  then $f^{\ast}g$ satisfies
  \begin{align*}
    |f^{\ast} g(x_{1},\dots ,x_{d})|&\prec \left |
      \prod_{i=1}^{k'}\left( \sum_{j=1}^{k} a_{i,j}\log (1/|x_{j}|)+
      \log (1/|u_{i}|)\right)\right | ^{N}\\
    &\prec \left |
      \prod_{j=1}^{k} \log (1/|x_{j}|)\right | ^{Nk'}.
  \end{align*}
  therefore $f^{\ast}g$ has logarithmic growth. It remains to bound the
  derivatives of $f^{\ast}g$. For ease of notation we will bound only
  the derivatives with respect to the holomorphic coordinates, being the
  general case analogous.

  For any multi-index $\alpha \in \mathbb{Z}_{\ge 0}^{d}$, the
  function $\partial^{|\alpha |}/\partial x^{\alpha }(f^{\ast}g)$ is a
  linear combination of the functions
  \begin{equation} \label{eq:12}
    \left\{\frac{\partial^{|\beta | }}{\partial z^{\beta }}g
    \prod_{i=1}^{|\beta |}\frac{\partial^{|\alpha ^{i}|}}
    {\partial x^{\alpha ^{i}}}z_{\Phi
        _{\beta }(i)}\right\}_{\beta ,\{\alpha _{i}\}},  
  \end{equation}
  where $\beta $ runs over all multi-indices $\beta \in
  \mathbb{Z}_{\ge 0}^{d'}$ such that
  $|\beta |\le |\alpha |$, and $\{\alpha _{i}\}$ runs over all
  families of multi-indices
  $\alpha ^{i}\in \mathbb{Z}_{\ge 0}^{d}$ such that 
  \begin{displaymath}
    \sum_{i=1}^{|\beta |}\alpha ^{i}=\alpha. 
  \end{displaymath}
  The function $\Phi _{\alpha }$ is the order function introduced in
  \ref{def:11}. 

  Then, since $g$ is a log function,
  \begin{align*}
    \left|\frac{\partial^{|\beta |} }{\partial z^{\beta }}g
    \prod_{i=1}^{|\beta |}\frac{\partial^{|\alpha ^{i}|}}
    {\partial x^{\alpha ^{i}}}z_{\Phi
        _{\beta }(i)} \right |&\prec
    \frac{\left|\prod_{j=1}^{k'} \log (1/|z _{j})|\right|^{N_{|\beta |}}}
    {|z^{\beta ^{\le k'}}|}\left|
    \prod_{i=1}^{|\beta |}\frac{\partial^{|\alpha ^{i}|}}
    {\partial x^{\alpha ^{i}}} z_{\Phi
        _{\beta }(i)}\right|\\
    &\prec 
    \left|
      \prod_{j=1}^{k} \log (1/|x_{j}|)\right | ^{N_{|\beta |}k'}
    \prod_{i=1}^{|\beta ^{\le k'}|}    \left | 
    \frac{1}{z_{\Phi _{\beta }(i)}}
    \frac{\partial^{|\alpha ^{i}|}}{\partial x^{\alpha ^{i}}}z_{\Phi
        _{\beta }(i)}\right|.
  \end{align*}
  But, by the assumption on the map $f$, it is easy to see that, for $1\le
  j\le k'$, we have
  \begin{equation} \label{eq:2}
    \left| \frac{1}{z_{j}}\frac{\partial ^{|\alpha
          ^{i}|}}{\partial x^{\alpha ^{i}}}z_{j}\right |\prec
    \frac{1}{|x^{(\alpha ^{i})^{\le k}}|}, 
  \end{equation}
  which implies the lemma.
\end{proof}

\nnpar{Polynomial growth in the local universal cover.} We can
characterize log forms as 
differential forms that have polynomial growth in a local universal
cover.  Let $M>1$ be a real number and let
$U_{M}\subset \mathbb{C}$ be the subset given by 
\begin{displaymath}
  U_{M}=\{x\in \mathbb{C}\mid \Im x >
  M\}.
\end{displaymath}
Let $K$ be an open subset of $\mathbb{C}^{d-k}$. We consider the
space $(U_{M})^{k}\times K$ with coordinates $(x_{1},\dots
,x_{d})$.

\begin{definition} \label{def:23}
  A function $f$ on $(U_{M})^{k}\times K$ is said to have
  \emph{imaginary polynomial growth } 
  if there is a sequence of integers $\{N_{n}\}_{n\ge 0}$ such that
  for every pair of multi-indices $\alpha, \beta  \in \mathbb{Z}_{\ge
  0}^{d}$,  the inequality
  \begin{equation}
    \left |\frac{\partial^{|\alpha |}}{\partial x^{\alpha }}
      \frac{\partial^{|\beta |}}{\partial \bar x^{\beta }}
      f(x_{1},\dots ,x_{d})\right | \prec \left
        |\prod_{i=1}^{k} \Im x_{i} 
      \right |^{N_{|\alpha |+|\beta |}}
  \end{equation}
  holds.
  The space of differential forms on $(U_{M})^{k}\times K$ with
  \emph{imaginary polynomial 
  growth} is generated by the functions with imaginary polynomial
  growth and the differentials 
 \begin{alignat*}{2}
   & \dd x_{i},\  \dd\bar {x}_{i},&
   \qquad \text{for } i&=1,\dots ,d.
 \end{alignat*}
\end{definition}

 Let $X$,
$\DD$, $U$ and $\iota $ be as 
in definition \ref{def:log}. 

\begin{definition} \label{def:24} Let $W$ be an open subset of $X$ and
  $\omega $ be a differential form in $\Gamma (W,\iota
  _{\ast}(\mathscr{E}_{U}^{\ast})) $. For every point $p\in W$, there is 
  an open coordinate neighborhood $V\subseteq W$, which is adapted to
  $\DD$ and such that the coordinates of 
  $V$ induce an identification $V\cap U= (\Delta ^{\ast}_{r})^{k}\times K$.
  We choose
  $M>\log (1/r)$ and denote by $\pi:(U_{M})^{k}\times 
  K\longrightarrow V$ the covering  map given by
  \begin{displaymath}
    \pi (x_{1},\dots ,x_{d})=
    (e^{2\pi i x_{1}},\dots ,e^{2\pi i x_{k}},x_{k+1},\dots ,x_{d}).
  \end{displaymath}
  We say that $\omega $ has
  \emph{polynomial growth in 
  the local universal cover}, if, for every $V$ and $\pi$ as above,
  $\pi ^{\ast}\omega $ has imaginary
  polynomial growth.
\end{definition}

It is easy to see that the differential forms with polynomial growth in
 the local universal cover form a sheaf of Dolbeault algebras.

\begin{theorem}\label{thm:15}
  A differential form has polynomial growth in
  the local universal cover if and only if it is a log form.
\end{theorem}
\begin{proof}
  We start with the case of a function. So let $f$ be a function with
  polynomial growth in the local universal cover and let $V$ be a
  coordinate neighborhood as in definition \ref{def:24}. Let $g=\pi
  ^{\ast}f$. By definition, $g$ satisfies
  \begin{equation}\label{eq:7}
    g(\dots ,x_{i}+1,\dots)=g(\dots ,x_{i},\dots ),\qquad\text{for
    }1\le i\le k.  
  \end{equation}
We write formally   
  \begin{displaymath}
    f(z_{1},\dots ,z_{d})=g(x_{1}(z_{1}),\dots ,x_{d}(z_{d})),
  \end{displaymath}
  with
  \begin{displaymath}
    x_{i}(z_{i})=
    \begin{cases}
      \frac{1}{2\pi i} \log z_{i},\ &\text{for }i\le k,\\
      z_{i},\ &\text{for }i>k. 
    \end{cases}
  \end{displaymath}
  Note that this makes sense because of the periodicity properties
  \eqref{eq:7}. Then we have
  \begin{multline}\label{eq:21}
    \frac{\partial^{|\alpha |}}{\partial z^{\alpha }}
      \frac{\partial^{|\beta |}}{\partial \bar z^{\beta }}
      f(z_{1},\dots ,z_{d})=
      \sum_{\substack{\alpha '\le\alpha \\ \beta ' \le\beta }}
      C_{\alpha ,\beta}^{\alpha ',\beta '}
      \frac{\partial^{|\alpha' |}}{\partial x^{\alpha' }}
      \frac{\partial^{|\beta' |}}{\partial \bar x^{\beta '}}
      g(x_{1},\dots ,x_{d}) \cdot\\
      \cdot \frac{\partial^{|\alpha -\alpha '|}}{\partial z^{\alpha-\alpha ' }}
      \left(
        \frac{\partial x}{\partial z}\right)^{\alpha '}
      \frac{\partial^{|\beta -\beta '|}}{\partial \bar z^{\beta -\beta
          '}}\left(
        \frac{\partial \bar x}{\partial \bar z}  
      \right)^{\beta'},
  \end{multline}
  for certain constants $C_{\alpha ,\beta}^{\alpha ',\beta '}$. But
  the estimates
  \begin{align*}
      \left |\frac{\partial^{|\alpha' |}}{\partial x^{\alpha' }}
        \frac{\partial^{|\beta' |}}{\partial \bar x^{\beta' }}
        g(x_{1},\dots ,x_{d})\right | &\prec 
      \left
          |\prod_{i=1}^{k} |x_{i}| 
        \right |^{N_{\alpha ',\beta '}} 
\prec
      \left
          |\prod_{i=1}^{k} \log (1/|z_i|)) 
        \right |^{N_{\alpha ',\beta '}}  
  \end{align*}
  and
  \begin{equation}\label{eq:22}
    \frac{\partial^{|\alpha -\alpha '|}}{\partial z^{\alpha-\alpha ' }}
      \left(
        \frac{\partial x}{\partial z}\right)^{\alpha '}
      \frac{\partial^{|\beta -\beta '|}}{\partial \bar z^{\beta -\beta
          '}}\left(
        \frac{\partial \bar x}{\partial \bar z}  
      \right)^{\beta'}\prec \frac{1}{|z^{\alpha^{\le k}}\bar z^{\beta
      ^{\le k}}|} 
  \end{equation}
  imply the bounds of $f$ and its derivatives. The converse is
  proved in the same way. 

 To prove the theorem for differential forms, observe that, 
  for $1\le i \le k$,  
 \begin{displaymath}
   \pi ^{\ast} \left( \frac{\dd z_{i}}{z_{i}}\right)=2\pi i \dd
   x_{i}.  
 \end{displaymath}
\end{proof}

\subsection{Log-log forms}

\nnpar{log-log growth forms.} Let $X$, $\DD$, $U$ and $\iota $ be as
in definition \ref{def:log}.

\begin{definition}
  \label{def:logloggrowth} Let $V$ be a coordinate neighborhood adapted
  to $\DD$. For every integer $K\ge 0$, we say
  that a smooth complex 
  function $f$ on $V\setminus \DD$  
  has \emph{log-log growth along $\DD$ of order $K$} if there
  exists an integer $N_{K}$ such that, 
  for every pair of multi-indices $\alpha, \beta  \in \mathbb{Z}_{\ge
  0}^{d}$, with $|\alpha +\beta |\le K$, it holds the inequality
  \begin{equation}
    \left |\frac{\partial^{|\alpha |}}{\partial z^{\alpha }}
      \frac{\partial^{|\beta |}}{\partial \bar z^{\beta }}
      f(z_{1},\dots ,z_{d})\right | \prec \frac{ \left
        |\prod_{i=1}^{k} \log (\log(1/r_{i})) 
      \right |^{N_{K}}}{|z^{\alpha^{\le k}}\bar z^{\beta ^{\le k}}|}.
  \end{equation}
  We say that $f$ 
  has \emph{log-log growth along $\DD$ of infinite order}, if  it
  has log-log growth along $\DD$ of order $K$ for all $K\ge 0$.
  The \emph{sheaf of differential forms on $X$ with
    log-log growth along $\DD$ of infinite order} is the subalgebra of
  $\iota _{\ast}\mathscr 
  {E}^{\ast}_{U}$ generated, in each coordinate neighborhood $V$
  adapted to $\DD$, by the functions with log-log growth along $\DD$
  and the differentials
  \begin{alignat*}{2}
    &\frac{\dd z_{i}}{z_{i}\log(1/r_{i})},\,\frac{\dd\bar{z}_{i}}
    {\bar{z}_{i}\log(1/r_{i})},&\qquad\text{for }i&=1,\dots,k, \\
    &\dd z_{i},\,\dd\bar{z}_{i},&\qquad\text{for }i&=k+1,\dots,d.
  \end{alignat*}
  A differential form with log-log growth along $\DD$ of infinite
  order will be called 
  a \emph{log-log growth form}. The sheaf of differential forms on $X$
  with log-log growth along $\DD$ of infinite order will be denoted 
  $\mathscr{E}^{\ast}_{X}
  \langle\langle \DD\rangle\rangle_{\growth}$.
\end{definition}

The following characterization of differential forms with log-log
growth of infinite order is left to the reader.
\begin{lemma} \label{lemm:1}
  Let $V$ be an open coordinate subset adapted to $\DD$ and let $I,J$ be
  two subsets of $\{1,\dots ,d\}$. Then the form
  $f\dd z_{I}\land \dd \bar z_{J}$ is a log-log
  growth form of infinite order if and only if,
  for every pair of multi-indices $\alpha, \beta  \in \mathbb{Z}_{\ge
  0}^{d}$, there is an integer 
  $N_{\alpha ,\beta }\ge 0$  such that 
  \begin{equation}\label{eq:14}
    \left |\frac{\partial^{|\alpha |}}{\partial z^{\alpha }}
      \frac{\partial^{|\beta |}}{\partial \bar z^{\beta }}
      f(z_{1},\dots ,z_{d})\right | \prec \frac{ \left
        |\prod_{i=1}^{k} \log (\log(1/r_{i})) 
      \right |^{N_{\alpha ,\beta }}}{r^{(\gamma ^{I}+\gamma
      ^{J}+\alpha+\beta )^{\le k}}(\log (1/r))^{(\gamma ^{I}+\gamma
      ^{J})^{\le k}}}. 
  \end{equation}
\hfill $\square$
\end{lemma}

\begin{definition} \label{def:16} A function that satisfies the bound
  \eqref{eq:14} 
  for any pair of multi-indices $\alpha , \beta $ with $\alpha +\beta \le K$
  will be called a \emph{$(I,J)$-log-log growth function of order
  $K$}. If it satisfies the bound \eqref{eq:14}
  for any pair multi-indices $\alpha , \beta $ it will be called a
  \emph{$(I,J)$-log-log growth function of infinite order}.
\end{definition}

\nnpar{Log-log forms.} Unlike the case of log growth forms, the fact
that $\omega $ is a log-log growth form does not imply that its
differential  $\partial
\omega $ is a log-log growth form.
\begin{definition}  \label{def:loglog} We say that a smooth complex
  differential form $\omega $  is \emph{log-log
  along $\DD$} if the differential forms $\omega $, $\partial \omega $,
  $\bar \partial \omega $ and $\partial \bar \partial \omega $  have
  log-log growth along $\DD$ of infinite 
  order. The sheaf of differential forms log-log
  along $\DD$ will be denoted by $\mathscr{E}^{\ast}_{X}
  \langle\langle \DD\rangle\rangle$. 
  For shorthand, if $\DD$ is clear from the context, a differential
  form which is log-log along $\DD$ will be called a log-log form.
\end{definition}

From the definition, it is clear that the sheaf of log-log forms is
contained in the sheaf of log forms.

Let $V$ be a coordinate subset adapted to $\DD$. For $i=1,\dots ,k$, the
 function $\log (\log (1/r_{i}))$ is a log-log function and the 
 differential forms
\begin{displaymath}
   \frac{\dd z_{i}}{z_{i}\log (1/r_{i})},\  
  \frac{\dd\bar {z}_{i}}{\bar
    {z}_{i}\log (1/r_{i})},
  \qquad \text{for } i=1,\dots ,k,
\end{displaymath}
are log-log forms.

\nnpar{The Dolbeault algebra of log-log forms.} 
As in the case of log forms, the sheaf $\mathscr{E}^{\ast }_{X}\langle
\langle \DD\rangle \rangle$ inherits  a real structure and a
bigrading. Moreover, we have forced the existence of operators
$\partial$ and $\bar \partial$. Therefore  $\mathscr{E}^{\ast }_{X}\langle
\langle \DD\rangle \rangle$ is a 
sheaf of Dolbeault algebras (see section
\ref{sec:dolb-algebr-deligne}). 
In
particular, there is a well
defined Hodge filtration, denoted by $F$.

\nnpar{Pre-log-log forms} Recall that, in \cite{BurgosKramerKuehn:cacg}
section 7.1 there is introduced the sheaf of pre-log-log forms,
denoted by  
$\mathscr{E}^{\ast }_{X}\langle \langle \DD\rangle \rangle_{\wlg} $. It is clear that
there is an inclusion of sheaves
\begin{displaymath}
 \mathscr{E}^{\ast }_{X}\langle\langle \DD\rangle \rangle \subset 
\mathscr{E}^{\ast }_{X}\langle\langle \DD\rangle \rangle _{\wlg}.
\end{displaymath}

\nnpar{The cohomology of the complex of log-log differential forms.}
Let
$\Omega^{\ast}_{X}$ be the sheaf of holomorphic forms.  Then theorem
\ref{thm:fq}, which will be proved later, implies:

\begin{theorem} \label{thm:loglogfq}  The inclusion 
\begin{displaymath}
  \begin{CD}
    \Omega ^{\ast }_{X} @>>>\mathscr{E}^{\ast }_{X}\langle \langle \DD
    \rangle \rangle  
  \end{CD}
\end{displaymath}
is a filtered quasi-isomorphism with respect to the Hodge filtration.
\end{theorem}

In other words, this complex is a resolution of $\Omega ^{\ast }_{X}$,
the sheaf of holomorphic differential forms on $X$. Therefore if $X$
is a compact K\"ahler manifold, the complex of global sections $\Gamma
(X,\mathscr{E}^{\ast }_{X}\langle\langle \DD\rangle \rangle )$ computes
the cohomology of $X$ with its Hodge
filtration. As in the case of log forms it also provides the usual
real structure of the cohomology of $X$. One may say that the
singularities of the log-log complex  are so 
mild that they don't change the cohomology.

\nnpar{Inverse images.} 
As in the case of pre-log-log forms, the sheaf of log-log forms is
functorial with respect to inverse images.
More precisely, we have the following result.

\begin{proposition} \label{prop:10}
Let $f:X\longrightarrow Y$ be a morphism of complex
manifolds of dimension $d$ and $d'$, let $\DD_{X}$, $\DD_{Y}$ be normal
crossing divisors on $X$ and 
$Y$ respectively, satisfying $f^{-1}(\DD_{Y})\subseteq \DD_{X}$.  If
$\eta$ is a section of $\mathscr{E}^{\ast}_{Y}\langle\langle
\DD_{Y}\rangle\rangle$, then $f^{\ast}\eta$ is a section of
$\mathscr{E}^{\ast}_{X}\langle\langle \DD_{X}\rangle\rangle$.
\end{proposition}
\begin{proof}
  Since the differential operators $\partial$ and $\bar
  \partial$ are compatible with inverse images, we have to show that
  the pre-image of a form with log-log 
  growth of infinite order has log-log growth of infinite order. We
  may assume that, locally, $f$ can be written as in equation \eqref{eq:11}. 
  If a function $g$ satisfies
  \begin{displaymath}
    |g(z_{1},\dots ,z_{d'})|\prec
    \left | \prod_{i=1}^{k'} \log (\log (1/|z_{i}|))\right |^{N},
  \end{displaymath}
  Then 
  \begin{align*}
    |(f^{\ast}g)(x_{1},\dots ,x_{d})|&\prec
    \left | \prod_{i=1}^{k'} f^{\ast}\log (\log (1/|z_{i}|))\right |^{N}\\
    &\prec
    \left | \prod_{i=1}^{k'} \sum_{j=1}^{k} \log (\log
    (1/|x_{j}|))\right |^{N}\\ 
    &\prec 
    \left | \sum_{j=1}^{k} \log (\log (1/|x_{j}|))\right |^{Nk'}.
  \end{align*}
  Therefore $f^{\ast}g $ has log-log growth.
  
  Next we have to bound the derivatives of $f^{\ast}g $. As in the
  proof of proposition \ref{prop:4} we will bound only the derivatives
  with respect to the holomorphic coordinates. Again we observe that,
  for any multi-index $\alpha \in \mathbb{Z}_{\ge 0}^{d}$, the
  function $\partial^{|\alpha |}/\partial x^{\alpha }(f^{\ast}g)$ is a
  linear combination of the functions \eqref{eq:12}.
  But, using that $g$ is a log-log growth function, we can bound
  \begin{align*}
    \left|\frac{\partial^{|\beta |} }{\partial z^{\beta }}g
      \prod_{i=1}^{|\beta |}\frac{\partial^{|\alpha ^{i}|}}{\partial
        x^{\alpha ^{i}}} z_{\Phi _{\beta }(i)}\right |&\prec
    \frac{\left|\prod_{j=1}^{k'} \log (\log (1/|z _{j}|))\right|^{N_{|\beta |}}}
    {|z^{\beta ^{\le k'}}|}\left| \prod_{i=1}^{|\beta
        |}\frac{\partial^{|\alpha ^{i}|}}{\partial x^{\alpha ^{i}}}
      z_{\Phi
        _{\beta }(i)}\right|\\
    &\prec \left| \prod_{j=1}^{k} \log (\log (1/|x_{j}|))\right |
    ^{N_{|\beta |}k'} \prod_{i=1}^{|\beta ^{\le k'}|} \left | \frac{1}{z_{\Phi
          _{\beta }(i)}} \frac{\partial^{|\alpha ^{i}|}}{\partial
        x^{\alpha ^{i}}}z_{\Phi _{\beta }(i)}\right|.
  \end{align*}
  By \eqref{eq:2}
  \begin{displaymath}
    \left|\frac{\partial^{|\beta |} }{\partial z^{\beta }}g
      \prod_{i=1}^{|\beta |}\frac{\partial^{|\alpha ^{i}|}}{\partial
        x^{\alpha ^{i}}} z_{\Phi _{\beta }(i)}\right |\prec
    \left|
      \prod_{j=1}^{k} \log (\log (1/|x_{j}|))\right | ^{N_{|\beta |}k'}
    \prod_{i=1}^{|\beta ^{\le k'}|}  \frac{1}{|x^{\alpha
        ^{\le k}}|}. 
  \end{displaymath}
  Thus $f^{\ast}g$ has log-log growth of infinite order.

  Finally we are led to study the inverse image  of the differential
  forms
  \begin{displaymath}
    \frac{\dd z_{i}}{z_{i} \log (1/|z_{i}|)}, \  
    \frac{\dd \bar z_{i}}{\bar z_{i} \log (1/|z_{i}|)}, \
    \text{for }i=1,\dots ,k'. 
  \end{displaymath}
  We have
  \begin{displaymath}
    f^{\ast}\left (
      \frac{\dd z_{i}}{z_{i}\log (1/z_{i}\bar z_{i})} 
      \right )
      =\frac{1}{\log (1/z_{i}\bar z_{i})} \left(
      \sum_{i=j}^{k} a_{i,j}\frac{\dd x_{j}}{x_{j}} + 
      \frac{\dd u_{i}}{u_{i}} \right ).  
  \end{displaymath}
  Since we have assumed that $u_{i}$ is a non-vanishing holomorphic
  function in a neighborhood of the adherence of $V$ (see the proof of
  proposition \ref{prop:4}), the function
  $1/u_{i}$ and all its derivatives are bounded. Therefore it only
  remains to show that the functions
  \begin{equation}\label{eq:13}
    f^{\ast}\left(\frac{1}{\log
  (1/|z_{i}|)}\right)\ \text{and}\  \log (1/|x_{j}|) f^{\ast}\left(
  \frac{1}{\log (1/|z_{i}|)}\right), \ 
  \text{for}\ a_{i,j}\not =0,  
  \end{equation}
  have log-log growth of infinite order, which is left to the reader.
\end{proof}

\nnpar{Integrability.} Since the sheaf of log-log forms is contained
in the sheaf of pre-log-log forms, then \cite{BurgosKramerKuehn:cacg}
proposition 7.6 implies
\begin{proposition}\label{prop:5}
\begin{enumerate}[(i)]
\item[(i)]
Any log-log form is locally integrable.
\item[(ii)]
If $\eta$ is a log-log form, and $[\eta ]_{X}$ is the
associated current, then
\begin{displaymath}
[\dd\eta]_{X}=\dd[\eta]_{X}.
\end{displaymath}
The same holds true for the differential operators $\partial$,
$\bar{\partial}$ and $\partial\bar{\partial}$.
\end{enumerate}
\hfill $\square$
\end{proposition}

\nnpar{Logarithmic growth in the local universal cover.} We will
define a new class of 
singular forms closely related to the log-log forms. The discussion
will be parallel to the one at the end of the 
previous section.

 Let $U_{M}$, $K$ and $(x_{1},\dots,x_{d})$ be as in
definition \ref{def:23}. 

\begin{definition}\label{def:25}
  A function $f$ on $(U_{M})^{k}\times K$ is said to have
  \emph{imaginary logarithmic growth }, 
  if there is a sequence of integers $\{N_{n}\}_{n\ge 0}$ such that
  for every pair of multi-indices $\alpha, \beta  \in \mathbb{Z}_{\ge
  0}^{d}$, the inequality
  \begin{equation}
    \left |\frac{\partial^{|\alpha |}}{\partial x^{\alpha }}
      \frac{\partial^{|\beta |}}{\partial \bar x^{\beta }}
      f(x_{1},\dots ,x_{d})\right | \prec \frac{ \left
        |\prod_{i=1}^{k} \log(\Im x_{i}) 
      \right |^{N_{|\alpha |+|\beta |}}}
    {|x^{\alpha^{\le k}}\bar x^{\beta ^{\le k}}|}
  \end{equation}
  holds.
  The space of differential forms with \emph{imaginary logarithmic
  growth} is generated by the functions with imaginary logarithmic
  growth and the differentials 
 \begin{alignat*}{2}
   & \frac{\dd x_{i}}{\Im x_{i}},\  \frac{\dd\bar {x}_{i}}{\Im
     {x}_{i}},&
   \qquad \text{for } i&=1,\dots ,k,\\
   & \dd x_{i},\  \dd\bar {x}_{i},&
   \qquad \text{for } i&=k+1,\dots ,d.
 \end{alignat*}
\end{definition}

Let $X$,
$\DD$, $U$ and $\iota $ be as 
in definition \ref{def:log}. 

\begin{definition} \label{def:22}Let $W$ be an open subset of $X$ and
  let $\omega $ be a differential form in $\Gamma (W,\iota
  _{\ast}(\mathscr{E}_{U})^{\ast}) $.  For every point $p\in W$, there is 
  an open coordinate neighborhood $V\subseteq W$, which is adapted to
  $\DD$ and such that the coordinates of $V$ induce an identification
  $V\cap U= (\Delta ^{\ast}_{r})^{k}\times K$. We choose
  $M>\log (1/r)$ and denote by $\pi:(U_{M})^{k}\times 
  K\longrightarrow V$ the covering  map given by  
  \begin{displaymath}
    \pi (x_{1},\dots ,x_{d})=
    (e^{2\pi i x_{1}},\dots ,e^{2\pi i x_{k}},x_{k+1},\dots ,x_{d}),
  \end{displaymath}
  We say that $\omega $ has
  \emph{logarithmic growth in 
  the local universal cover}, if, for every $V$ and $\pi$ as above,
  $\pi ^{\ast}\omega $ has imaginary
  logarithmic growth.
\end{definition}

It is easy to see that the differential forms with logarithmic growth in
 the local universal cover form a sheaf of Dolbeault algebras.

\begin{theorem}\label{thm:14}
  The sheaf of differential forms with logarithmic growth in
  the local universal cover is contained in the sheaf of log-log forms.
\end{theorem}
\begin{proof}
  Since the forms with logarithmic growth in the local universal cover
  form a Dolbeault algebra, it is enough to check that a differential
  form with logarithmic growth in
  the local universal cover has log-log growth of infinite order. We
  start with the case of a function. So let $f$ and $g$ be as in the
  proof of theorem \ref{thm:15}. To bound the derivatives of $f$ we
  use equation \eqref{eq:21}. But in this case
  \begin{equation}\label{eq:24}
    \begin{aligned}
      \left |\frac{\partial^{|\alpha' |}}{\partial x^{\alpha' }}
        \frac{\partial^{|\beta' |}}{\partial \bar x^{\beta' }}
        g(x_{1},\dots ,x_{d})\right | &\prec 
      \frac{ \left
          |\prod_{i=1}^{k} \log(|x_{i}|) 
        \right |^{N_{\alpha ',\beta '}}}{|x^{\alpha'{}^{\le k}}\bar
        x^{\beta'{}^{\le k}}|}\\ 
      &\prec
      \frac{ \left
          |\prod_{i=1}^{k} \log(\log (1/|z_i|)) 
        \right |^{N_{\alpha ',\beta '}}}{|\log(1/|z|)
        ^{\alpha'{}^{\le k}+\beta'{}^{\le k}}|}.  
    \end{aligned}
  \end{equation}
  Note that now the different terms of equation \eqref{eq:21} have
  slightly different bounds. If we combine the worst bounds of
  \eqref{eq:24} with \eqref{eq:22} then we obtain 
  \begin{equation} \label{eq:25}
        \left |\frac{\partial^{|\alpha |}}{\partial z^{\alpha }}
      \frac{\partial^{|\beta |}}{\partial \bar z^{\beta }}
      f(z_{1},\dots ,z_{d})\right | \prec \frac{ \left
        |\prod_{i=1}^{k} \log (\log(1/r_{i})) 
      \right |^{N_{K}}}{|z^{\alpha^{\le k}}\bar z^{\beta ^{\le k}}|
    \prod_{i=1}^{k} |\log(1/r_{1})|^{\min(\alpha _{i},1)+\min(\beta
      _{i},1)} } .
  \end{equation}
  which imply the bounds of $f$ and its derivatives. 

  To prove the statement for differential forms,
 we observe that for $1\le i \le k$,  
 \begin{displaymath}
   \pi ^{\ast} \frac{\dd z_{i}}{z_{i}\log (1/|z_{i}|)}=\frac{i \dd
   x_{i}}{\Im x_{i}}.  
 \end{displaymath}
\end{proof}

\begin{remark} \label{rem:4}
  The differential forms that interest us are the forms with
  logarithmic growth in the local universal cover. We have introduced
  the log-log forms because it is easier to work with bounds of the
  function and its derivatives in usual coordinates than with the
  condition of logarithmic growth in the local universal cover. This
  is particularly true in the proof of the Poincar\'e lemma. Note
  however that theorem \ref{thm:14} provides us only with an inclusion
  of sheaves and does not give us a characterization of differential
  forms with logarithmic growth in the local universal cover. This can
  be seen by the fact that the bounds  \eqref{eq:25} are sharper than
  the bounds of definition \ref{def:16}. We have choosed the bounds of
  definition \ref{def:16} because the sharper bounds \eqref{eq:25} are
  not functorial. Moreover, they do not characterize forms with logarithmic
  growth in the 
  local universal cover. In fact it does not exist
  a characterization of forms with logarithmic growth in the
  local universal cover in terms of bounds of the function and its
  derivatives in usual coordinates. 
\end{remark}

\subsection{Log and log-log mixed forms}
For the general situation which
we are interested in, we need a combination of the concepts of log
and log-log forms.

\nnpar{Mixed growth forms.} Let $X$, $\DD$, $U$ and $\iota$ be as in
 the previous 
 section. 
Let $\DD_1$
and $\DD_2$ be normal crossing divisors, which may have common
components, and such
that $\DD=\DD_1\cup \DD_2$. We denote
by $\DD_2'$ the union of the components of $\DD_2$ that are not contained
in $\DD_1$. We say that an
open coordinate subset $V$, with coordinates $z_{1}, \dots , z_{d}$,
 is adapted to $\DD_1$ and $\DD_2$ if $\DD_1\cap V$ has 
equation $z_{1}\dots z_{k}=0$, $\DD'_2\cap V$ has equation
 $z_{k+1}\dots z_{l}=0$ 
and $r_{i}:=|z_{i}|<1/e^{e}$ for
$i=1,\dots ,d$. 

\begin{definition} \label{def:mixedgrowthform}
  Let $V$ be a coordinate neighborhood adapted to $\DD_{1}$ and
  $\DD_{2}$.  For every integer $K\ge 0$, we say that a smooth complex
  function $f$ on $V\setminus \DD$ has \emph{log growth along
    $\DD_{1}$ and log-log growth along $\DD_{2}$ of order $K$}, if
  there exists an integer $N_{K}\ge 0$ such that, for every pair of
  multi-indices $\alpha, \beta \in \mathbb{Z}_{\ge 0}^{d}$, with
  $|\alpha +\beta |\le K$, it holds the inequality
  \begin{equation}
    \left |\frac{\partial^{|\alpha |}}{\partial z^{\alpha }}
      \frac{\partial^{|\beta |}}{\partial \bar z^{\beta }}
      f(z_{1},\dots ,z_{d})\right | \prec \frac{ \left
        | \prod_{i=1}^{k}\log(1/r_{i})
      \prod_{j=k+1}^{l}\log(\log(1/r_{j}))
      \right |^{N_{K}}}{|z^{\alpha^{\le l}}\bar z^{\beta ^{\le l}}|}.
  \end{equation}
  We say that $f$ 
  has \emph{log growth along $\DD_{1}$ and log-log growth along
  $\DD_{2}$ of infinite order}, if it has log growth along
    $\DD_{1}$ and log-log growth along $\DD_{2}$ of order $K$ for all
  $K\ge 0$.
  The \emph{sheaf of differential forms on $X$ with log growth along
    $\DD_{1}$ and 
    log-log growth along $\DD_{2}$ of infinite order} is the subalgebra of
  $\iota _{\ast}\mathscr 
  {E}^{\ast}_{U}$ generated, in each coordinate neighborhood $V$
  adapted to $\DD_{1}$ and $\DD_{2}$, by the functions with log growth
  along $\DD_{1}$ and log-log growth along $\DD_{2}$,
  and the differentials
  \begin{alignat*}{2}
    &\frac{\dd z_{i}}{z_{i}},\,\frac{\dd\bar{z}_{i}}
    {\bar{z}_{i}},&\qquad\text{for }i&=1,\dots,k, \\
    &\frac{\dd z_{i}}{z_{i}\log(1/r_{i})},\,\frac{\dd\bar{z}_{i}}
    {\bar{z}_{i}\log(1/r_{i})},&\qquad\text{for }i&=k+1,\dots,l, \\
    &\dd z_{i},\,\dd\bar{z}_{i},&\qquad\text{for }i&=l+1,\dots,d.
  \end{alignat*}
  When the normal crossing divisors $\DD_{1}$ and $\DD_{2}$ are clear
  from the context, a differential form with log growth along
  $\DD_{1}$ and log-log 
  growth along $\DD_{2}$ of infinite 
  order will be called 
  a \emph{mixed growth form}. The sheaf of differential forms on $X$
  with log growth along $\DD_{1}$ and log-log growth along $\DD_{2}$
  of infinite order will be denoted  
  $\mathscr{E}^{\ast}_{X}
  \langle\DD_{1} \langle \DD_{2}\rangle\rangle_{\growth}$.
\end{definition}

It is clear that
  \begin{equation}
  \label{eq:deflll}
  \mathscr{E}^{\ast
  }_{X}\langle \DD_1 \rangle \wedge \mathscr{E}^{\ast
  }_{X}\langle \langle \DD_2\rangle \rangle_{\growth} \subseteq
  \mathscr{E}^{\ast
  }_{X}\langle \DD_1\langle \DD_2\rangle \rangle_{\growth}. 
\end{equation}
Observe moreover that, by definition,
\begin{displaymath}
\mathscr{E}^{\ast
}_{X}\langle \DD_1\langle \DD_2\rangle \rangle_{\growth}=\mathscr{E}^{\ast
}_{X}\langle \DD_1\langle \DD_2'\rangle \rangle_{\growth}.
\end{displaymath}

We leave it to the reader to state the analogue of lemma \ref{lemm:1}.

\nnpar{Partial differentials.} Let $V$ be an open coordinate system
  adapted to $\DD_{1}$ and $\DD_{2}$.
  In this coordinate system we may decompose the operators $
  \partial$  and $\bar \partial$ as
  \begin{equation}
    \label{eq:9}
    \partial=\sum_{j} \partial_{j}\quad\text{and}\quad
    \bar \partial=\sum_{j} \bar \partial_{j} 
  \end{equation}
  where $\partial_{j}$ and $\bar \partial_{j}$ contain only the
  derivatives with 
  respect to $z_{j}$. 

The following lemma follows directly from the definition.
\begin{lemma} \label{lemm:13}
  Let $E_{j}$ denote the divisor $z_{j}=0$. If $\omega \in
  \mathscr{E}^{\ast}_{X}\langle \DD_{1}\langle
  \DD_{2}\rangle\rangle_{\growth}(V)$, then 
  \begin{displaymath}
    \partial _{j}\omega \in 
    \begin{cases}
      \mathscr{E}^{\ast}_{X}\langle \DD_{1}\langle
      \DD_{2}\rangle\rangle_{\growth}(V),\ &\text{if }j\le k\text{ or
      }j>l,\\   
      \mathscr{E}^{\ast}_{X}\langle \DD_{1}\cup E_{j}\langle
      \DD_{2}\rangle\rangle_{\growth}(V),\ &\text{if }k<j\le l,
    \end{cases}
  \end{displaymath}
  and the same is true for the operator $\bar \partial_{j}$.
\end{lemma}
\hfill $\square$

\nnpar{Mixed Forms.}
\begin{definition}  \label{def:mixedforms} We say that a section
  $\omega $ of $\iota _{\ast}\mathscr{E}^{\ast}_{U}$
  is \emph{log along $\DD_{1}$ and log-log
  along $\DD_{2}$} if the differential forms $\omega $, $\partial \omega $,
  $\bar \partial \omega $ and $\partial \bar \partial \omega $ are
  sections of $\mathscr{E}^{\ast }_{X}\langle \DD_1\langle
  \DD_2\rangle \rangle_{\growth}$.
  The sheaf of differential forms log along $\DD_{1}$ and log-log
  along $\DD_{2}$ will be denoted by $\mathscr{E}^{\ast}_{X}
  \langle\DD_{1} \langle \DD_{2}\rangle\rangle$. 
  For shorthand, a differential
  form which is log along $\DD_{1}$ and log-log along $\DD_{2}$ will
  be called a mixed form. 
\end{definition}

As the complexes we have defined in the previous sections, the complex
$\mathscr{E}^{\ast 
}_{X}\langle \DD_{1}\langle \DD_{2}\rangle \rangle$ is a sheaf of Dolbeault
algebras.

\nnpar{Inverse images.} 
We can generalize propositions \ref{prop:4} and \ref{prop:10}, with
the same technique, to the case of mixed forms.

\begin{proposition}
\label{prop:invimagell} 
Let $f:X\longrightarrow Y$ be a morphism of complex
manifolds. Let $\DD_{1}$, $\DD_{2}$ and $E_{1}$, $E_{2}$ be normal 
crossing divisors on $X$ and  $Y$ respectively, such that $\DD_{1} 
\cup \DD_{2}$ and $E_{1}\cup E_{2}$ are also normal crossing divisors.
Furthermore, assume that $f^{-1}(E_{1})\subseteq \DD_{1}$ and $f^
{-1}(E_{2})\subseteq \DD_{1}\cup \DD_{2}$. If $\eta$ is a section of
$\mathscr{E}^{\ast}_{Y}\langle E_{1}\langle E_{2}\rangle\rangle$, 
then $f^{\ast}\eta$ is a section of $\mathscr{E}^{\ast}_
{X}\langle \DD_{1}\langle \DD_{2}\rangle\rangle$.  
\hfill $\square$
\end{proposition}

\nnpar{Integrability.} Let $X$ be a complex manifold and $\DD$ a
normal crossing divisor.
Let $y$ be a $p$-codimensional cycle of $X$ and let $Y=\supp y$. 
Let $\pi:\widetilde{X}\longrightarrow X$ be 
an embedded resolution of singularities of $Y$, with normal 
crossing divisors $\DD_{Y}=\pi^{-1}(Y)$ and $\widetilde{\DD}=\pi^
{-1}(\DD)$ and such that $\DD_{Y}\cup \widetilde {\DD}$ is also a normal
crossing divisor.

\begin{lemma}
\label{lemm:3} 
Assume that $g\in\Gamma(\widetilde{X},\mathscr{E}^{n}_{\widetilde
{X}}\langle \DD_{Y}\langle\widetilde{\DD}\rangle\rangle)$. Then,
the following statements hold: 
\begin{enumerate}
\item[(i)] 
If $n<2p$, then $g$ is locally integrable on the whole of $X$. 
We denote by $[g]_{X}$ the current associated to $g$.
\item[(ii)]
If $n<2p-1$, then $\dd [g]_{X}=[\dd g]_{X}$. 
\end{enumerate}
\end{lemma}
\begin{proof}
This is a particular case of \cite{BurgosKramerKuehn:cacg} lemma 7.13.
\end{proof}

\nnpar{The cohomology of the complex of mixed forms.}
We are now in position to state the 
main result of this section:

\begin{theorem} \label{thm:fq}  The inclusion 
\begin{displaymath}
  \begin{CD}
    \Omega ^{\ast }_{X}(\log D_1) @>>>\mathscr{E}^{\ast }_{X}\langle \DD_1
    \langle \DD_2 \rangle\rangle  
  \end{CD}
\end{displaymath}
is a filtered quasi-isomorphism with respect to the Hodge filtration.
\end{theorem}
\begin{proof}
  To prove the theorem we will use 
  the classical argument for proving the Poincar\'e lemma in many
  variables. We will state here the general argument and we will delay
  the specific analytic lemmas that we need until the next section. 

  The theorem is equivalent to the exactness of the sequence of
  sheaves 
  \begin{displaymath}
    \begin{CD}
      0@>>> \Omega ^{p}_{X}(\log \DD_1) @ > i >>
      \mathscr{E}^{p,0 }_{X}\langle \DD_1 \langle \DD_2 \rangle\rangle
      @> \bar \partial >>
      \mathscr{E}^{p,1 }_{X}\langle \DD_1 \langle \DD_2 \rangle\rangle
      @> \bar \partial >> \dots 
    \end{CD}
  \end{displaymath}
  The exactness in the first step is clear because a holomorphic form
  on $X\setminus (\DD_1\cup \DD_2)$ that satisfies the growth conditions
  imposed in the definitions can only have logarithmic poles along
  $\DD_1$. 

  For the exactness in the other steps we choose a point $x\in X$. Let
  $V$ be 
  a coordinate neighborhood of $x$ adapted to $\DD_1$ and $\DD_2$, and such
  that $x$ has coordinates $(0,\dots ,0)$.

  Let $0<\epsilon \ll 1$, we denote by $\Delta ^{d}_{x,\epsilon }$ the
  poly-cylinder 
  \begin{displaymath}
    \Delta _{x,\epsilon }^{d}=\{(z_{1},\dots ,z_{d})\in V\mid
    r_{i}<\epsilon ,\ i=1,\dots ,d\}.
  \end{displaymath}

  In the next section we will prove that, for $j=1,\dots ,d$ and
  $0<\epsilon '<\epsilon \ll 1$, there 
  exist operators
  \begin{displaymath}
    \begin{matrix}
      K^{\epsilon ',\epsilon }_{j}&:
      \mathscr{E}^{p,q }_{X}\langle \DD_1 \langle \DD_2
      \rangle\rangle(\Delta ^{d}_{x,\epsilon }) &\longrightarrow
      \mathscr{E}^{p,q-1 }_{X}\langle \DD_1 \langle \DD_2 
      \rangle\rangle(\Delta ^{d}_{x,\epsilon '}),\\
      P^{\epsilon ',\epsilon }_{j}&:\mathscr{E}^{p,q }_{X}
      \langle \DD_1 \langle \DD_2
      \rangle\rangle(\Delta ^{d}_{x,\epsilon }) &\longrightarrow
      \mathscr{E}^{p,q 
        }_{X}\langle \DD_1 \langle \DD_2 \rangle\rangle(\Delta ^{d}
      _{x,\epsilon '}), 
    \end{matrix}    
  \end{displaymath}
  that satisfy the following conditions
  \begin{itemize}
  \item[(A)] \label{prop:1} If the form $\omega $ does not contain any term
    with $\dd\bar z_{i}$ for $i>j$ then $K^{\epsilon ',\epsilon
    }_{j}\omega $ 
    and $P^{\epsilon ',\epsilon }_{j}\omega $ do not contain any term
    with $\dd\bar z_{i}$, for $i\ge j$.
  \item[(B)] 
$\bar \partial K^{\epsilon ',\epsilon }_{j}
    +K^{\epsilon ',\epsilon }_{j}\bar
    \partial+P^{\epsilon ',\epsilon }_{j}=\Id$.
  \end{itemize}
  Let $q>0$ and let $\omega \in \mathscr{E}^{p,q }_{X}\langle \DD_1
  \langle \DD_2 \rangle\rangle_{x}$ be a germ of a closed form. Assume
  that $\omega $ is defined in a poly-cylinder $\Delta
  ^{d}_{x,\epsilon 
  }$. By abuse of notation we will not distinguish between sections
  and germs. Therefore $\omega $ will denote also a closed
  differential form over $\Delta ^{d}_{x,\epsilon
  }$ that represents this germ.  Moreover, as the open set of
  definition of each section will be clear from
  the context we will not make it explicit. We choose real numbers
  $0<\epsilon _{1}<\dots <\epsilon 
  _{d}<\epsilon $. 
  Then, by property  
(B), we have
  \begin{displaymath}
    \omega =
    \bar \partial K^{\epsilon ,\epsilon _{d}}_{d}(\omega )+P^{\epsilon
    ,\epsilon _{d}}_{d}(\omega ). 
  \end{displaymath}
  We write $\omega _{1}=P^{\epsilon ,\epsilon _{d}}_{d}(\omega )$. 
  Then $\omega_{1} $ does not
  contain $d\bar z_{d}$ and $\omega -\omega _{1}$ is a boundary.
  We define inductively $\omega _{j+1}=P^{\epsilon_{d-j+1},\epsilon
  _{d-j}}_{d-j}(\omega 
  _{j})$. Then, for all $j$, $\omega -\omega _{j}$ is a boundary and
  $\omega _{j}$ does not contain $\dd\bar z_{i}$, for $i>
  d-j$. Therefore $\omega _{d-q+1}=0$ and $\omega $ is a boundary.
\end{proof}
 
\subsection{Analytic lemmas}
\label{sec:al}

In this section we will prove the analytic lemmas needed to prove
theorem \ref{thm:fq} and we will define the operators $K$ and $P$
that appear in the proof of this theorem.

\nnpar{Primitive functions with growth conditions.}
Let $f$ be a smooth function on $\Delta ^{\ast}_{\epsilon }$, which is
integrable on any compact subset of $\Delta _{\epsilon }$. Then, for
$\epsilon '<\epsilon $ and $z\in \Delta ^{\ast}_{\epsilon '}$, we
write 
\begin{displaymath}
  I_{\epsilon '}(f)(z)=\frac{1}{2\pi \sqrt{-1}} \int_{\overline
  {\Delta} _{\epsilon '}}  
  f(w)\frac{\dd w\land \dd\bar w} 
  {w-z}.  
\end{displaymath}
We denote $r=| z|$.

\begin{lemma}\label{lemm:14}
  \begin{enumerate}
  \item 
If $f$ is a smooth function on $\Delta
    ^{\ast}_{\epsilon }$, 
    such that 
    \begin{displaymath}
      |f(z)|\prec \frac{|\log(\log (1/r))|^{N}}{(r\log (1/r))^{2}} 
    \end{displaymath}
    then $f$ is integrable in each compact subset of $\Delta _{\epsilon
    }$ and
    \begin{displaymath}
      \frac{\partial}{\partial \bar z} I_{\epsilon '}(f)(z)=f(z). 
    \end{displaymath}
  \item 
If $f$ is a smooth function on $\Delta
    ^{\ast}_{\epsilon }$, such that
    \begin{displaymath}
      |f(z)|\prec \frac{|\log(\log (1/r))|^{N}}{r\log (1/r)}\text{ and }
      \left |\frac{\partial }{\partial \bar z} f(z) \right|\prec
      \frac{|\log(\log (1/r))|^{N}}{(r\log (1/r))^{2}},
    \end{displaymath}
    then
    \begin{displaymath}
      2\pi \sqrt{-1} f(z)=\int_{\partial \Delta _{\epsilon '}}
      f(w)\frac{\dd w}{w-z}
      +\int_{\overline
        {\Delta} _{\epsilon '}}  
      \frac{\partial}{\partial\bar w} f(w)\frac{\dd w\land \dd\bar w} 
      {w-z}.  
    \end{displaymath}
  \item 
If $f$ is a smooth function on $\Delta
    ^{\ast}_{\epsilon }$, such that
    \begin{displaymath}
      |f(z)|\prec \frac{|\log(\log (1/r))|^{N}}{r\log (1/r)}\text{ and }
      \left |\frac{\partial}{\partial z} f(z) \right|\prec
      \frac{|\log(\log (1/r))|^{N}}{(r\log (1/r))^{2}},
    \end{displaymath}
    then
    \begin{displaymath}
      \frac{\partial}{\partial z}  \int_{\overline
        {\Delta} _{\epsilon '}}  
      f(w)\frac{\dd w\land \dd\bar w} 
      {w-z}
      =
      -\int_{\partial \Delta _{\epsilon '}} f(w)\frac{\dd \bar w}{w-z} 
      + \int_{\overline
        {\Delta} _{\epsilon '}}  
      \frac{\partial}{\partial w} f(w)\frac{\dd w\land \dd\bar w} 
      {w-z}.  
    \end{displaymath}
  \end{enumerate}
\end{lemma}
\begin{proof}
     We start by proving the integrability of
  $f$.  Viewed as a function of $\epsilon $, 
  \begin{align*}
    \left |\int_{\Delta _{\epsilon }}\frac{(\log (\log
    (1/r)))^{N}}{r^{2}(\log (1/r))^{2}}\dd z\land \dd\bar z\right
    |&\prec 
    \left | \int_{0}^{\epsilon } \frac{(\log (\log
        (1/r)))^{N}}{r^{2}(\log (1/r))^{2}}r \dd r \right |\\
    &\prec \left | \int_{0}^{\epsilon } \frac{1}{r(\log (1/r))^{3/2}}
    \dd r
    \right | \\
    &\prec \frac{1}{(\log (1/\epsilon) )^{1/2}}, 
  \end{align*}
  Which proves the integrability. Then the formulas
  are proved as in \cite{GunningRossi:afscv} pag. 24-26. The only new point one
  has to care about is that the singularities at $z=0$ do not contribute to
  the Stokes theorem.
\end{proof}

\begin{lemma}
\label{lemm:growth} Let $0<\epsilon \ll 1$ be a real
number and let $f$ be a smooth function on $\Delta
^{\ast}_{\epsilon }$. Let $\epsilon '<\epsilon $. 
\begin{enumerate}
\item 
If $\omega =f\,\dd\bar z\in
   \mathscr{E}^{0,1}_{\Delta _{\epsilon }}\langle 0 \rangle(\Delta
   _{\epsilon })$, 
  then the function $f$ is integrable on any compact subset of
  $\Delta _{\epsilon }$. We write $g=I_{\epsilon '}(f)$. 
  Then $g\in \mathscr{E}^{0,0}_{\Delta _{\epsilon '}}\langle 0
   \rangle(\Delta _{\epsilon '})$
  and
  \begin{equation}\label{eq:dg}
    \bar {\partial} g=\omega .  
  \end{equation}
\item  
If moreover $\omega \in
      \mathscr{E}^{0,1}_{\Delta _{\epsilon 
      }}\langle 
  \langle 0 \rangle \rangle_{\growth}(\Delta
   _{\epsilon })$, then $g\in \mathscr{E}^{0,0}_{\Delta
    _{\epsilon '}}\langle \langle 0 \rangle \rangle_{\growth}(\Delta
   _{\epsilon '})$.
 \item \label{lemm:1c} If $\omega =f\,\dd\bar z\land \dd z\in
   \mathscr{E}^{1,1}_{\Delta _{\epsilon }}\langle \langle 0 \rangle
   \rangle_{\growth}(\Delta _{\epsilon })$, then the function $f$ is
   integrable on any compact subset of $\Delta _{\epsilon }$. If we
   write $g=I_{\epsilon '}(f)$ as before, then $g\,\dd z\in
   \mathscr{E}^{1,0}_{\Delta _{\epsilon '}}\langle\langle 0
   \rangle\rangle_{\growth}(\Delta
   _{\epsilon '})$ and
  \begin{equation}\label{eq:dgb}
    \bar {\partial} (g \land \dd z)=\omega.  
  \end{equation}
\end{enumerate}
\end{lemma}
\begin{proof}
  The integrability in the three cases and equations \eqref{eq:dg} and
  \eqref{eq:dgb} are in lemma \ref{lemm:14}. Therefore it remains only
  to prove the necessary bounds.

  \emph{Proof of 
(i).} The condition
  on $\omega $ is equivalent to the inequalities  
  \begin{equation}
    \left | \frac{\partial^{\alpha +\beta }}{\partial z^{\alpha
    }\partial \bar z^{\beta }} f(z)\right | \prec
    \frac{\left| \log (1/r)\right|^{N_{\alpha +\beta }}}
    {r^{\alpha +\beta +1}},
  \end{equation}
  for certain family of  integers $\{N_{n}\}_{n\in \mathbb{Z}_{\ge
  0}}$. We may assume that these integers satisfy that, if $a\le b$
  then $N_{a}\le N_{b}$.
  We can apply \cite{HarrisPhong:cdcli}
  lemma 1 to conclude that $g$ is smooth on
  $\Delta ^{\ast}_{\epsilon '}$ and that 
  \begin{displaymath}
    \left | g(z) \right | \prec \left | \log (1/r)\right |^{N'_{0}}, 
  \end{displaymath}
  for some integer $N'_{0}$. 

  Thus to prove the statement 
(i) it remains to bound the
  derivatives of $g$. The equation \eqref{eq:dg} implies the bound of
  the derivatives 
  \begin{displaymath}
    \frac{ \partial^{\alpha +\beta }}{\partial z^{\alpha
    }\partial\bar z^{\beta }}g
  \end{displaymath}
  when $\beta \ge 1$. Therefore we may assume $\beta =0$ and $\alpha
  \ge 1$.

  Let $\rho :\mathbb{C}\longrightarrow [0,1]$ be a smooth function such that  
  \begin{displaymath}
    \rho |_{B(0,1 )}=1,\qquad \rho |_{\mathbb{C}\setminus B(0,2)} =0,
  \end{displaymath}
  where $B(p,\delta )$ is the open ball of center $p$ and radius $\delta $.
  Fix $z_{0}\in \Delta ^{\ast}_{\epsilon '}$.
  Since we want to bound the derivatives of $g(z)$ as $z$ goes to
  zero, we may assume $z_{0}\in \Delta ^{\ast}_{\epsilon '/2}$.
  Write $r_{0}=|z_{0}|$, and put 
  \begin{displaymath}
    \rho _{z_{0}}(z)=\rho \left (3\frac{z-z_{0}}{r_{0}} \right ).
  \end{displaymath}
  Then
  \begin{displaymath}
    \rho |_{B(z_{0},r_{0}/3 )}=1,\qquad \rho
    |_{\mathbb{C}\setminus B(z_{0},2r_{0}/3)} =0.
  \end{displaymath}
  Moreover 
  \begin{equation}\label{eq:rho}
    \frac{\partial^{\alpha } }
    {\partial z^{\alpha }}
    \rho _{z_{0}}(z) \le \frac{C_{\alpha }}{r_{0}^{\alpha}},
  \end{equation}
  for some constants $C_{\alpha }$.

  By the choice of $z_{0}$, we have that $\supp(\rho _{z_{0}})\subset
  \Delta ^{\ast}_{\epsilon '}$. 
  We write $f_{1}=\rho_{z_{0}} f$ and $f_{2}=(1-\rho_{z_{0}}
  )f$. Then, for $z\in B(z_{0},r_{0}/3)$, we introduce the auxiliary
  functions 
  \begin{displaymath}
    g_{1}(z)=\frac{1}{2\pi \sqrt {-1}}
    \int_{\Delta _{\epsilon }}f_{1}(w)
    \frac{\dd w\land \dd\bar w}{w-z},\quad      
    g_{2}(z)=\frac{1}{2\pi \sqrt {-1}}
    \int_{\Delta _{\epsilon }}f_{2}(w)
    \frac{\dd w\land \dd\bar w}{w-z}.      
  \end{displaymath}
  These functions satisfy
  \begin{displaymath}
    g=g_{1}+g_{2}.
  \end{displaymath}
  Therefore we can  bound separately the derivatives of $g_{1}$ and
  $g_{2}$. We first bound the derivatives of $g_{1}$.
  \begin{align*}
    \frac{\partial^{\alpha}}
      {\partial z^{\alpha }}
      \int_{\Delta _{\epsilon }}f_{1}(w)
      \frac{\dd w\land \dd\bar w}{w-z} &=
    \frac{\partial^{\alpha }}
      {\partial z^{\alpha }}
      \int_{\mathbb{C}}f_{1}(w)
      \frac{\dd w\land \dd\bar w}{w-z} \\
    &=\frac{\partial^{\alpha}}
      {\partial z^{\alpha }}
      \int_{\mathbb{C}}f_{1}(u+z)
      \frac{\dd u\land \dd\bar u}{u}\\
    &=\int_{\mathbb{C}} \frac{\partial^{\alpha}}
      {\partial z^{\alpha }}
      f_{1}(u+z)
      \frac{\dd u\land \dd\bar u}{u}\\
    &=\int_{\mathbb{C}}\frac{\partial^{\alpha }}
      {\partial w^{\alpha }}
      f_{1}(w)
      \frac{\dd w\land \dd\bar w}{w-z}\\
    &=\int_{B(z_{0},2r_{0}/3)}\frac{\partial^{\alpha }}
      {\partial w^{\alpha }}
      f_{1}(w)
      \frac{\dd w\land \dd\bar w}{w-z}.
  \end{align*}
  Hence using the bound of the derivatives of $f$ and equation
  \eqref{eq:rho}, it holds the inequality
  \begin{align*}
    \left|\frac{\partial^{\alpha } }
      {\partial z^{\alpha }}g_{1}(z_{0})\right| 
    &\prec \frac{\left | \log (1/r_{0})\right |^{N_{\alpha }
          }}{r_{0}^{\alpha +1}} 
    \left| \int_{B(z_{0},2r_{0}/3)}
      \frac{\dd w\land \dd\bar w}{w-z_{0}} \right | \\
    &\prec \frac{\left | \log (1/r_{0})\right |^{N_{\alpha }
          }}{r_{0}^{\alpha }}.
  \end{align*}
  Now we bound the derivatives of $g_{2}$. Since for $z\in
  B(z_{0},r_{0}/3)$, the function $f_{2}(w)$ is
  identically
  zero in a neighborhood of the point $w=z$, we have that
  \begin{displaymath}
    \frac{\partial^{\alpha }}
      {\partial z^{\alpha }}
      \int_{\Delta _{\epsilon }}f_{2}(w)
      \frac{\dd w\land \dd\bar w}{w-z} =
      \int_{\Delta _{\epsilon }}f_{2}(w) \alpha !
      \frac{\dd w\land \dd\bar w}{(w-z)^{\alpha  +1}}.
    \end{displaymath}
  Let $A=B(0,r_{0}/2)$. Then, for $w\in A$ we have $|w-z_{0}|\ge
  r_{0}/2$. Thus
  \begin{align*}
    \left| \int_{A} f_{2}(w) \frac{\dd w\land \dd \bar
        w}{(w-z_{0})^{\alpha +1}} \right |&\prec \frac{1}{r_{0}^{\alpha 
        +1}}\int_{0}^{r_{0}/2} \frac{|\log (1/\rho) |^{N_{0}}}{\rho } \rho
        \dd\rho \\
        &\prec \frac{1}{r_{0}^{\alpha }}|\log (1/r_{0}) |^{N_{0}}. 
  \end{align*}
  Here we use that 
  \begin{equation*}\label{eq:intlog}
    \int (\log x)^{N}\dd x=x\sum_{i=0}^{N} (-1)^{i}\frac{N !}{(N-i)!} 
    (\log x)^{N-i}.  
  \end{equation*}
  We write $B=\Delta _{\epsilon '}\setminus (A\cup B(z_{0},r_{0}/3))$. In
  this region $|w-z_{0}|\ge |w/4 |$. Therefore 
  \begin{align*}
    \left| \int_{B} f_{2}(w) \frac{\dd w\land \dd \bar
        w}{(w-z_{0})^{\alpha +1}} \right |\prec
    \int_{r_{0}/2}^{\epsilon '}\frac{|\log (1/\rho) |^{N_{0}}}{\rho }
    \frac{\rho \dd\rho }{\rho ^{\alpha +1}} \prec
        \frac{1}{r_{0}^{\alpha }}|\log (1/r_{0}) |^{N_{0}+1}.   
  \end{align*}
  Here we use that 
  \begin{equation*}
    \label{eq:intlog2}
    \int(\log x)^{n}\frac{1}{x^{m}}\dd x=
    \begin{cases}
      \frac{1}{n+1} (\log x)^{n+1}&\text{ if }m=1,\\
      \frac{1}{x^{m-1}}P_{n,m}(\log x)&\text{ if }m>1,  
    \end{cases}
  \end{equation*}
  where $P_{n,m}$ is a polynomial of degree $n$. Summing up, we obtain
  \begin{displaymath}
    \left| \frac{\partial^{\alpha }}{\partial z^{\alpha }}g(z_{0}) \right|
    \prec
    \frac{|\log (1/r_{0})|^{N_{\alpha }+1}}{r_{0}^{\alpha }} 
  \end{displaymath}
 Observe that, for $\alpha =0$, this
  is the proof of \cite{HarrisPhong:cdcli} lemma 1.
  
  \emph{Proof of (ii)
.} In this case, by lemma \ref{lemm:1},
  the condition on $\omega 
  $ is equivalent to the inequalities 
  \begin{equation}
    \left | \frac{\partial^{\alpha +\beta }}{\partial z^{\alpha
    }\partial \bar z^{\beta }} f(z)\right | \prec
    \frac{\left| \log (\log(1/ r))\right|^{N_{\alpha +\beta }}}
    {r^{\alpha +\beta +1}\log (1/r)},
  \end{equation}
  for certain increasing family of
  integers $\{N_{n}\}_{n\in \mathbb{Z}_{\ge 0}}$ . Again by lemma
  \ref{lemm:1}, to prove statement (ii)
, 
  we have to show  
  \begin{equation}
    \left | \frac{\partial^{\alpha +\beta }}{\partial z^{\alpha
    }\partial \bar z^{\beta }} g(z)\right | \prec
    \frac{\left| \log (\log(1/ r))\right|^{N'_{\alpha +\beta }}}
    {r^{\alpha +\beta}},
  \end{equation}
  for certain  family of integers $\{N'_{n}\}_{n\in \mathbb{Z}_{\ge 0}}$. 
  
  By \eqref{eq:dg} the functions  
  \begin{displaymath}
    \frac{ \partial^{\alpha +\beta }}{\partial z^{\alpha
    }\partial\bar z^{\beta }}g,
  \end{displaymath}
  when $\beta \ge 1$, satisfy the required bounds.
  Thus it remains to
  bound $\partial ^{\alpha }/\partial z^{\alpha }g$ for $\alpha \ge 0$. 
  As in the proof of statement (i) 
  we fix $z_{0}$ and we write $g=g_{1}+g_{2}$. For $g_{1}$ we work as
  before and we get, for $\alpha \ge 0$,
  \begin{displaymath}
    \left|\frac{\partial^{\alpha } }
      {\partial z^{\alpha }}g_{1}(z_{0})\right| 
    \prec \frac{\left | \log (\log (1/r_{0}))\right |^{N_{\alpha } 
          }}{r_{0}^{\alpha }\log(1/r_{0})}.
  \end{displaymath}

  To bound $g_{2}$ we integrate over the regions $A$ and
  $B$ as before. We first bound the integral over the region
  $A=B(0,r_{0}/2)$. 
  \begin{displaymath}
        \left| \int_{A} f_{2}(w) \frac{\dd w\land \dd \bar
        w}{(w-z_{0})^{\alpha +1}} \right |\prec \frac{1}{r_{0}^{\alpha 
        +1}}\int_{0}^{r_{0}/2} \frac{|\log (\log (1/\rho)) |^{N_{0}
      }}{\log (1/\rho )} \dd\rho .
  \end{displaymath}
  Since
  for $\rho <1/e^{e^{N_{0}}}$ the function 
  \begin{displaymath}
    \frac{(\log (\log (1/\rho)) )^{N_{0}}}{\log (1/\rho) }  
  \end{displaymath}
  is an increasing function, then 
  \begin{align*}
    \frac{1}{r_{0}^{\alpha 
        +1}}\int_{0}^{r_{0}/2} \frac{|\log (\log (1/\rho)) |^{N_{0}
      }}{\log (1/\rho )} \dd\rho 
    &\prec \frac{|\log (\log (1/r_{0} ))|^{N_{0}
      }}{r_{0}^{\alpha 
        +1}\log (1/r_{0})}
    \int_{0}^{r_{0}/2}  \dd\rho \\
    &\prec  \frac{|\log (\log (1/r_{0})) |^{N_{0}
      }}{r_{0}^{\alpha } \log (1/r_{0} )}. 
  \end{align*}
  in the domain $0<r_{0}\le 2/e^{e^{N_{0}}}$.

  If $f$ and $g$ are two continuous functions, with $g$
  strictly positive,
  defined on a compact set, then $f\prec g$. Therefore the above
  inequality extends to the domain $0\le r_{0}\le \epsilon '/2$.

  We now bound the integral over the region $B=\Delta _{\epsilon
  '}\setminus (A\cup B(z_{0},r_{0}/3))$.  By the bound of the
  function $f$, we have  
    \begin{align*}
    \left| \int_{B} f_{2}(w) \frac{\dd w\land \dd \bar
        w}{(w-z_{0})^{\alpha +1}} \right |\prec
    \int_{r_{0}/2}^{\epsilon '}\frac{|\log (\log (1/\rho)) |^{N_{0}
          }}
    {\rho ^{\alpha +1}\log (1/\rho) }
    \,\dd\rho.   
  \end{align*}
  Thus, in the case $\alpha =0$,  
  \begin{align*}
    \int_{r_{0}/2}^{\epsilon '}\frac{|\log (\log (1/\rho)) |^{N_{0}
          }}
    {\rho \log (1/\rho )}
    \,\dd\rho
    \prec
    |\log (\log (1/r_{0})) |^{N_{0}+1}.  
  \end{align*}
  In the case $\alpha >0$, since for $\rho
  <1/e^{e}$ the function 
  \begin{displaymath}
    \frac{(\log (\log (1/\rho)) )^{N_{0}}}{\rho ^{1/2}\log (1/\rho )} 
  \end{displaymath}
  is a decreasing function, we have
  \begin{align*}
    \int_{r_{0}/2}^{\epsilon '}\frac{|\log (\log (1/\rho ))|^{N_{0}
          }}
    {\rho ^{\alpha +1}\log (1/\rho) }
    \,\dd\rho &\prec
    \frac{|\log (\log (1/r_{0})) |^{N_{0}
          }}
    {r_{0} ^{1/2}\log (1/r_{0}) }
    \int_{r_{0}/2}^{\epsilon '}
    \frac{1}{\rho ^{\alpha +1/2}} \,\dd\rho\\
    &\prec \frac{|\log (\log (1/r_{0})) |^{N_{0}
        }}{r_{0}^{\alpha } \log (1/r_{0} )}.
  \end{align*}

  Summing up, we obtain
  \begin{equation*}
    \left |\frac{\partial^{\alpha }}{\partial z^{\alpha }}g(z_{0})\right
    |
    \prec \frac{|\log (\log (1/r _{0}))|^{N_{\alpha }+1
        }}{r_{0}^{\alpha }}.
  \end{equation*}
  This finishes the proof of the second statement. 
  
  \emph{Proof of (iii)
.} In this case, by lemma
  \ref{lemm:1} once more, the condition on $\omega 
  $ is equivalent to the conditions
  \begin{equation*}
    \left | \frac{\partial^{\alpha +\beta }}{\partial z^{\alpha
    }\partial \bar z^{\beta }} f(z)\right | \prec
    \frac{\left| \log (\log(1/ r))\right|^{N_{\alpha +\beta }}}
    {r^{\alpha +\beta +2}(\log (1/r))^{2}},
  \end{equation*}
 for certain increasing family of
  integers $\{N_{n}\}_{n\in \mathbb{Z}_{\ge 0}}$, and the inequalities
  we have to prove are 
  \begin{equation*}
    \left | \frac{\partial^{\alpha +\beta }}{\partial z^{\alpha
    }\partial \bar z^{\beta }} g(z) \right | \prec
    \frac{\left| \log (\log(1/ r))\right|^{N'_{\alpha +\beta }}}
    {r^{\alpha +\beta+1} \log (1/r)},
  \end{equation*}
  for certain family of integers $\{N'_{n}\}_{n\in \mathbb{Z}_{\ge
  0}}$.  

  First we note that, by equation \eqref{eq:dgb}, for $\beta \ge 1$,
  the functions 
  \begin{displaymath}
    \frac{ \partial^{\alpha +\beta }}{\partial z^{\alpha
    }\partial\bar z^{\beta }}g
  \end{displaymath}
  satisfy the required bounds. Thus it remains to bound the functions
  $\partial^{\alpha }g/\partial z^{\alpha }$, for $\alpha \ge 0$. The
  proof is similar as before. We decompose again $g=g_{1}+g_{2}$. In
  this case 
  \begin{displaymath}
    \left|\frac{\partial^{\alpha } }
      {\partial z^{\alpha }}g_{1}(z_{0})\right| 
    \prec \frac{\left | \log (\log (1/r_{0}))\right |^{N_{\alpha } 
          }}{\log(1/r_{0})^{2}r_{0}^{\alpha +1}}.
  \end{displaymath}
  Whereas the integral of $g_{2}$ over $A$ is bounded as
  \begin{align*}
    \left| \int_{A} f_{2}(w) \frac{\dd w\land \dd \bar
        w}{(w-z_{0})^{\alpha +1}} \right |&\prec \frac{1}{r_{0}^{\alpha 
        +1}}\int_{0}^{r_{0}/2} \frac{|\log (\log (1/\rho)) |^{N_{0}
      }}{\rho \log (1/\rho )^{2}} \dd\rho \\
    &\prec  \frac{|\log (\log (1/r_{0})) |^{N_{0}}}{r_{0}^{\alpha +1} \log
        (1/r_{0})},  
  \end{align*}
  and the integral of $g_{2}$ over $B$ is bounded as
  \begin{align*}
    \left| \int_{B} f_{2}(w) \frac{\dd w\land \dd \bar
        w}{(w-z_{0})^{\alpha +1}} \right |&\prec
    \int_{r_{0}/2}^{\epsilon '}\frac{|\log (\log (1/\rho)) |^{N_{0}
      }}
    {\rho ^{\alpha +2}\log (1/\rho)^{2} }
    \,\dd\rho\\
    &\prec
     \frac{|\log (\log (1/r_{0})) |^{N_{0}
      }}{r_{0}^{\alpha +1} \log (1/r_{0} )^{2}}.
  \end{align*}
  Summing up, we obtain that, for $\alpha \ge 0$, 
  \begin{align*}
    \left |\frac{\partial^{\alpha }}{\partial z^{\alpha }}g(z_{0})\right
    |
    &\prec \frac{|\log (\log (1/r_{0})) |^{N_{\alpha }
        }}{r_{0}^{\alpha+1 } \log (1/r _{0})}.
  \end{align*}
  This finishes the proof of the lemma. 
\end{proof}

\begin{remark}\label{rem:2} 
  Observe that, in general, a section of $\mathscr{E}^{1,1}_{\Delta _{\epsilon
  }}\langle 0 \rangle (\Delta _{\epsilon
  })$ is not locally integrable (see remark
  \ref{rem:6}). Therefore the 
  analogue of lemma \ref{lemm:growth} (iii) 
is not true for
  log forms.
\end{remark}

\nnpar{The operators $K$ and $P$.}
Let $X$, $U$, $\DD$, $\iota $, $\DD_1$ and $\DD_2$ be as in definition
\ref{def:mixedforms}. 

\begin{notation}
Let 
$x\in X$.  
Let $V$ be a open coordinate neighborhood of $x$, with coordinates
$z_{1},\dots  ,z_{d}$,  adapted to $\DD_1$ and 
$\DD_2$ and such that $x$ has coordinates $(0,\dots ,0)$. Thus $\DD_1$ has
equation $z_{1}\dots 
z_{k}=0$ and $\DD'_2$ has equation 
$z_{k+1}\dots z_{l}$. 
Once this coordinate neighborhood is chosen we
put 
\begin{alignat*}{2}
  \zeta _{i}&=\frac{\dd z_{i}}{z_{i}},\quad &\text{ if } &1\le i\le k,\\
  \zeta _{i}&=\dd z_{i},\quad &\text{ if } &i> k. 
\end{alignat*}
For any subset $I\subset \{1,\dots ,d\}$ we denote 
\begin{displaymath}
  \zeta _{I}=\bigwedge_{i\in I}\zeta _{i},\quad
  \dd\bar z _{I}=\bigwedge_{i\in I}\dd\bar z _{i}.
\end{displaymath}
Given any differential form $\omega $, let 
\begin{displaymath}
  \omega =\sum_{I,J} f_{I,J}\zeta _{I}\land \dd \bar z_{J} 
\end{displaymath}
be the decomposition of $\omega $ in monomials. Then we write
\begin{displaymath}
  \omega _{I,J}=f_{I,J}\zeta _{I}\land \dd \bar z_{J}.
\end{displaymath}
For any subset $I\subset \{1,\dots ,d\}$ and
$i\in I$ we will write 
$$\sigma (I,i)=\sharp \{j\in I\mid j<i\}\qquad \text{ and }\qquad
I_{i}=I\setminus \{i\}.$$
\end{notation}

\begin{definition} \label{def:12} Let $0<\epsilon '<\epsilon \ll
  1$. Let $\Delta 
  _{x,\epsilon }^{d}$ be the poly-cylinder centered at $x$ of radius
  $\epsilon $. 
Let $\omega \in \mathscr{E}^{p,q }_{X}\langle
\DD_1 \langle \DD_2\rangle\rangle_{\growth}(\Delta ^{d}_{x,\epsilon })$,
  and let     
\begin{equation} \label{eq:1}
  \omega =\sum_{I,J}f_{I,J}\zeta _{I}\land \dd\bar{z}_{J}
\end{equation}
be the decomposition of $\omega $ in monomials.
We define
\begin{multline*}
  K_{j}^{\epsilon ',\epsilon }(\omega )= \sum_{I} (-1)^{\# I}
  \zeta_{I}\land \\ 
\sum_{J\,\mid \,j\in J}
  \frac{(-1)^{\sigma (J,j)}}{2\pi \sqrt{-1}} \int_{\Delta _{\epsilon '}} 
  f_{I,J}(\dots ,z_{j-1},w,z_{j+1},\dots )\frac{\dd w\land \dd\bar w} 
  {w-z_{j}}
  \dd\bar z_{J_{j}},
\end{multline*}
\begin{multline*}
  P_{j}^{\epsilon ',\epsilon }(\omega )= \sum_{I}
  \zeta_{I}\land \\ 
\sum_{J\,\mid\, j\not \in J} 
  \frac{1}{2\pi \sqrt{-1}}
  \int_{\partial\Delta _{\epsilon '}}  
  f_{I,J}(\dots ,z_{j-1},w,z_{j+1},\dots )\frac{\dd w} 
  {w-z_{j}}
  \dd \bar z_{J}.
\end{multline*}
To ease notation, if $\epsilon $ and $\epsilon '$ are clear from
the context, we will drop them and write $K_{j}$ and $P_{j}$ instead
of $K^{\epsilon ',\epsilon }_{j}$ and $P^{\epsilon ',\epsilon }_{j}$.
\end{definition}

\begin{remark} \label{rem:6}
  The reason why we use the differentials $\zeta _{I}$ instead of $\dd
  z_{I}$ in the definition of $K$ and $P$, is that, in general, a log form is
  not locally integrable. For instance, if $d=k=1$ and $\omega =
  f\,\dd z\land \dd \bar z$ is a section of $\mathscr{E}^{1,1}_{\Delta
  _{\epsilon }}\langle 0 \rangle(\Delta
  _{\epsilon })$, then $f$ satisfies
  \begin{displaymath}
    |f(z)|\prec \frac{|\log (1/r)|^{N}}{r^{2}} 
  \end{displaymath}
  and the integral
  \begin{displaymath}
    \int_{\overline {\Delta} _{\epsilon '}}  
  \frac{|\log (1/|w|)|^{N}}{|w|^{2}} \frac{\dd w\land \dd\bar w} 
  {w-z}
  \end{displaymath}
  does not converge. But, by the definition we have adopted,
  $K^{\epsilon ',\epsilon 
  }(\omega )=g\,\dd z$, where
  \begin{displaymath}
    g(z)=\frac{1}{z} I_{\epsilon '}(z.f)= \frac{1}{2\pi \sqrt{-1}}
  \frac{1}{z}  \int_{\overline {\Delta} _{\epsilon '}}  
  w f(w)\frac{\dd w\land \dd\bar w} 
  {w-z}.  
  \end{displaymath}
  This integral is absolutely convergent and
  \begin{displaymath}
    \frac{\partial }{\partial\bar
    z}g(z)=\frac{1}{z}\frac{\partial}{\partial\bar z} I_{\epsilon
    '}(z.f)(z)=\frac{zf(z)}{z}=f(z).     
  \end{displaymath}
  This trick will force us to be careful when studying the compatibility of
  $K$ with the operator $\partial$ because, for a log form $\omega $,
  the definitions of $K(\omega )$ and of $K(\partial \omega )$ use
  different kernels in the integral operators. 
\end{remark}

\begin{theorem}
Let $\omega \in \mathscr{E}^{p,q }_{X}\langle
\DD_1 \langle \DD_2\rangle\rangle_{\growth}(\Delta _{x,\epsilon })$. Then
\begin{align*}
 K^{\epsilon ',\epsilon }_{j}(\omega) &\in 
\mathscr{E}^{p,q -1}_{X}\langle \DD_1 \langle \DD_2\rangle\rangle_{\growth}(\Delta
_{x,\epsilon '}),\ \text{and} \\ 
P_{j}^{\epsilon ',\epsilon }(\omega) &\in \mathscr{E}^{p,q }_{X}\langle
\DD_1 \langle \DD_2\rangle\rangle_{\growth}(\Delta _{x,\epsilon '}).
\end{align*}
These operators 
satisfy 
\begin{enumerate}
\item 
 If the form $\omega $ does not contain any term
  with $\dd\bar z_{i}$ for $i>j$ then $K_{j}\omega $ and
  $P_{j}\omega $ do not contain any term with $\dd\bar z_{i}$,
  for $i\ge j$.
\item 
If $\omega \in \mathscr{E}^{p,q }_{X}\langle
\DD_1 \langle \DD_2\rangle\rangle(\Delta _{x,\epsilon })$, then
\begin{align*}
 K^{\epsilon ',\epsilon }_{j}(\omega) &\in 
\mathscr{E}^{p,q-1 }_{X}\langle \DD_1 \langle \DD_2\rangle\rangle(\Delta
_{x,\epsilon '}),\text{ and }\\
P_{j}^{\epsilon ',\epsilon }(\omega) &\in \mathscr{E}^{p,q }_{X}\langle
\DD_1 \langle \DD_2\rangle\rangle(\Delta _{x,\epsilon '}).
\end{align*}
\item 
In this case, $\bar \partial K_{j}+K_{j}\bar
  \partial +P_{j}=\Id$.
\end{enumerate}
\end{theorem}
\begin{proof}
  By lemma \ref{lemm:growth} and the theorem of derivation under the
  integral sign we have that $K_{j}(\omega)\in \mathscr{E}^{p,q-1
  }_{X}\langle \DD_1 \langle \DD_2 \rangle\rangle_{\growth}(\Delta
  _{x,\epsilon '})$ and it is clear that $P_{j}(\omega )\in
  \mathscr{E}^{p,q}_{X}\langle \DD_1\langle \DD_2
  \rangle\rangle_{\growth}(\Delta _{x,\epsilon '})$. Then property
  (i) 
follows from the definition and it is easy to see
  that, if $\bar \partial\omega $, $\partial \omega $ and $\partial
  \bar \partial \omega $ belong to $\mathscr{E}^{\ast}_{X}\langle \DD_1
  \langle \DD_2\rangle\rangle_{\growth}(\Delta _{x,\epsilon })$, the
  same is true for $\bar \partial P_{j}(\omega )$, $\partial
  P_{j}(\omega )$ and $\partial \bar \partial P_{j}(\omega )$.

  In the sequel of the proof, we will denote by
$E_{m}$ the divisor $z_{m}=0$. 
Assume now that $\bar \partial \omega \in \mathscr{E}^{\ast}_{X}\langle
\DD_1 \langle \DD_2\rangle\rangle_{\growth}(\Delta _{x,\epsilon })$. We
will prove 
property (iii). 
We write
\begin{align*}
  \omega &=\sum_{I,J}f_{I,J} \zeta _{I}\land \dd \bar z_{J}\\
  \omega_{1} &=\sum_{I,j\in J}f_{I,J} \zeta _{I}\land \dd \bar z_{J}\\
  \omega _{2}&=\sum_{I,j\not \in J}f_{I,J} \zeta _{I}\land \dd \bar z_{J}.\\
\end{align*}
Recall that we have introduced the operator $\bar \partial_{j}$ in
equation \eqref{eq:9}. We write $\bar \partial_{\not = j}=\bar \partial - \bar
\partial_{j}$ and we decompose
\begin{displaymath}
  \bar \partial K_{j}(\omega )=\bar \partial K_{j}(\omega _{1})=
  \bar \partial _{\not = j} K_{j}(\omega _{1}) + \bar \partial_{j}
  K_{j}(\omega _{1}), 
\end{displaymath}
and 
\begin{align*}
  K_{j}(\bar \partial \omega )=
  K_{j}(\bar \partial_{\not =j }\omega _{1}+\bar \partial_{j} \omega _{2})
\end{align*}
The difficulty at this point is that, when $k<j\le l$, the form
$\omega $ is log-log along $E_{j}$ but, according to lemma
\ref{lemm:13}, $\bar \partial_{j}\omega $ only needs lo be log along
$E_{j}$, and the integral operator $K_{j}$ for log-log forms
may diverge when applied to log forms.   
The key point is to observe that the extra hypothesis about $\bar
\partial \omega $ allow us to apply the operator $K_{j}$
to the differential forms $\bar
\partial_{\not =j }\omega 
_{1}$ and $\bar \partial_{j} \omega _{2}$ individually: Fix $I$ and
$J$ with $j \in  
J$ and $m\not =j$. We consider first the problematic case $k<j\le l$. 
By lemma \ref{lemm:13},
\begin{displaymath}
  \bar \partial_{m} \omega _{I,J_{m}}\in 
  \begin{cases}
    \mathscr{E}^{\ast}_{X}\langle
    \DD_{1}\cup E_{m}\langle \DD_{2}\rangle\rangle_{\growth}(\Delta
    ^{d}_{x,\epsilon }), 
    &\text{ if } k<m\le l,\\ 
    \mathscr{E}^{\ast}_{X}\langle
    \DD_{1}\langle \DD_{2}\rangle\rangle_{\growth}(\Delta
    ^{d}_{x,\epsilon }), 
    &\text{ otherwise.} 
  \end{cases}
\end{displaymath}
Therefore, if we denote by $\DD'$ the union of all the components of
$\DD$ different from $E_{j}$, then
\begin{displaymath}
  (\bar \partial_{\not = j}\omega _{1})_{I,J}\in \mathscr{E}^{\ast}_{X}\langle
  \DD'\langle E_{j}\rangle\rangle_{\growth}(\Delta
  ^{d}_{x,\epsilon }). 
\end{displaymath}
Since, by hypothesis, $(\bar \partial
\omega)_{I,J} \in \mathscr{E}^{\ast}_{X}\langle 
\DD_1 \langle \DD_2\rangle\rangle_{\growth}(\Delta _{x,\epsilon })$ and
$(\bar \partial_{j} \omega _{2})_{I,J}=(\bar \partial \omega -\bar
\partial_{\not = j}\omega _{1})_{I,J}$, 
 then
\begin{displaymath}
  (\bar \partial_{j} \omega _{2})_{I,J}\in \mathscr{E}^{\ast}_{X}\langle
  \DD'\langle E_{j}\rangle\rangle_{\growth}(\Delta
  ^{d}_{x,\epsilon }),
\end{displaymath}
and we can apply the operator $K_{j}$ for log-log forms to the
differential forms $\bar 
\partial_{\not =j }\omega 
_{1}$ and $\bar \partial_{j} \omega _{2}$ individually.
If $j\le k$ then $\omega $ is log along $E_{j}$ and the same is true for the
differential forms $\bar 
\partial_{\not =j }\omega 
_{1}$ and $\bar \partial_{j} \omega _{2}$. But in this case the
operator $K_{j}$ is the operator for log forms and can be applied to 
$\bar 
\partial_{\not =j }\omega 
_{1}$ and $\bar \partial_{j} \omega _{2}$ individually.
The case $j>l$ is similar.
Thus we can write
\begin{displaymath}
  K_{j}(\bar \partial_{\not =j }\omega _{1}+\bar \partial_{j} \omega
  _{2})=
  K_{j}(\bar \partial_{\not =j }\omega _{1})+K_{j}(\bar \partial_{j}
  \omega _{2}). 
\end{displaymath}
But by the theorem of derivation under the integral sign
\begin{displaymath}
  \bar \partial_{\not =j }K_{j}(\omega _{1})+
  K_{j}(\bar \partial_{\not =j }\omega _{1})=0.
\end{displaymath}
By lemma \ref{lemm:growth} 
\begin{displaymath}
  \bar \partial_{j} K_{j}(\omega _{1})=\omega _{1},
\end{displaymath}
and by the generalized Cauchy integral formula (lemma
\ref{lemm:14} (ii)) 
\begin{displaymath}
  K_{j}(\bar \partial_{j} \omega _{2})=\omega _{2}-P_{j}(\omega _{2})=
  \omega _{2}-P_{j}(\omega).
\end{displaymath}
Summing up we obtain
\begin{equation} \label{eq:8}
  \bar \partial K_{j}(\omega )+K_{j}(\bar \partial \omega )=\omega -
  P_{j}(\omega ). 
\end{equation}
By \eqref{eq:8} and the fact that $K_{j}(\bar \partial \omega
),P_{j}(\omega )\in \mathscr{E}^{\ast}_{X}\langle \DD_1 \langle 
\DD_2\rangle\rangle_{\growth}(\Delta _{x,\epsilon '})$, we obtain that
$$\bar \partial 
K_{j}(\omega )\in \mathscr{E}^{\ast}_{X}\langle \DD_1 \langle
\DD_2\rangle\rangle_{\growth}(\Delta _{x,\epsilon '}).$$

Assume now that $\partial \omega \in \mathscr{E}^{\ast}_{X}\langle
\DD_1 \langle \DD_2\rangle\rangle_{\growth}(\Delta _{x,\epsilon })$. 
We fix $I,J\subset \{1,\dots ,d\}$, with $j\in J$. If $j\not \in I$ then
\begin{displaymath}
  (\partial K_{j}(\omega ))_{I,J_{j}}=
  \sum_{m\not= j} \partial_{m} K_{j}(\omega _{I_{m},J})=
  K_{j}\left(\sum_{m\not= j} \partial_{m} \omega _{I_{m},J}\right)=
  K_{j}((\partial \omega)_{I,J} ).
\end{displaymath}
Therefore it belongs to $\mathscr{E}^{\ast}_{X}\langle 
  \DD_{1}\langle \DD_{2}\rangle\rangle_{\growth}(\Delta ^{d}_{x,\epsilon })$.
If $j\in I$,
we write
\begin{equation}\label{eq:10}
  (\partial K_{j}(\omega ))_{I,J_{j}}=
  \sum_{m\not= j} \partial_{m} K_{j}(\omega
  _{I_{m},J})+\partial_{j}K_{j}(\omega 
  _{I_{j},I}).  
\end{equation}

The theorem of derivation under the integral sign implies that, for $m\not=j$
\begin{displaymath}
  \partial_{m}K_{j}(\omega _{I_{m},J})=-K_{j}(\partial_{m}\omega _{I_{m},J}),
\end{displaymath}
Note that the right term is well defined by lemma \ref{lemm:13}.  We
treat first the case $j\le k$. We have to be careful because
the integral kernel that appear in the expressions
$\partial_{j}K_{j}(\omega _{I_{j},J})$ and $K_{j}(\partial_{j} \omega
_{I_{j},J})$ is different in each term.

Again by lemma \ref{lemm:13}, 
\begin{displaymath}
  \partial_{j} \omega _{I_{j},J}\in \mathscr{E}^{\ast}_{X}\langle
  \DD_{1}\langle \DD_{2}\rangle\rangle_{\growth}(\Delta ^{d}_{x,\epsilon }.
\end{displaymath}
Since moreover $\partial \omega \in \mathscr{E}^{\ast}_{X}\langle
  \DD_{1}\langle \DD_{2}\rangle\rangle_{\growth}(\Delta
  ^{d}_{x,\epsilon })$, 
\begin{displaymath}
  \sum_{m\not = j} \partial_{m} \omega _{I_{m},J}=(\partial \omega
  )_{I,J}- \partial_{j} \omega _{I_{j},J}
\in \mathscr{E}^{\ast}_{X}\langle
  \DD_{1}\langle \DD_{2}\rangle\rangle_{\growth}(\Delta ^{d}_{x,\epsilon }).
\end{displaymath}
Hence, by lemma \ref{lemm:growth}
\begin{displaymath}
  K_{j}\left(\sum_{m\not = j} \partial_{m} \omega _{I_{m},J}\right)\in
  \mathscr{E}^{\ast}_{X}\langle 
  \DD_{1}\langle \DD_{2}\rangle\rangle_{\growth}(\Delta ^{d}_{x,\epsilon }).
\end{displaymath}
By the same lemma it follows that
\begin{displaymath}
  \partial_{j}K_{j}(\omega
  _{I_{j},I})\in
  \mathscr{E}^{\ast}_{X}\langle 
  \DD_{1}\langle \DD_{2}\rangle\rangle_{\growth}(\Delta ^{d}_{x,\epsilon }).
\end{displaymath}
Now we treat the case $j>k$. In this case the expressions
$\partial_{j}K_{j}(\omega _{I_{j},J})$ and $K_{j}(\partial_{j} \omega
_{I_{j},J})$
use the same integral kernel. By lemma
\ref{lemm:14} (iii)
\begin{displaymath}
  \partial_{j}K_{j}(\omega _{I_{j},J})=-K_{j}(\partial_{j} \omega
  _{I_{j},J})
  +\frac{(-1)^{\# I+\sigma (J,j)+\sigma (I,j)}}{2\pi \sqrt{-1}} 
  \int_{\gamma _{\epsilon '}}f_{I_{j},J}\frac{\dd \bar w}{w-z} \zeta
  _{I}\land \dd \bar z_{J_{j}}.
\end{displaymath}
Hence 
\begin{displaymath}
  (\partial K_{j}(\omega ))_{I,J_{j}}=-(K_{j}(\partial \omega
  ))_{I,J_{j}}
    +\frac{(-1)^{\# I+\sigma (J,j)+\sigma (I,j)}}{2\pi \sqrt{-1}} 
    \int_{\gamma _{\epsilon '}}f_{I_{j},J}\frac{\dd \bar w}{w-z} \zeta
  _{I}\land \dd \bar z_{J_{j}}.
\end{displaymath}
Thus it belongs to $\mathscr{E}^{\ast}_{X}\langle 
  \DD_{1}\langle \DD_{2}\rangle\rangle_{\growth}(\Delta
  ^{d}_{x,\epsilon })$.

  Finally assume that $\partial
  \omega, \bar \partial \omega, \partial \bar \partial \omega   
  \in \mathscr{E}^{\ast}_{X}\langle 
  \DD_{1}\langle \DD_{2}\rangle\rangle_{\growth}(\Delta
  ^{d}_{x,\epsilon })$. By equation \eqref{eq:8}
  \begin{displaymath}
    \partial \bar \partial K_{j}(\omega )=
    -\partial K_{j}(\bar \partial \omega )+
    \partial \omega -\partial P_{j}(\omega )
  \end{displaymath}
  therefore the result follows from the previous cases.
\end{proof}

\subsection{Good forms.}
\label{sec:good-forms-poincare}

In this section we recall the definition of good forms in the
sense of \cite{Mumford:Hptncc}. We introduce the complex of Poincar\'e
singular forms  
 that is contained in both, the
complex of good forms and the complex of log-log forms.

\nnpar{Poincar\' e growth.} Let $X$, $\DD$, $U$ and $\iota $ be as
in definition \ref{def:log}.

\begin{definition}
  \label{def:poincaregrowth} Let $V$ be a coordinate neighborhood adapted
  to $\DD$.
  We say that a smooth complex function $f$ on $V\setminus \DD$ 
  has \emph{Poincar\'e growth (along  $\DD$)} if it is bounded. We say
  that it
  has \emph{Poincar\'e growth (along $\DD$) of infinite order}, if
  for all multi-indices $\alpha, \beta  \in \mathbb{Z}_{\ge 0}^{d}$,
  \begin{equation}
    \left |\frac{\partial^{|\alpha |}}{\partial z^{\alpha }}
      \frac{\partial^{|\beta |}}{\partial \bar z^{\beta }}
      f(z_{1},\dots ,z_{d})\right | \prec \frac{ 1}{|z^{\alpha^{\le
      k}}\bar z^{\beta ^{\le k}}|}. 
  \end{equation}
  The \emph{sheaf of differential forms on $X$ with
    Poincar\'e growth (resp. of infinite order)} is the
    subalgebra of 
  $\iota _{\ast}\mathscr 
  {E}^{\ast}_{U}$ generated, in each coordinate neighborhood $V$
  adapted to $\DD$, by the functions with Poincar\'e growth
    (resp. of infinite order) and the differentials
  \begin{alignat*}{2}
    &\frac{\dd z_{i}}{z_{i}\log(1/r_{i})},\,\frac{\dd\bar{z}_{i}}
    {\bar{z}_{i}\log(1/r_{i})},&\qquad\text{for }i&=1,\dots,k, \\
    &\dd z_{i},\,\dd\bar{z}_{i},&\qquad\text{for }i&=k+1,\dots,d.
  \end{alignat*}
\end{definition}

\nnpar{Good forms.}   We recall that a smooth form $\omega $ on
  $X\setminus \DD$ is 
  \emph{good  (along  $\DD$)} if $\omega $ and
  $\dd \omega $ 
  have Poincar\'e 
  growth along $\DD$ (\cite{Mumford:Hptncc}). Observe that, since the
  operator $\dd$ is not 
  bi-homogeneous, the sheaf of good forms is not bigraded. Although
  good forms are very similar to pre-log-log forms, there is no
  inclusion between both sheaves. Nevertheless we have the following
  easy 
  \begin{lemma}
    If $\omega $ is a good form of pure bidegree, then it is a
    pre-log-log form if and only if $\partial \bar \partial \omega $
    has log-log growth of order 0. \hfill $\square$
  \end{lemma}

  \nnpar{Poincar\'e singular forms.} 
  \begin{definition}\label{def:15}
    We will say that $\omega $ is \emph{Poincar\'e singular (along
    $\DD$)}, if 
    $\omega $,  
    $\partial\omega $, $\bar \partial \omega $ and $\partial \bar
    \partial \omega $ have Poincar\'e growth of infinite
    order.  
  \end{definition}

Note that the sheaf of Poincar\'e singular forms is contained in both,
the sheaf of 
good forms and the sheaf of log-log
forms. Observe moreover that we cannot expect to have a Poincar\'e
lemma for the complex of Poincar\'e singular forms, due precisely to
the absence of 
the functions $\log (\log (1/r_{i}))$.

\nnpar{Functoriality.} The complex of Poincar\'e singular forms share
some of the 
properties of the complex of log-log forms. For instance we have the
following compatibility with respect inverse images that is proved as
in proposition \ref{prop:10}.

\begin{proposition} \label{prop:invimggood}
Let $f:X\longrightarrow Y$ be a morphism of complex
manifolds of dimension $d$ and $d'$, let $\DD_{X}$, $\DD_{Y}$ be normal
crossing divisors on $X$ and 
$Y$ respectively, satisfying $f^{-1}(\DD_{Y})\subseteq \DD_{X}$.  If
$\eta$ is a Poincar\'e singular form on $Y$, then $f^{\ast}\eta$ is a
Poincar\'e singular form on
$X$.  \hfill$\square$
\end{proposition}

\section{Arithmetic Chow rings with log-log growth conditions}
\label{sec:arithmetic-chow-ring-log-log}

In this section we use the theory of abstract cohomological arithmetic Chow
rings developed in \cite{BurgosKramerKuehn:cacg} to obtain a theory of
arithmetic Chow rings with log-log  forms.
Since we have computed the
cohomology of the complex of log-log forms, we have a more
precise knowledge of the size of these arithmetic Chow rings than of
 the arithmetic Chow
rings with pre-log-log  forms considered in 
\cite{BurgosKramerKuehn:cacg}.

\subsection{Dolbeault algebras and Deligne algebras }
\label{sec:dolb-algebr-deligne}

In this section we recall the notion of Dolbeault algebra and the
properties of the associated Deligne algebra.

\nnpar{Dolbeault algebras.}
\begin{definition}\label{def:5}
  A \emph{Dolbeault algebra} $A=(A_{\mathbb{R}}^{\ast},\dd_{A},\land)$
  is a real differential graded commutative algebra which is bounded
  from below and equipped with a bigrading on $A_{\mathbb{C}}:=
  A_{\mathbb{R}}\otimes{\mathbb{C}}$,
\begin{displaymath}
A_{\mathbb{C}}^{n}=\bigoplus_{p+q=n}A^{p,q},
\end{displaymath}  
satisfying the following properties:
\begin{enumerate}
\item[(i)]
The differential $\dd_{A}$ can be decomposed as the sum of operators
$\dd_{A}=\partial+\bar{\partial}$ of type $(1,0)$, resp. $(0,1)$.
\item[(ii)] 
It satisfies the symmetry property $\overline{A^{p,q}}=A^{q,p}$,
where $\overline{\phantom{M}}$ denotes complex conjugation.
\item [(iii)]  The product induced on $A_{\mathbb{C}}$ is 
compatible with the bigrading:
\begin{displaymath}  
A^{p,q}\land A^{p',q'}\subseteq A^{p+p',q+q'}.
\end{displaymath}
\end{enumerate}
By abuse of notation, we will denote also by $A^{\ast}$
the complex differential graded commutative algebra
$A^{\ast}_{\mathbb{C}}$.
\end{definition}

\begin{notation}
\label{def:13}
Given a Dolbeault algebra $A$ we will use the following notations.
The Hodge filtration $F$ of $A^{\ast}$ is the decreasing filtration
given by
\begin{displaymath}
F^{p}A^{n}_{\mathbb{C}}=\bigoplus_{p'\geq p}A^{p',n-p'}.
\end{displaymath}
The filtration $\overline F$ is the complex conjugate of $F$, i.e.,
\begin{displaymath}
\overline F^{p}A^{n}=\overline{F^{p}A^{n}}.
\end{displaymath}
For an element $x\in A$, we write $x^{i,j}$ for its 
component in $A^{i,j}$. For $k,k'\geq 0$, we define an 
operator $F^{k,k'}: A\longrightarrow A$
by the rule 
\begin{displaymath}
F^{k,k'}(x):=\sum_{l\geq k,l'\geq k'}x^{l,l'}.
\end{displaymath}
We note that the operator $F^{k,k'}$ is the projection of $A^{\ast}$
onto the subspace $F^{k}A^{\ast}\cap 
\overline{F}^{k'}A^{\ast}$. We will write $F^{k}=F^{k,
-\infty}$. 

We denote by $A^n_{\mathbb{R}}(p)$ the subgroup $(2\pi i)^{p}\cdot
A^{n}_{\mathbb{R}}\subseteq  
A^{n}$, and we define the operator
$$ 
\pi_{p}:A\longrightarrow A_{\mathbb{R}}(p)
$$
by setting $\pi_{p}(x):=\frac{1}{2}(x+(-1)^{p}\bar{x})$.
\end{notation}

\nnpar{The Deligne complex.} 

\begin{definition}
Let $A$ be a Dolbeault algebra. Then, the \emph{Deligne complex 
$(\mathcal{D}^{\ast}(A,\ast),\dd_{\mathcal{D}})$ associated to $A$} 
is the graded complex given by   
\begin{align*}
&\mathcal{D}^{n}(A,p)=
\begin{cases}
A_{\mathbb{R}}^{n-1}(p-1)\cap F^{n-p,n-p}A^{n-1},
&\qquad\text{if}\quad n\leq 2p-1, \\
A_{\mathbb{R}}^n(p)\cap F^{p,p}A^{n},
&\qquad\text{if}\quad n\geq 2p,  
\end{cases}  
\intertext{with differential given by ($x\in\mathcal{D}^{n}(A,p)$)}
&\dd_{\mathcal{D}}x=
\begin{cases}
-F^{n-p+1,n-p+1}\dd_{A}x,
&\qquad\text{if}\quad n<2p-1, \\
-2\partial\bar{\partial}x, 
&\qquad\text{if}\quad n=2p-1, \\
\dd_{A}x, 
&\qquad\text{if}\quad n\geq 2p.
\end{cases}  
\end{align*}  
\end{definition}

\nnpar{The Deligne algebra.}

\begin{definition}
\label{def:17}
Let $A$ be a Dolbeault algebra. The \emph{Deligne algebra associated 
to $A$} is the Deligne complex $\mathcal{D}^{\ast}(A,\ast)$ together 
with the graded commutative product
$\bullet:\mathcal{D}^{n}(A,p)\times\mathcal{D}^{m} 
(A,q)\longrightarrow\mathcal{D}^{n+m}(A,p+q)$ given by   
\begin{align*}
&x\bullet y= \\
&\begin{cases}
(-1)^{n}r_{p}(x)\land y+x\land r_{q}(y),
&\text{if }n<2p,\,m<2q, \\
F^{l-r,l-r}(x\land y),
&\text{if }n<2p,\,m\geq 2q,\,l<2r, \\
F^{r,r}(r_{p}(x)\land y)+2\pi_{r}(\partial(x\land y)^{r-1,l-r}),
&\text{if }n<2p,\,m\geq 2q,\,l\geq 2r, \\
x\land y,
&\text{if }n\geq 2p,\,m\geq 2q,
\end{cases}
\end{align*}
where we have written $l=n+m$, $r=p+q$, and $r_{p}(x)=2\pi_{p}
(F^{p}\dd_{A}x)$.  
\end{definition}

\nnpar{Specific degrees.}
In the sequel we will be interested in some specific degrees
where we can give simpler formulas. Namely, we consider
\begin{align*}
&\mathcal{D}^{2p}(A,p)=A^{2p}_{\mathbb{R}}(p)\cap A^{p,p}, \\
&\mathcal{D}^{2p-1}(A,p)=A^{2p-2}_{\mathbb{R}}(p-1)\cap A^{p-1,p-1}, \\
&\mathcal{D}^{2p-2}(A,p)=A^{2p-3}_{\mathbb{R}}(p-1)\cap
(A^{p-2,p-1}\oplus A^{p-1,p-2}).
\end{align*}
The corresponding differentials are given by
\begin{alignat*}{2}
&\dd_{\mathcal{D}}x=\dd_{A}x, 
&\quad &\text{if }x\in\mathcal{D}^{2p}(A,p), \\
&\dd_{\mathcal{D}}x=-2\partial\bar{\partial}x, 
&\quad &\text{if }x\in\mathcal{D}^{2p-1}(A,p), \\
&\dd_{\mathcal{D}}(x,y)=-\partial x-\bar{\partial}y, 
&\quad &\text{if }(x,y)\in\mathcal{D}^{2p-2}(A,p). 
\end{alignat*}
Moreover, the product is given as follows: for $x\in\mathcal{D}^
{2p}(A,p)$, $y\in\mathcal{D}^{2q}(A,q)$ or $y\in\mathcal{D}^{2q-1}
(A,q)$, we have
\begin{align*}
x\bullet y&=x\land y, \\
\intertext{and for $x\in\mathcal{D}^{2p-1}(A,p)$, $y\in
\mathcal{D}^{2q-1}(A,q)$, we have}
x\bullet y&=-\partial x\land y+\bar{\partial}x\land y+
x\land\partial y-x\land\bar{\partial}y.
\end{align*}

\nnpar{Deligne complexes and Deligne-Beilinson cohomology.}
The main interest in Deligne complexes is expressed by the 
following theorem which is proven in \cite{Burgos:CDB} in a 
particular case, although the proof is valid in general.

\begin{theorem} \label{thm:4} 
Let $X$ be a complex algebraic manifold, $\overline{X}$ a 
smooth compactification of $X$ with $D=\overline{X}\setminus 
X$ a normal crossing divisor, and denote by $j:X\longrightarrow
\overline{X}$ the natural inclusion. Let $\mathscr{A}^{\ast}$ 
be a sheaf of Dolbeault algebras over $\overline{X}^
{\an}$ such that, for every $n,p$ the sheaves $\mathscr{A}^{\ast}$ and
$F^{p}\mathscr{A}^{\ast}$ are acyclic,
$\mathscr{A}_{\mathbb{R}}^{\ast}$ is a multiplicative  
resolution of $Rj_{\ast}\mathbb{R}$ and $(\mathscr{A}^{\ast},F)$ is a
multiplicative filtered resolution of 
$(\Omega^{\ast}_{\overline{X}}(\log D),F)$. Putting $A^{\ast}=
\Gamma(\overline{X},\mathscr{A}^{\ast})$, we have a natural 
isomorphism of graded algebras  
\begin{displaymath}
H^{\ast}_{\mathcal{D}}(X,\mathbb{R}(p))\cong H^{\ast}
(\mathcal{D}(A,p)).
\end{displaymath}
\hfill $\square$
\end{theorem}

\nnpar{Notation.} In the sequel we will use the following notation.
The sheaves of differential forms will be denoted the italic letter
$\mathscr{E}$, and the 
corresponding spaces of global sections will be denoted with the same
letter in roman typography $E$. For instance
\begin{displaymath}
  E_X^{n}\langle D_1\langle D_2 \rangle \rangle =
  \Gamma( X,\mathscr{E}_X^{n}\langle D_1\langle D_2 \rangle \rangle). 
\end{displaymath}

\nnpar{The Deligne complex with logarithmic singularities.}
Let $X$ be a quasi-projective complex manifold. Let $E_{\log}(X)$ be
the Dolbeault algebra of differential forms with logarithmic
singularities at infinity (see \cite{BurgosKramerKuehn:cacg}  \S5). We
will denote
\begin{displaymath}
  \mathcal{D}^{\ast}_{\log}(X,p)=\mathcal{D}^{\ast}(E_{\log}(X),p).
\end{displaymath}
Then theorem \ref{thm:4} implies that 
\begin{displaymath}
  H^{\ast}_{\mathcal{D}}(X,\mathbb{R}(p))\cong H^{\ast}
  (\mathcal{D}_{\log}(X,p)).
\end{displaymath}

\subsection{The $\mathcal {D}_{\log}$-complex of log-log forms}
\label{sec:mathc-d_log-compl}

\nnpar{$\mathcal{D}_{\log}$-complexes.}
Recall that, to define the arithmetic Chow groups of an arithmetic
variety $X$ as in
\cite{BurgosKramerKuehn:cacg} we need first an auxiliary complex of
graded abelian sheaves on the Zariski site of smooth real
schemes, that
satisfies Gillet axioms. As in \cite{BurgosKramerKuehn:cacg} we will
use the complex of sheaves
$\mathcal{D}_{\log}$. This sheaf is given, for any smooth real scheme
$U_{\mathbb{R}}$, by
\begin{displaymath}
  \mathcal{D}_{\log}(U_{\mathbb{R}},p)=
  \mathcal{D}_{\log}(U_{\mathbb{C}},p)^{\sigma  }, 
\end{displaymath}
where $\sigma $ is the involution that acts as complex conjugation on
the space and on the coefficients.
(see \cite{BurgosKramerKuehn:cacg} \S 5.3). 

Then, we need to choose a
$\mathcal{D}_{\log}$-complex over $X_{\mathbb{R}}$. That is, a
complex, $\cc^{\ast}_{X_{\mathbb{R}}}(\ast)$, of graded abelian
sheaves on the 
Zariski topology of $X_{\mathbb{R}}$ together with a morphism  
\begin{displaymath}
  \mathcal{D}_{\log,X_{\mathbb{R}}}
  \longrightarrow \cc_{X_{\mathbb{R}}},
\end{displaymath}
such that all the sheaves $\cc^{n}_{X_{\mathbb{R}}}(p)$ are totally
acyclic (see 
\cite{BurgosKramerKuehn:cacg} definition 3.1 and definition 3.4). The
$\mathcal{D}_{\log}$-complex $\cc$ plays the role of the fiber over the
Archimedean places of the arithmetic ring $A$. The aim of this section
is to construct a $\mathcal{D}_{\log}$-complex by mixing log and
log-log forms.  

\nnpar{Varieties with a fixed normal crossing divisor.} We will follow
the notations of \cite{BurgosKramerKuehn:cacg} \S 7.4 that we recall
shortly.  
Let $X$ be a complex algebraic manifold of dimension $d$, and $D$  
a normal crossing divisor. We will denote by $\underline{X}$ the 
pair $(X,D)$. If $W\subseteq X$ is an open subset, we will write 
$\underline{W}=(W,D\cap W)$.

In the sequel we will consider all operations adapted to the pair 
$\underline{X}$. For instance, if $Y\subsetneq X$ is a closed 
algebraic subset and $W=X\setminus Y$, then an embedded resolution 
of singularities of $Y$ in $\underline{X}$ is a proper modification 
$\pi:\widetilde{X}\longrightarrow X$ such that $\pi\big|_{\pi^{-1}
(W)}:\pi^{-1}(W)\longrightarrow W$ is an isomorphism, and
\begin{displaymath}
\pi^{-1}(Y),\,\pi^{-1}(D),\,\pi^{-1}(Y\cup D)  
\end{displaymath}
are normal crossing divisors on $\widetilde{X}$. Using Hironaka's 
theorem on the resolution of singularities \cite{Hironaka:rs}, one 
can see that such an embedded resolution of singularities exists.

Analogously, a normal crossing compactification of $\underline{X}$
will be a smooth compactification $\overline{X}$ such that the 
adherence $\overline{D}$ of $D$, the subsets $B_{\overline{X}}=
\overline{X}\setminus X$ and  $B_{\overline{X}}\cup
\overline{D}$  are normal 
crossing divisors. 

\nnpar{Logarithmic growth along infinity.}
Given a diagram of normal crossing compactifications of $\underline{X}$ 
\begin{displaymath}
\xymatrix{
\overline{X}'\ar[r]^{\varphi}&\overline{X} \\
&X\ar[ul]\ar[u],}
\end{displaymath}
with divisors $B_{\overline{X}'}$ and $B_{\overline{X}}$ at 
infinity, respectively, proposition \ref{prop:invimagell}
gives rise to an induced morphism
\begin{displaymath}
\varphi^{\ast }:\mathscr{E}^{\ast}_{\overline{X}}\langle 
B_{\overline{X}}\langle\overline{D}\rangle\rangle
\longrightarrow\mathscr{E}^{\ast}_{\overline{X}'}\langle 
B_{\overline{X}'}\langle\overline{D}'\rangle\rangle.
\end{displaymath}
In order to have a complex that is independent of the choice 
of a particular compactification we take the limit over all possible 
compactifications.

\begin{definition} 
Let $\underline{X}=(X,D)$ be as above. Then we define the
\emph{complex $E^{\ast}_{\lgi}(\underline{X})$ of differential 
forms on $X$ log along infinity and log-log
along $D$} as
\begin{displaymath}
E^{\ast}_{\lgi}(\underline{X})=\lim_{\longrightarrow}\Gamma
(\overline{X},\mathscr{E}^{\ast}_{\overline{X}}\langle 
B_{\overline{X}}\langle\overline{D}\rangle\rangle),
\end{displaymath}
where the limit is taken over all normal crossing compactifications 
$\overline{X}$ of $\underline{X}$.
\end{definition}

\nnpar{A $\mathcal{D}_{\log}$-complex.}
Let $X$ be a smooth real variety and $D$ a normal crossing divisor 
defined over $\mathbb{R}$; as before, we write $\underline{X}=(X,D)$. 
For any $U\subseteq X$, the complex
$E^{\ast}_{\lgi}(\underline{U}_{\mathbb{C}})$  
is a Dolbeault algebra with respect to the wedge product.

\begin{definition}\label{def:1}
For any Zariski open subset $U\subseteq X$, we put
\begin{displaymath}
\mathcal{D}^{\ast}_{\lgi,\underline{X}}(U,p)=(\mathcal{D}^{\ast}_
{\lgi,\underline{X}}(U,p),\dd_{\mathcal{D}})=(\mathcal{D}^{\ast}
(E_{\lgi}(\underline{U}_{\mathbb{C}}),p)^{\sigma},\dd_{\mathcal{D}}),
\end{displaymath}
where the operator $\mathcal{D}$ is as in definition \ref{def:17} and
$\sigma$ 
is the involution that acts as complex conjugation in 
the space and in the coefficients (see \cite{BurgosKramerKuehn:cacg}
5.55). When the pair $\underline X$ is understood then we write
$\mathcal{D}^{\ast}_{\lgi}$ instead of
$\mathcal{D}^{\ast}_{\lgi,\underline{X}}$. The complex
$\mathcal{D}^{\ast}_{\lgi}$ will be called \emph{the
  $\mathcal{D}_{\log}$-complex of log-log forms} or just \emph{the
  complex of log-log forms}.   
\end{definition}

Then it holds the analogue of \cite{BurgosKramerKuehn:cacg} theorem
7.18.
\begin{theorem} 
\label{thm:1}
The complex $\mathcal{D}_{\lgi,\underline{X}}$ is a $\mathcal{D}_
{\log}$-complex on $X$. Moreover, it is a pseudo-associative and
commutative $\mathcal{D}_{\log}$-algebra.
\hfill $\square$
\end{theorem}

\nnpar{The cohomology of the complex $\mathcal{D}_{\lgi,\underline
    X}$.} The main 
advantage of the complex $\mathcal{D}_{\lgi,\underline X}$ over the
complex $\mathcal{D}_{\wlg,\underline X}$ of
\cite{BurgosKramerKuehn:cacg} is the following result, that is a
consequence of theorem \ref{thm:fq} and theorem \ref{thm:4}. (see
    \cite{BurgosKramerKuehn:cacg} 
theorem 5.19 and \cite{Burgos:CDB}). 

\begin{theorem} \label{thm:12}
  The inclusion $\mathcal{D}_{\log,X}\longrightarrow
  \mathcal{D}_{\lgi,\underline X}$ is a quasi-isomorphism. Therefore
  the hyper-cohomology over $X$ of the complex of sheaves
  $\mathcal{D}_{\lgi,\underline X}$, as well as the cohomology of its
  complex of global sections, is
  naturally isomorphic to the Deligne-Beilinson cohomology of $X$.
  \hfill $\square$
\end{theorem}

\subsection{Properties of Green objects with values in $\mathcal{D}_{\lgi}$.}
We start by noting that theorem \ref{thm:1} together with
\cite{BurgosKramerKuehn:cacg} section 3
provides us with a theory of Green objects with values in
$\mathcal{D}_{\lgi,\underline X}$. 

\nnpar{Mixed forms representing the class of a cycle.} Since we know
the cohomology of the complex of mixed forms, we obtain the analogue
of proposition \cite{BurgosKramerKuehn:cacg} 5.48 which is more precise than
the analogue of proposition \cite{BurgosKramerKuehn:cacg} 7.20. In
particular,
\begin{proposition} 
\label{prop:6}
Let $X$ be a smooth real variety and $D$ a normal crossing divisor. Put
$\underline X=(X,D)$. Let $y$ be a $p$-codimensional cycle 
on $X$ with support $Y$. Then, we have that 
the class of the cycle $(\omega,g)$ in $H^{2p}_{\mathcal{D}_
{\lgi},Y}(X,p)$ is equal to the class of $y$, if and only if  
\begin{align}
\label{eq:eqgreenwlog}
-2\partial\bar{\partial}[g]_{X}=[\omega]-\delta_{y}. 
\end{align}
\end{proposition}  
\begin{proof}
The proof is completely analogous to the proof of
\cite{BurgosKramerKuehn:cacg} 5.48, using theorem \ref{thm:12} and
lemma \ref{lemm:3}.
\end{proof}

\nnpar{Inverse images.} 
\begin{proposition}
\label{prop:invimageloglog}
Let $f:X\longrightarrow Y$ be a morphism of smooth real varieties, 
let $D_{X}$, $D_{Y}$ be normal crossing divisors on $X$, $Y$ 
respectively, satisfying $f^{-1}(D_{Y})\subseteq D_{X}$. Put 
$\underline{X}=(X,D_{X})$ and $\underline{Y}=(Y,D_{Y})$. Then,
there exists a contravariant $f$-morphism
\begin{displaymath}
f^{\#}:\mathcal{D}_{\lgi,\underline{Y}}\longrightarrow f_{\ast}  
\mathcal{D}_{\lgi,\underline{X}}.
\end{displaymath}
\end{proposition} 
\begin{proof} 
By proposition \ref{prop:invimagell}, the pull-back of differential
forms induces a morphism of the corresponding Dolbeault algebras of
mixed forms. This morphism is compatible with the involution $\sigma$. 
Thus, this morphism gives rise to an induced morphism between the 
corresponding Deligne algebras. 
\end{proof}

\nnpar{Push-forward.}  
We will only state the most basic property concerning direct images,
which is necessary to define arithmetic degrees. Note however that 
we expect that the complex of log-log forms will be useful in the
study of non smooth,  
proper, surjective morphisms. By proposition \ref{prop:5}, we have 

\begin{proposition}
\label{prop:pushforward-point}
Let $\underline{X}=(X,D)$ be a proper, smooth real variety with 
fixed normal crossing divisor $D$. Let $f:X\longrightarrow\Spec
(\mathbb{R})$ denote the structural morphism. Then, there exists
a covariant $f$-morphism
\begin{displaymath}
f_{\#}:f_{\ast}\mathcal{D}_{\lgi,\underline{X}}\longrightarrow
\mathcal{D}_{\log,\Spec(\mathbb{R})}.
\end{displaymath}
\end{proposition}
\noindent
In particular, if $X$ has dimension $d$, we obtain a well defined 
morphism  
\begin{displaymath}
f_{\#}:\widehat{H}^{2d+2}_{\mathcal{D}_{\lgi},\mathcal{Z}^{d+1}}
(X,d)\longrightarrow\widehat{H}^{2}_{\mathcal{D}_{\log},\mathcal
{Z}^{1}}(\Spec(\mathbb{R}),1)=\mathbb{R}.
\end{displaymath}
Note that, by dimension reasons, we have $\mathcal{Z}^{d+1}=
\emptyset$, and 
\begin{displaymath}
\widehat{H}^{2d+2}_{\mathcal{D}_{\lgi},\mathcal{Z}^{d+1}}(X,d)=
H^{2d+1}(\mathcal{D}_{\lgi}(X,d+1))=
H^{2d+1}_{\mathcal{D}}(X,\mathbb{R}(d+1). 
\end{displaymath}
Thus, every element of $\widehat{H}^{2d+2}_{\mathcal{D}_{\lgi},
\mathcal{Z}^{d+1}}(X,d)$ is represented by a pair $\mathfrak{g}=
(0,\widetilde{g})$. The morphism $f_{\#}$ mentioned above, is 
then given by
\begin{displaymath}
\mathfrak{g}=(0,\widetilde{g})\longmapsto\left(0,\frac{1}{(2\pi i)^
{d}}\int_{X}g\right).
\end{displaymath}

\subsection{Arithmetic Chow rings with log-log forms}
\label{sec:acgwllg}

\nnpar{Arithmetic Chow groups.} We are now in position to apply the
machinery of \cite{BurgosKramerKuehn:cacg}. Let $(A,\Sigma
,F_{\infty})$ be an arithmetic ring and let $X$ be a regular 
arithmetic variety over $A$. Let $\DD$ be a fixed
normal crossings divisor of $X_{\Sigma }$ stable under
$F_{\infty}$. As in the previous section we will denote by
$\underline {X}$ the pair $(X_{\mathbb{R}},\DD)$.  The natural
inclusion $\mathcal{D}_{\log}\longrightarrow \mathcal{D}_{\lgi}$
induces a $\mathcal{D}_{\log}$-complex structure in
$\mathcal{D}_{\lgi}$. Then,
$(X,\mathcal{D}_{\lgi})$ is a
$\mathcal{D}_{\log}$-arithmetic variety. 
Therefore, applying the
theory of \cite{BurgosKramerKuehn:cacg}, section 4, we define the
arithmetic Chow groups 
$\cha^{\ast}(X,\mathcal{D}_{\lgi})$. These groups will
be called log-log arithmetic Chow groups.

\nnpar{Exact sequences.} We start the study of these arithmetic Chow
groups by writing the exact
sequences of \cite{BurgosKramerKuehn:cacg}, theorem 4.13. Observe
that, since we have better control on the cohomology of
$\mathcal{D}_{\lgi}$, we obtain better results than in
\cite{BurgosKramerKuehn:cacg} \S 7.

\begin{theorem}\label{thm:16}
The following sequences are exact:
\begin{align}
&\CH^{p-1,p}(X) \stackrel{\rho}{\longrightarrow}
\widetilde{\mathcal{D}}^{2p-1}_{\lgi}(X,p) 
\stackrel{\amap}{\longrightarrow} 
\cha^p(X,\mathcal{D}_{\lgi}) \stackrel{\zeta}{\longrightarrow}
\CH^p(X) \longrightarrow 0,\notag \\[3mm]
&\CH^{p-1,p}(X)\stackrel{\rho}{\longrightarrow}
H^{2p-1}_{\mathcal{D}}(X_{\mathbb{R}},\mathbb{R}(p)) 
\stackrel{\amap}{\longrightarrow} 
\cha^p(X,\mathcal{D}_{\lgi})\stackrel{(\zeta,-\omega)}{\longrightarrow}  
\notag \\ 
&\phantom{CH^{p-1,p}(X)\stackrel{\rho}{\longrightarrow}}
\CH^p(X) \oplus {\rm Z}\mathcal{D}_{\lgi}^{2p}(X,p) 
\stackrel{\cl +h}{\longrightarrow} 
H^{2p}_{\mathcal{D}}(X_{\mathbb{R}},\mathbb{R}(p)) 
\longrightarrow 0,\notag \\[3mm]
&\CH^{p-1,p}(X) \stackrel{\rho}{\longrightarrow}
H^{2p-1}_{\mathcal{D}}(X_{\mathbb{R}},\mathbb{R}(p))  
\stackrel{\amap}{\longrightarrow} 
\cha^p(X,\mathcal{D}_{\lgi})_0 \stackrel{\zeta}{\longrightarrow}
\CH^p(X)_0 \longrightarrow 0.\notag
\end{align}
\hfill $\square$
\end{theorem} 

\nnpar{Multiplicative properties.} 
Since $\mathcal{D}_{\lgi}$ is a pseudo-associative and commutative
$\mathcal{D}_{\log}$-algebra, we have
\begin{theorem} \label{thm:ll-product}
  The abelian group
  \begin{displaymath}
    \cha^{\ast}(X,\mathcal{D}_{\lgi})_{\mathbb{Q}}=
    \bigoplus_{p\ge 0}\cha^{p}(X,\mathcal{D}_{\lgi})\otimes \mathbb{Q}
  \end{displaymath}
  is an associative commutative $\mathbb{Q}$-algebra with unit.
\hfill $\square$
\end{theorem}

\nnpar{Inverse images.}
 By proposition \ref{prop:invimagell}, there are some cases
 where we can define the inverse image for the log-log arithmetic Chow
 groups. 
 
\begin{theorem} 
\label{thm:wll-inverse-image}
  Let $f:X\longrightarrow Y$ be a morphism of arithmetic varieties
  over $A$. Let $E$ be a normal crossing 
  divisor on $Y_{\mathbb{R}}$ and $D$ a normal crossing divisor on
  $X_{\mathbb{R}}$ such that $f^{-1}(E)\subset D$. Write $\underline
  X=(X_{\mathbb{R}},D)$ and $\underline Y=(Y_{\mathbb{R}},E)$. Then
  there is defined an inverse image morphism  
  \begin{displaymath}
    f^{\ast}:\cha^{\ast}(Y,\mathcal{D}_{\lgi})
    \longrightarrow 
    \cha^{\ast}(X,\mathcal{D}_{\lgi}).
  \end{displaymath}
  Moreover, it is a morphism of rings after tensoring with $\mathbb{Q}$.
  \hfill $\square$
\end{theorem}

\nnpar{Push-forward.}
We will state only the consequence of the integrability of log-log
forms. 

\begin{theorem}
\label{thm:pushforward-point}
  If $X$ is projective over $A$, then there is a direct image
  morphism of groups
  \begin{displaymath}
    f_{\ast}:\cha^{d+1}(X,\mathcal{D}_{\lgi})
    \longrightarrow 
    \cha^{1}(\Spec A),
  \end{displaymath}
  where $d$ is the relative dimension of $X$.\hfill $\square$
\end{theorem}

\nnpar{Relationship with other arithmetic Chow groups.}

Since we know the cohomology of the complex $\mathcal{D}_{\lgi}$, we
can make a comparison statement more precise than
\cite{BurgosKramerKuehn:cacg}, theorem 6.23.
\begin{theorem} \label{thm:7}
The structural morphism 
\begin{displaymath}
  \mathcal{D}_{\log,X}\longrightarrow
  \mathcal{D}_{\lgi,\underline{X}}
\end{displaymath}
induces a morphism 
\begin{displaymath}
  \cha^{\ast}(X,\mathcal{D}_{\log})
  \longrightarrow 
  \cha^{\ast}(X,\mathcal{D}_{\lgi})
\end{displaymath}
that is compatible with inverse images, intersection products and
arithmetic degrees. If $X$ is projective the isomorphism between
$\cha^{\ast}(X,\mathcal{D}_{\log})$ and the arithmetic Chow groups
defined by Gillet and Soul\'e (denoted $\cha^{\ast}(X)$) induce
morphisms 
\begin{equation}
\label{eq:34}
   \cha^{\ast}(X)
  \longrightarrow 
  \cha^{\ast}(X,\mathcal{D}_{\lgi})  
\end{equation}
also compatible with inverse images, intersection products and
arithmetic degrees. Moreover, if $D$ is empty and $X$ is projective
then the above morphisms are isomorphisms.  \hfill $\square$
\end{theorem}

\subsection{The $\mathcal {D}_{\log}$-complex of log-log forms with arbitrary singularities at infinity }

The arithmetic Chow groups defined by Gillet and Soul\'e for 
quasi-projective varieties use differential forms with arbitrary
singularities in the boundary. Therefore,in order to be able to
recover Gillet and Soul\'e arithmetic Chow groups we have to
 to introduce another variant of arithmetic Chow
groups, where we allow the differential forms to have arbitrary 
singularities in certain directions.

\nnpar{Mixing log, log-log and arbitrary singularities.} Let $X$ be a
complex algebraic manifold and $D$ a fixed normal crossing divisor of 
$X$. We write $\underline{X}=(X,D)$.   
\begin{definition} \label{def:20}
  For every Zariski open subset $U$ of $X$ we write
  \begin{displaymath}
    E^{\ast}_{\as,\underline
    X}(U)=\lim_{\substack{\longrightarrow\\\overline U} }\Gamma
    (\overline U,\mathscr{E}^{\ast}_{\overline U}\left 
      <B_{\overline U}\left<\overline D\right>\right>),
  \end{displaymath}
  where the limit is taken over all diagrams
  \begin{displaymath}
    \xymatrix{U\ar[r]^{\overline \iota}\ar[rd]^{\iota}
      &\overline U \ar[d]^{\beta }\\
      &X
    }
  \end{displaymath}
  such that $\overline \iota$ is an open immersion, $\beta $ is a
  proper morphism and $B_{\overline
  U}=\overline U\setminus U$, $\overline D=\beta ^{-1}(D)$ and
  $B_{\overline U}\cup \overline D$ are normal crossing divisors. 
\end{definition}

\begin{definition}
\label{def:19} 
Let $X$ be a complex algebraic manifold and $D$ a fixed normal 
crossing divisor of $X$. We write $\underline{X}=(X,D)$ as before.  
For any Zariski open subset $U\subseteq X$, we put
  \begin{displaymath}
    \mathcal{D}^{\ast}_{\as,\underline{X}}(U,p)=(\mathcal{D}^{\ast}_
    {\as,\underline{X}}(U,p),\dd_{\mathcal{D}})=(\mathcal{D}^{\ast}
    (E_{\as,\underline X}(U_{\mathbb{C}}),p),\dd_{\mathcal{D}}).
  \end{displaymath}
If $X$ is a smooth algebraic variety over $\mathbb{R}$, and $D$, $U$ 
are defined over $\mathbb{R}$, we then put
  \begin{displaymath}
    \mathcal{D}^{\ast}_{\as,\underline{X}}(U,p)=(\mathcal{D}^{\ast}_
    {\as,\underline{X}}(U,p),\dd_{\mathcal{D}})=(\mathcal{D}^{\ast}
    (E_{\as,\underline X}(U_{\mathbb{C}}),p)^{\sigma},\dd_{\mathcal{D}}),
  \end{displaymath}
  where $\sigma $ is as in section \ref{sec:mathc-d_log-compl}.
\end{definition}

Note that, when $X$ is quasi-projective, the varieties $\overline U$ of
definition \ref{def:20}, are not compactifications of $U$ but only
\emph{partial compactifications.} Therefore the sections of 
$\mathcal{D}^{\ast}_{\as,\underline{X}}(U,p)$ have three different
kinds of singularities. We can see this more concretely as follows. 
Let $Y$ be a closed subset of $X$ with $U=X\setminus Y$, and
let $\overline X$ be a 
smooth compactification of $X$ with $Z=\overline X\setminus X$. Let
$\eta $ be a section of
$\mathcal{D}^{\ast}_{\as,\underline{X}}(U,p)$. If we consider $\eta$
as a singular form on $\overline X$, then $\eta$ is log along $Y$ (in
the sense that it is log along a certain resolution of singularities
of $Y$), log-log along $D$ and has arbitrary
singularities along $Z$. Therefore, in general we have
\begin{displaymath}
  \mathcal{D}^{\ast}_{\as,\underline{X}}(U,p)\not =
  \mathcal{D}^{\ast}_{\as,\underline{U}}(U,p).
\end{displaymath}
Nevertheless, when $\underline X$ is clear from the context we will
drop it from the notation.

\begin{remark} \label{rem:3}
  If $X$ is projective, then the complexes of sheaves
  $\mathcal{D}^{\ast}_{\as,\underline X}$ and
  $\mathcal{D}^{\ast}_{\lgi,\underline X}$ agree. By contrast they do
  not agree when $X$ is quasi-projective. 
  Note moreover that when $X$ is quasi-projective the complex
  $\mathcal{D}^{\ast}_{\as,\underline X}$  does not compute the
  Deligne-Beilinson cohomology of $X$, but a mixture between
  Deligne-Beilinson cohomology and analytic Deligne
  cohomology. Nevertheless, as we will see, the local nature of the
  purity property of Deligne-Beilinson
  cohomology implies also a purity property for these complexes. 
\end{remark}

\nnpar{Logarithmic singularities and Blow-ups.} Let $X$ be a complex
manifold, $D\subseteq X$ a normal crossing divisor and $Y\subseteq X$
a codimension $e$ smooth subvariety such that the pair $(D,Y)$ has normal
crossings. Let $\pi :\widetilde X\longrightarrow X$ be the blow-up of
$X$ along $Y$. Write $\widetilde D=\pi ^{-1}(D)$ and $\widetilde Y=\pi
^{-1}(Y)$. Let $i:Y\longrightarrow X$ and $j:\widetilde Y
\longrightarrow \widetilde X$ denote the inclusions, and let
$g:\widetilde Y\longrightarrow Y$ denote the induced morphism. Observe
that $g$ is a projective bundle.

\begin{proposition} \label{prop:3} Let $p\ge 0$ be an integer.
  \begin{enumerate}
  \item If $Y\subset D$ then the morphism $\Omega ^{p}_{X}(\log
    D)\longrightarrow R\pi _{\ast}\Omega ^{p}_{\widetilde X}(\log
    \widetilde 
    D)$ is a quasi-isomorphism, i.e.,
    \begin{alignat*}{2}
      \pi _{\ast}\Omega ^{p}_{\widetilde X}(\log \widetilde
      D)&\cong \Omega ^{p}_{X}(\log
      D),\qquad &&\text{and} \\
      R^{q}\pi _{\ast}\Omega ^{p}_{\widetilde X}(\log \widetilde
      D)&=0,\qquad&&\text{for $q>0$.} 
    \end{alignat*}
  \item If $Y\not\subset D$ and $e>1$, then
    \begin{align*}
      \pi _{\ast}\Omega ^{p}_{\widetilde X}(\log \widetilde
      D\cup \widetilde Y)&\cong \Omega ^{p}_{X}(\log
      D),\\
      R^{q}\pi _{\ast}\Omega ^{p}_{\widetilde X}(\log \widetilde
      D\cup \widetilde Y)&=0,\qquad\text{for $q\not=0,e-1$, and}\\ 
      R^{e-1}\pi _{\ast}\Omega ^{p}_{\widetilde X}(\log \widetilde
      D\cup \widetilde Y)&\cong i_{\ast}(R^{e-1}g_{\ast}\Omega
      ^{p-1}_{\widetilde Y}(\log \widetilde D\cap \widetilde Y))\\
      &\cong
      i_{\ast}(\Omega
      ^{p-e}_{Y}(\log D\cap Y)\otimes 
      R^{e-1}g_{\ast}\Omega ^{e-1}_{\widetilde
      Y/Y}).
    \end{align*}
  \item If $Y\not\subset D$ and $e=1$, then $\pi=\Id$ and there is a
  short exact sequence
  \begin{displaymath}
    0\longrightarrow \Omega ^{p}_{X}
    \longrightarrow \Omega ^{p}_{X}(\log Y) 
    \longrightarrow i_{\ast} \Omega ^{p-1}_{Y}\longrightarrow 0. 
  \end{displaymath}
  \end{enumerate}
\end{proposition}
\begin{proof}
  The third statement is standard.

  The first statement is \cite{GuillenNavarro:cef} Proposition 4.4,
  (ii). 

  Using \cite{GuillenNavarro:cef} Proposition 4.4,
  (i), the fact that $i_{\ast}$ is an exact functor and that $g$ is a
  projective bundle we obtain
  \begin{align*}
    \pi _{\ast}\Omega ^{p}_{\widetilde X}(\log \widetilde
    D)&\cong \Omega ^{p}_{X}(\log
    D),\\
    \pi _{\ast} j_{\ast} \Omega ^{p-1}_{\widetilde Y}
    (\log \widetilde D\cap \widetilde Y)&\cong
    i _{\ast} g_{\ast} \Omega ^{p-1}_{\widetilde Y}
    (\log \widetilde D\cap \widetilde Y)\\
    &\cong i _{\ast} \Omega ^{p-1}_{Y}
    (\log D\cap  Y),\\
    R^{q}\pi _{\ast} \Omega ^{p}_{\widetilde X}(\log \widetilde
    D)&\cong R^{q} (\pi \circ j)_{\ast} \Omega ^{p}_{\widetilde Y}
    (\log \widetilde D\cap \widetilde Y)\\
    &\cong i_{\ast} R^{q} g_{\ast} \Omega ^{p}_{\widetilde Y}
    (\log \widetilde D\cap \widetilde Y)\\
    &\cong
    \begin{cases}
      i_{\ast}(\Omega
      ^{p-q}_{Y}(\log D\cap Y)\otimes 
      R^{q}g_{\ast}\Omega ^{q}_{\widetilde
      Y/Y}), & \text{if }1\le q < e,\\
    0, & \text{if }g\ge e. 
    \end{cases}
  \end{align*}
  Let $\mathscr{O}(1)$ be the ideal sheaf of $\widetilde Y$ in
  $\widetilde X$. We consider the exact sequence 
  \begin{displaymath}
    0\longrightarrow  \Omega ^{p}_{\widetilde X}(\log \widetilde
    D) \longrightarrow 
    \Omega ^{p}_{\widetilde X}(\log \widetilde
    D\cup \widetilde Y )\overset{\Res }{\longrightarrow }
    j_{\ast} \Omega ^{p-1}_{\widetilde Y}(\log \widetilde
    D\cap \widetilde Y ) \longrightarrow  0
  \end{displaymath}
  and the corresponding long exact sequence obtained by applying the
  functor $R\pi _{\ast}$. The connecting morphism of this long 
  exact sequence:
  \begin{multline*}
    R^{q-1}\pi _{\ast} j_{\ast} \Omega ^{p-1}_{\widetilde Y}(\log \widetilde
    D\cap \widetilde Y ) \cong
    i_{\ast}(\Omega
      ^{p-q}_{Y}(\log D\cap Y)\otimes 
      R^{q-1}g_{\ast}\Omega ^{q-1}_{\widetilde
      Y/Y})\\
    \longrightarrow
    R^{q}\pi _{\ast} \Omega ^{p}_{\widetilde X}(\log \widetilde
    D)
    \cong i_{\ast}(\Omega
      ^{p-q}_{Y}(\log D\cap Y)\otimes 
      R^{q}g_{\ast}\Omega ^{q}_{\widetilde
      Y/Y})
  \end{multline*}
  can be identified with the product by $c_{1}\mathscr{O}_{\widetilde
    X}(1)$, which is an isomorphism for $0<q\le e-1$.
  Then the result follows from this exact sequence.
\end{proof}

This proposition has the following consequence.

\begin{corollary}\label{cor:2}
  Let $X$ be a complex algebraic manifold and $Y$ a complex
  subvariety of 
  codimension $e$. Let $\widetilde X\longrightarrow X$ be an embedded
  resolution of singularities of $Y$ obtained as in
  \cite{Hironaka:rs}. Then 
  \begin{displaymath}
    R^{q}\pi _{\ast}\Omega ^{p}_{\widetilde X}(\log D)\cong
    \begin{cases}
      \Omega ^{p}_{X},& \text{ if }q=0,\\
      0,\text{ if $p<e$ or $0<q<e-1$}. 
    \end{cases}
  \end{displaymath}
\end{corollary}
\begin{proof}
  According to \cite{Hironaka:rs}, $\widetilde X$ is obtained by
  a series of elementary steps
  \begin{displaymath}
    \widetilde X=\widetilde X_{N}\longrightarrow \widetilde
    X_{N-1}\longrightarrow  \dots \longrightarrow \widetilde X_{0}=X, 
  \end{displaymath}
  where $\widetilde X_{k}$ is the blow-up of $\widetilde X_{k-1}$ 
  along a smooth
  subvariety $W_{k-1}$, contained in the strict transform of $Y$, therefore of
  codimension greater or equal than $e$. Moreover, if $D_{k}$ is the union of
  exceptional divisors up to the step $k$, then the pair
  $(D_{k},W_{k})$ has normal crossings. The result follows by
  applying proposition \ref{prop:3} to each blow-up. 
\end{proof}

The following theorem implies in particular the weak purity condition
for the complex $\mathcal{D}_{\as,\underline X}$.

\begin{theorem} \label{thm:9}
  Let $\underline X=(X,D)$ be as above. Let $Y\subset X$ be a Zariski
  closed 
  subset of codimension  greater or equal than $p$. Let $c$ be the number of connected
  components of $Y$ of codimension $p$. Then the natural morphisms 
  \begin{displaymath}
    H^{n}_{\mathcal{D}_{\as},Y}(X,p)\longrightarrow
    H^{n}_{Y}(X,\mathbb{R}(p)),
  \end{displaymath}
  is an isomorphism 
  for all integers $n$. Therefore
  \begin{align*}
    H^{n}_{\mathcal{D}_{\as},Y}(X,p)&= 0,\qquad
    \text{for } 
    n<2p,\text{ and} \\
    H^{2p}_{\mathcal{D}_{\as},Y}(X,p)&\cong
    \mathbb{R}(p)^{c}.  
  \end{align*}
\end{theorem}
\begin{proof}
  We fix a diagram 
  \begin{displaymath}
    \xymatrix{U\ar[r]^{\overline \iota}\ar[rd]^{\iota}
      &\overline U \ar[d]^{\beta }\\
      &X
    }
  \end{displaymath}
  such that $\overline \iota$ is an open immersion, $\beta $ is a
  proper morphism and $B=\overline U\setminus U$, $\overline D=\beta
  ^{-1}(D)$ and 
  $B \cup \overline D$ are normal crossing
  divisors. Hence $\overline U$ is an embedded resolution of
  singularities of $Y$. We assume moreover that $\overline U$ is
  obtained from $X$ as $\widetilde X$ is obtained from $X$ in
  corollary \ref{cor:2}. 

  By theorem \ref{thm:fq}, the complexes
  $\mathcal{D}^{\ast}_{\as,\underline{X}}(X,p)$ and
  $\mathcal{D}^{\ast}_{\as,\underline{X}}(U,p)$ are quasi-isomorphic
  to the complexes 
  $\mathcal{D}^{\ast}(E^{\ast}_{X},p)$ and 
  $\mathcal{D}^{\ast}(E^{\ast}_{\overline
  U}\langle B\rangle,p)$, respectively.

  By the definition of the Deligne complex, and theorem
  \cite{Burgos:CDB} 2.6.2 there are
  quasi-isomorphisms 
  \begin{align*}
    \mathcal{D}^{\ast}(E^{\ast}_{X},p)&\longrightarrow 
    s\left(E_{X,\mathbb{R}}^{\ast}(p)\rightarrow E^{\ast}_{X}\left /
      F^{p}E^{\ast}_{X}\right .\right)\\
    \mathcal{D}^{\ast}(E^{\ast}_{\overline
      U}\langle B\rangle,p)&\longrightarrow 
    s\left(E_{\overline U}^{\ast}\langle B \rangle_{\mathbb{R}}(p)
    \rightarrow E^{\ast}_{\overline U}\langle B \rangle
    \left / F^{p}E^{\ast}_{\overline U}\langle B \rangle\right.\right)
  \end{align*}
  By corollary \ref{cor:2} and theorem \ref{thm:logfq} the natural
  morphism 
  \begin{displaymath}
    E^{\ast}_{X}\left /
      F^{p}E^{\ast}_{X}\right . \longrightarrow 
    E^{\ast}_{\overline U}\langle B \rangle
    \left / F^{p}E^{\ast}_{\overline U}\langle B \rangle\right.
  \end{displaymath}
  is a quasi-isomorphism. Hence the morphism 
  \begin{displaymath}
    s\left(\mathcal{D}^{\ast}(E^{\ast}_{X},p)\rightarrow 
    \mathcal{D}^{\ast}(E^{\ast}_{\overline
      U}\langle B\rangle,p)\right)
    \longrightarrow 
    s\left(E_{X,\mathbb{R}}^{\ast}(p) \rightarrow 
    E_{\overline U}^{\ast}\langle B \rangle_{\mathbb{R}}(p)
    \right)
  \end{displaymath}
  is a quasi-isomorphism. Since the left hand complex computes 
  $H^{n}_{\mathcal{D}_{\as},Y}(X,p)$ and the right hand
  complex computes 
    $H^{n}_{Y}(X,\mathbb{R}(p))$ we obtain the first statement of the
  theorem. The second statement follows form the purity of singular
  cohomology. 
\end{proof}

Summing up the properties of the complex
$\mathcal{D}_{\as,\underline{X}}$ we obtain

\begin{theorem} 
\label{thm:6}
The complex $\mathcal{D}_{\as,\underline{X}}$ is a $\mathcal{D}_
{\log}$-complex on $X$. Moreover, it is a pseudo-associative and
commutative $\mathcal{D}_{\log}$-algebra and satisfies the weak purity
condition (\cite{BurgosKramerKuehn:cacg} Definition 3.1).
\hfill $\square$
\end{theorem}

\subsection{Arithmetic Chow rings with arbitrary singularities at
  infinity}
Let $A$, $X$, $\DD$ and $\underline X$ be as at the beginning of
section \ref{sec:acgwllg}. Applying \cite{BurgosKramerKuehn:cacg},
section 4, we define the arithmetic Chow groups
$\cha^{\ast}(X,\mathcal{D}_{\as})$. Then the theorems
\ref{thm:ll-product}, \ref{thm:wll-inverse-image} and \ref{thm:7} are
also true for these groups. For theorem \ref{thm:pushforward-point} to
be true we need $X$ to be projective, but in this case there is no
difference between $\cha^{\ast}(X,\mathcal{D}_{\as})$
and $\cha^{\ast}(X,\mathcal{D}_{\lgi})$.

  Since
  $\mathcal{D}_{\as,\underline{X}}$ satisfies the weak purity
  property, the analogue of
  theorem \ref{thm:16} is reads as follows.
  \begin{theorem}\label{thm:10}
    The following sequence is exact:
    \begin{displaymath}
      \CH^{p-1,p}(X) \stackrel{\rho}{\longrightarrow}
      \widetilde{\mathcal{D}}^{2p-1}_{\as}(X,p) 
      \stackrel{\amap}{\longrightarrow} 
      \cha^p(X,\mathcal{D}_{\as}) \stackrel{\zeta}{\longrightarrow}
      \CH^p(X) \longrightarrow 0.
    \end{displaymath}
    \hfill $\square$
  \end{theorem} 

Another consequence of theorem \ref{thm:9} is the analogue of proposition
\ref{prop:6}, that is proved in the same way. 

\begin{proposition} \label{prop:9}
Let $X$ be a smooth real variety and $D$ a normal crossing divisor. Put
$\underline X=(X,D)$. Let $y$ be a $p$-codimensional cycle 
on $X$ with support $Y$. Then,
the class of the cycle $(\omega,g)$ in $H^{2p}_{\mathcal{D}_
{\as},Y}(X,p)$ is equal to the class of $y$, if and only if  
\begin{align}
-2\partial\bar{\partial}[g]_{X}=[\omega]-\delta_{y}. 
\end{align}
\hfill $\square$
\end{proposition}  

From this proposition and theorem \ref{thm:10} we obtain the analogue
of \cite{BurgosKramerKuehn:cacg} theorem 6.23:

\begin{theorem}\label{thm:8}  Let
$\cha^{p}(X)$ be the arithmetic Chow groups defined by Gillet and
Soul\'e.  
If $\DD=\emptyset$ then
the assignment $$[y,(\omega_{y},\widetilde{g}_{y})]\mapsto
[y,2(2\pi i)^{d-p+1}[g_{y}]_{X}]$$ induces a well defined 
isomorphism
\begin{displaymath}
\Psi:\cha^{p}(X,\mathcal{D}_{\as})\longrightarrow\cha^{p}(X),
\end{displaymath}
which is compatible with products and pull-backs.
\end{theorem}

\begin{remark}
Note that, if $f:X\longrightarrow Y$ is a proper morphism between
arithmetic varieties over $A$ and such that
$f_{\mathbb{R}}:X_{R}\longrightarrow Y_{R}$ is smooth, then there is a
covariant $f$-pseudo morphism (see \cite{BurgosKramerKuehn:cacg}
definition 3.71) that induces a push-forward morphism   
\begin{displaymath}
  f_{\ast}:\cha^{p}(X,\mathcal{D}_{\as})\longrightarrow
  \cha^{p}(Y,\mathcal{D}_{\as}).
\end{displaymath}
This push-forward is compatible with the push-forward defined by Gillet and
Soul\'e.
\end{remark}

\begin{remark}
  We can define $\mathcal{D}_{\as, \wlg,\underline X}$ as
  $\mathcal{D}_{\as,\underline X}$ replacing pre-log and pre-log-log
  forms for log and log-log forms. We then obtain a theory of
  arithmetic Chow groups $\cha^{p}(X,\mathcal{D}_{\as,\wlg})$ with
  analogous properties. Note however that since we have not established
  the weak purity property of pre-log forms, we do not have the
  analogue of theorem \ref{thm:8}. 
\end{remark}

\section{Bott-Chern forms of log-singular hermitian 
vector bundles}
\label{sec:chern-bott-chern}

The arithmetic intersection theory of Gillet and Soul\'e is complemented
by an arithmetic $K$-theory and a theory of characteristic classes. A
main ingredient of the theory of arithmetic characteristic classes are
the Chern forms and Bott-Chern forms of hermitian vector bundles. In 
this section after defining the class of singular metrics considered
in this paper, we will generalize the theory of Chern forms and
Bott-Chern forms to include this class of singular metrics.

\subsection{Chern forms for hermitian metrics}
\label{sec:chern-forms-bott}
Here we recall the Chern-Weil theory of characteristic classes of 
hermitian vector bundles. By a hermitian metric we will always mean a 
smooth hermitian metric.

\nnpar{Chern forms.} Let $B\subset \mathbb{R}$ be a subring, let $\phi
\in B[[T_{1},\dots ,T_{n}]]$ be any symmetric power series in $n$
variables and let $M_{n}(\mathbb{C})$ be the algebra of $n\times n$
complex matrices. For every $k\ge 0$, let $\phi ^{(k)}$ be the
homogeneous component of $\phi $ of degree $k$. We will denote also by
$\phi ^{(k)}:M_{n}(\mathbb{C})\longrightarrow \mathbb{C}$ the unique
polynomial map which is invariant under conjugation by
$\GL_{n}(\mathbb{C})$ and whose value in the diagonal matrix
$\diag(\lambda _{1},\dots,\lambda _{n})$, $\lambda _{i}\in
\mathbb{C}$, 
is $\phi ^{(k)}(\lambda _{1},\dots ,\lambda _{n})$. More generally, if
$A$ is any $B$-algebra, $\phi ^{(k)}$ defines a map $\phi
^{(k)}:M_{n}(A)\longrightarrow A$, and if $I\subset A$ is a nilpotent
subalgebra we can define $\phi =\sum_{k}\phi
^{(k)}:M_{n}(I)\longrightarrow A$. 

Let $\overline E=(E,h)$ be a hermitian vector bundle of rank $n$ on a
complex manifold $X$. Let $\xi =\{\xi _{1},\dots ,\xi _{n}\}$ be a frame
for $E$ in an open subset $V\subset X$. We denote by $h(\xi
)=(h_{ij}(\xi ))$ the matrix of $h$ in the frame $\xi $. Let $K(\xi )$
be the curvature matrix $K(\xi )=\overline \partial(\partial h(\xi
)\cdot  h(\xi
)^{-1})$. The Chern form associated to $\phi $ and $\overline E$ is
the form  
\begin{displaymath}
  \phi (\overline E)=\phi (-K(\xi))\in E^{\ast}_{V}.
\end{displaymath}

\nnpar{Basic properties.} The following properties of the Chern forms
are well known.  
\begin{theorem}\label{thm:20}
  \begin{enumerate}
  \item By the invariance of the $\phi ^{(k)}$, the Chern form $\phi
    (\overline E)$ is independent of the choice of frame $\xi
    $. Therefore, it globalizes to a differential form $\phi (\overline
    E)\in E^{\ast}_{X}$. 
  \item  The Chern forms are closed.
  \item The component $\phi ^{(k)}$ belongs to 
    $\mathcal{D}^{2k}(E_{X},k)=E^{k,k}_{X}\cap
    E^{2k}_{X,\mathbb{R}}(k)$.
  \item If $(X_{\mathbb{R}})=(X,F_{\infty})$ is a real manifold, the
  vector bundle $E$ is defined over $\mathbb{R}$ and the hermitian
  metric 
  $h$ is invariant under $F_{\infty}$, then $\phi ^{(k)}(E,h)\in
  \mathcal{D}^{2k}(E_{X},k)^{\sigma }$, where $\sigma $ is as in
  definition \ref{def:1}.
  \end{enumerate}
  \hfill$\square$
\end{theorem}

\nnpar{Chern classes.} Since the Chern forms are closed, they
represent cohomology classes $\phi (E)=[\phi (E,h)]\in \bigoplus_{k}
H^{2k}(\mathcal{D}(E_{X},k))$. If $X$ is projective, then
$\bigoplus_{k}H^{2k}(\mathcal{D}(E_{X},k))= \bigoplus_{k}
H^{2k}_{\mathcal{D}}(X,\mathbb{R}(k))$, hence we obtain classes in
Deligne-Beilinson cohomology $\phi (E)\in \bigoplus_{k}
H^{2k}_{\mathcal{D}}(X,\mathbb{R}(k))$.  

Note that, to simplify the notation, the function $\phi $ will have
different meanings according to its arguments. For instance $\phi
(E,h)=\phi (\overline E)$ will mean the Chern form that depends on
the bundle and the metric, whereas $\phi (E)$ will mean the Chern
class that depends only on the bundle.

When $X$ is algebraic, not necessarily projective, by means of
smooth at infinity hermitian metrics (see
\cite{BurgosWang:hBC} for a precise definition), the Chern-Weil theory
also allows us to 
construct Chern classes in Deligne-Beilinson cohomology .

Recall that there are Chern classes defined 
in the Chow ring $\phi (E)_{\CH}\in \CH^{\ast}(X)$ and they are
compatible with the Chern classes in cohomology. More precisely, we have
\begin{proposition}
  The composition 
  \begin{displaymath}
    \CH^{k}(X)\overset{\cl}{\longrightarrow }
    H^{2k}_{\mathcal{D}}(X,\mathbb{R}(p))\longrightarrow 
    H^{2k}(\mathcal{D}(E_{X},k))
  \end{displaymath}
  sends $\phi ^{(k)}(E)_{\CH}$ to $\phi ^{(k)}(E)$.
  \hfill $\square$
\end{proposition}

\subsection{Bott-Chern forms for hermitian metrics}
Here we recall the theory of Bott-Chern forms.  For more details we
refer to \cite{Soule:lag}, \cite{BurgosWang:hBC}, \cite{Burgos:hvbcc}.

\nnpar{Bott-Chern forms.} 
Let 
\begin{displaymath}
  \overline {{\mathcal{E}}}: 0\longrightarrow (E',h')\longrightarrow (E,h)
\longrightarrow (E'',h'')\longrightarrow 0
\end{displaymath}
be a short exact sequence of hermitian vector bundles; by this we mean
a short exact sequence of vector bundles, where each vector bundle is
equipped with an arbitrarily chosen hermitian metric. Let $\phi $ be
as in \ref{sec:chern-forms-bott} and assume $E$ has rank $n$.

The Chern 
classes behave additively with respect to exact sequences, i.e.,
\begin{displaymath}
  \phi (E)=\phi (E'\oplus E'').
\end{displaymath}
In general this is not true for the Chern forms. This lack of
additivity at the level of Chern forms is measured by the Bott-Chern
forms. 

The fundamental result of the theory of Bott-Chern forms is the
following theorem (see \cite{BottChern:hvb}, \cite{BismutGilletSoule:at},
\cite{GilletSoule:vbhm}). 
\begin{theorem}\label{thm:2}
  There is a unique way to attach to every sequence $\overline
  {{\mathcal{E}}}$ as above, a form $\widetilde {\phi} (\overline
  {{\mathcal{E}}})$ in
  \begin{displaymath}
    \bigoplus_{k}\widetilde {\mathcal{D}}^{2k-1}(E_{X},k)
    =\bigoplus_{k}\mathcal{D}^{2k-1}(E_{X},k)/\Img(\dd_{\mathcal{D}})
  \end{displaymath}
  satisfying the following properties
  \begin{enumerate}
  \item $\dd_{\mathcal{D}}\widetilde {\phi }(\overline {\mathcal{E}})=
    \phi (E'\oplus E'',h'\oplus h'')-\phi (E,h)$.
  \item $f^{\ast}\widetilde {\phi }(\overline{{\mathcal{E}} })=\widetilde {\phi
  }(f^{\ast}\overline {\mathcal{E}} )$, for every holomorphic map
  $f:X\longrightarrow Y$. 
  \item If $(E,h)=(E',h')\overset{\perp}{\oplus}(E'',h'')$, then
  $\widetilde {\phi }(\overline {{\mathcal{E}}})=0$.
  \end{enumerate}
\end{theorem}

There are different methods to construct Bott-Chern forms. We will
introduce a variant of the method used in \cite{GilletSoule:vbhm} and
that is the dual of the construction used in \cite{BurgosWang:hBC}.



\nnpar{The first transgression bundle.} Let $\mathcal{O}(1)$ be the
dual of the tautological bundle of $\mathbb{P}^{1}$ with the standard
metric. If $(x:y)$ are projective coordinates of
$\mathbb{P}^{1}_{\mathbb{C}}$ then 
$x$ and $y$ are generating global sections of $\mathcal{O}(1)$ with
norm 
\begin{displaymath}
  \|x\|^{2}=\frac{x\overline x}{x\overline x+y\overline y}\quad
  \text{and}\quad  
  \|y\|^{2}=\frac{y\overline y}{x\overline x+y\overline y}.
\end{displaymath}

Let 
\begin{displaymath}
  \overline {{\mathcal{E}}}: 0\longrightarrow (E',h')\longrightarrow (E,h)
\longrightarrow (E'',h'')\longrightarrow 0
\end{displaymath}
be a short exact sequence of hermitian vector bundles such that $h'$ is
induced by $h$.

Let $p_{1},p_{2},$ be the first and
the second projection of $X\times \mathbb{P}^{1}_{\mathbb{C}}$,
respectively. We write $
E(n)=p_{1}^{\ast}E\otimes p_{2}^{\ast}\mathcal{O}(n)$. On this vector bundle
we consider the metric
induced by $h$ and the standard metric of $\mathcal{O}(n)$, and we
denote by $\overline E(n)$ this hermitian vector bundle. Analogously,
we write $E''(n)=p_{1}^{\ast}E''\otimes p_{2}^{\ast}\mathcal{O}(n)$
and  denote by $\overline E''(n)$ the corresponding hermitian vector
bundle.

\begin{definition}\label{def:8}
  The \emph{first transgression bundle} $\tr_{1}(\overline {\mathcal{E}})$ is
  the kernel of the morphism
  \begin{displaymath}
    \begin{matrix}
      \overline E(1)\oplus \overline E''(1)&\longrightarrow &E''(2)\\
      (s,t)&\longmapsto&s\otimes x-t\otimes y
    \end{matrix}
  \end{displaymath}
with the induced metric. 
\end{definition}
Note that the definition of
$\tr_{1}(\overline {\mathcal{E}})$ includes the metric; therefore, the
expression $\phi (\tr_{1}(\overline {\mathcal{E}}))$ means the Chern
form of the hermitian vector bundle $\tr_{1}(\overline {\mathcal{E}})$
and not its Chern class.

The key property of the first transgression bundle is the following.
We denote by $i_{0}$ and $i_{\infty}$ the morphisms
$X\longrightarrow X\times \mathbb{P}^{1}$ given by
\begin{align*}
  i_{0}(p)&=(p,(0:1)),\\
  i_{\infty}(p)&=(p,(1:0)).
\end{align*}
Then, $i_{0}^{\ast}(\tr_{1}(\overline {\mathcal{E}}))$ is isometric to
$(E,h)$ and $i_{\infty}^{\ast}(\tr_{1}(\overline {\mathcal{E}}))$ is
isometric to $(E',h')\overset{\perp}{\oplus}(E'',h'')$.

\nnpar{The construction of Bott-Chern forms.} Let $t=x/y$ be the
absolute coordinate of 
$\mathbb{P}^{1}$. Let us consider the current $W_{1}=\left[-\frac{1}{2}\log
(t\overline t )\right]$ on $\mathbb{P}^{1}$ given by 
\begin{displaymath}
  W_{1}(\eta)=\left[-\frac{1}{2}\log( t\overline t) \right](\eta)=
  -\frac{1}{2 \pi i}\int_{\mathbb{P}^{1}}  \frac{\eta}{2} \log( t
  \overline t). 
\end{displaymath}
By the Poincar\'e Lelong equation
\begin{equation}\label{eq:31}
  -2\partial\overline \partial \left[-\frac{1}{2}\log( t\overline t) \right]=
  \delta _{(1:0)}-\delta _{(0:1)}.
\end{equation}
\begin{definition}\label{def:2}
  Let $X$ be a complex manifold, $\overline {\mathcal{E}}$ an exact sequence of
  hermitian vector bundles 
  \begin{displaymath}
    \overline {{\mathcal{E}}}: 0\longrightarrow (E',h')\longrightarrow (E,h)
    \longrightarrow (E'',h'')\longrightarrow 0,
  \end{displaymath}
  such that the metric $h'$ is induced by the metric $h$.
  The \emph{Bott-Chern form}
  associated to the exact sequence $\overline {\mathcal{E}}$ is the
  differential
  form over $X$ 
  \begin{displaymath}
    \phi(\overline {\mathcal{E}})= W_{1}
    (\phi(\tr_{1}(\overline {\mathcal{E}}))) = 
    -\frac{1}{2 \pi i}\int_{\mathbb{P}^{1}} 
    \phi(\tr_{1}(\overline {\mathcal{E}})) \frac{1}{2} \log (t \overline t).
  \end{displaymath}
\end{definition}
Note that we use also the letter $\phi $ to denote the Bott-Chern form
associated to a power series $\phi $ because the meaning of $\phi
(\overline {\mathcal{E}})$ is determined again by the argument
$\overline {\mathcal{E}}$, which, in this case, is an exact sequence
of hermitian vector bundles.

\begin{definition}\label{def:3} 
If $\overline {\mathcal{E}}$ is an exact sequence as above, but such that
  $h'$ is not the metric induced by $h$, then we consider the exact
  sequences
\begin{displaymath}
  \lambda^{1}\overline {\mathcal{E}}:
  0\longrightarrow (E',\widetilde {h'})\longrightarrow (E,h)
\longrightarrow (E'',h'')\longrightarrow 0,
\end{displaymath}
where $\widetilde {h'}$ is the hermitian metric induced by $h$, and
\begin{displaymath}
  \lambda^{2}\overline {\mathcal{E}}:
  0\longrightarrow 0 \longrightarrow (E'\oplus E'',\widetilde
  {h'}\oplus h'')
  \longrightarrow (E'\oplus E'',h'\oplus h'') \longrightarrow 0.
\end{displaymath}
The \emph{Bott-Chern form} associated to the exact sequence
  $\overline {\mathcal{E}}$ is
  \begin{displaymath}
    \phi (\overline {\mathcal{E}})= \phi (\lambda ^{1}\overline {\mathcal{E}})+
    \phi (\lambda ^{2}\overline {\mathcal{E}}).
  \end{displaymath}
\end{definition}

\begin{proposition}
  If $\overline {\mathcal{E}}$ is an exact sequence as above with $h'$ induced
  by $h$, then the Bott-Chern forms obtained from definition
  \ref{def:2} and definition \ref{def:3} agree.
\end{proposition}
\begin{proof}
  In this case we have $\lambda ^{1}\overline {\mathcal{E}}=\overline
  {\mathcal{E}}$. Thus 
  we have to show that $\phi (\lambda ^{2}\overline {\mathcal{E}})=0$. But
  $\tr_{1}(\lambda ^{2}(\overline {\mathcal{E}} ))$ is the bundle
  $p_{1}^{\ast}(E'\oplus E'')$ with the hermitian metric $h'\oplus
  h''$ that does not depend 
  on the coordinate of $\mathbb{P}^{1}$. Therefore, we have
  \begin{displaymath}
    \phi (\lambda ^{2}\overline {\mathcal{E}})= -\frac{1}{2 \pi
  i}\int_{\mathbb{P}^{1}} \phi (E'\oplus E'',h'\oplus h'') \frac{1}{2}
  \log( t\overline t)= 0.
  \end{displaymath}
\end{proof}

It is easy to see that the forms $\phi (\overline {\mathcal{E}})$ belong to 
$\bigoplus_{k}\mathcal{D}^{2k-1}(E_{X},k)$. 
We will denote by $\widetilde \phi (\overline {\mathcal{E}})$ the
class of $\phi 
(\overline {\mathcal{E}})$ in the group
\begin{displaymath}
      \bigoplus_{k}\widetilde {\mathcal{D}}^{2k-1}(E_{X},k)
    =\bigoplus_{k}\mathcal{D}^{2k-1}(E_{X},k)/\Img(\dd_{\mathcal{D}}).
\end{displaymath}

\begin{proposition}
  The classes $\widetilde \phi (\overline {\mathcal{E}})$ satisfy the
  properties of theorem \ref{thm:2}.
\end{proposition}
\begin{proof}
  The first property follows from the Poincar\'e lemma (see for
  instance \cite{Soule:lag}). The second property is clear because all
  the ingredients of the
  construction are functorial. We prove the third property. If
  $\overline {\mathcal{E}}$ is a split exact sequence with $(E,h)=
(E',h') \overset{\perp}{\oplus}
  (E'',h'')$ and the obvious morphisms, then   
  \begin{displaymath}
    \tr_{1}(\overline {\mathcal{E}})=\overline E'(1)\oplus \overline E''(0)
  \end{displaymath}
  with the induced metrics. Let $\omega $ be the first Chern form of
  the line bundle $\overline{ \mathcal{O}}_{\mathbb{P}^{1}}(1)$. Then, we find 
  \begin{displaymath}
    \phi(\overline E'(1)\oplus \overline E''(0))
=p_{1}^{\ast}(a)+p_{1}^{\ast}(b)\land p_2^{\ast}(\omega),
  \end{displaymath}
  where $a$ and $b$ are suitable forms on $X$. Now we get
  \begin{align*}
    \phi (\overline {\mathcal{E}})&= -\frac{1}{2 \pi
      i}\int_{\mathbb{P}^{1}} (p_{1}^{\ast}(a)+p_{1}^{\ast}(b)\land
    p_2^{\ast} (\omega))
    \frac{1}{2} \log( t  \overline t)\\
    &=-\frac{1}{2 \pi i}\,a\land\int_{\mathbb{P}^{1}} \frac{1}{2} \log(
    t \overline t) -\frac{1}{2 \pi i}\, b\land\int_{\mathbb{P}^{1}}
    \frac{\omega }{2} \log( t \overline t)=0.
  \end{align*}
\end{proof}

\nnpar{Change of metrics.} Of particular importance is the Bott-Chern
form associated to a change of hermitian metrics. Let $E$ be a
holomorphic vector bundle of rank $n$ with two hermitian metrics $h$
and $h'$.   We
denote by $\tr_{1}(E,h,h')$ the first transgression bundle associated
to the short exact sequence
\begin{displaymath}
  0\longrightarrow 0\longrightarrow (E,h)\longrightarrow (E,h')
\longrightarrow 0.
\end{displaymath}

Explicitly, $\tr_{1}(E,h,h')$ is isomorphic to $p_{1}^{\ast}E$ with the
embedding   
\begin{displaymath}
  \begin{matrix}
    p_{1}^{\ast}E&\longrightarrow &\overline E(1)\oplus \overline E'(1)\\
    s&\longmapsto&(s\otimes y,s\otimes x);
  \end{matrix}
\end{displaymath}
here $\overline E = (E, h)$ and $\overline E'=(E,h')$. Therefore, if
$\xi$ is a local frame for $E$ on an open set $U$, then it determines
a local frame for $\tr(E,h,h')$, also denoted $\xi$, on $U\times
\mathbb{P}^{1}$. In this frame the metric is given by the matrix
\begin{equation}\label{eq:32}
  \frac{y\overline y h(\xi)+ x\overline x h'(\xi)}
  {x\overline x+ y \overline y}.  
\end{equation}

\begin{definition}
  Let $X$ be a complex manifold, $E$ be a complex vector bundle of
  rank $n$,
  $h$, $h'$ two hermitian metrics on $E$, and $\phi $ as in section
  \ref{sec:chern-forms-bott}. The \emph{Bott-Chern form}
  associated to the change of metric $(E,h,h')$ is the Bott-Chern form 
  associated to the short exact sequence
  \begin{displaymath}
    0\longrightarrow 0\longrightarrow 
    (E,h)\longrightarrow (E,h')\longrightarrow 0.
  \end{displaymath}
   We will denote this form by $\phi (E,h,h')$ or, if $E$ is
   understood, $\phi (h,h')$. This form satisfies
  \begin{equation}\label{eq:19}
    \dd_{\mathcal{D}}\phi (E,h,h')=-2\partial\bar \partial \phi
    (E,h,h')=\phi (E,h')-\phi (E,h).
  \end{equation}
\end{definition}

\subsection{Iterated Bott-Chern forms for hermitian metrics}
 
The theory of Bott-Chern forms can
be iterated defining higher Bott-Chern forms for exact $k$-cubes of
hermitian vector bundles. This theory provides explicit
representatives of 
characteristic classes for higher $K$-theory (see
\cite{BurgosWang:hBC}, \cite{Burgos:hvbcc}).

\nnpar{Exact squares.} Let $\left<-1,0,1\right>$ be the category
associated to the ordered set $\{-1,0,1\}$.
\begin{definition}
  A \emph{square of vector bundles} over $X$ is a functor from the
  category $\left<-1,0,1\right>^{2}$ to the category of vector bundles
  over $X$. Given a square of vector bundles $\mathcal{F}$ and numbers
  $i\in \{1,2\}$, $j\in \{-1,0,1\}$, then the \emph{$(i,j)$-face} of
  $\mathcal{F}$, denoted by $\partial^{j}_{i}\mathcal{F}$, is the sequence 
  \begin{align*}
    \partial^{j}_{1}\mathcal{F}&: \mathcal{F}_{j,-1}\longrightarrow
    \mathcal{F}_{j,0}\longrightarrow
    \mathcal{F}_{j,1},\\
    \partial^{j}_{2}\mathcal{F}&: \mathcal{F}_{-1,j}\longrightarrow
    \mathcal{F}_{0,j}\longrightarrow
    \mathcal{F}_{1,j}.\\
  \end{align*}
  A square of vector bundles is called \emph{exact}, if all the faces
  are short exact sequences. A \emph{hermitian exact square}
  $\overline {\mathcal{F}}$ is an exact square $\mathcal{F}$ such that
  the vector bundles $\mathcal{F}_{i,j}$ are equipped with arbitrarily
  chosen hermitian
  metrics. If $\overline{\mathcal{F}}$ is a hermitian exact square,
  then the faces of $\overline{\mathcal{F}}$ are provided with the
  induced hermitian metrics.
The reader is referred to \cite{BurgosWang:hBC} for the definition of
exact $n$-cubes.
\end{definition}

Let $\phi $ as before and
let $\overline {\mathcal{F}}$ be a hermitian exact square of vector
bundles over $X$ such that $\mathcal{F}_{0,0}$ has rank $n$. Then the
form 
\begin{displaymath}
    \phi (\partial_{1}^{-1} 
    \overline{\mathcal{F}}\oplus \partial_{1}^{1} 
    \overline{\mathcal{F}})-\phi (\partial_{1}^{0} 
    \overline{\mathcal{F}})-\phi (\partial_{2}^{-1} 
    \overline{\mathcal{F}}\oplus \partial_{2}^{1} 
    \overline{\mathcal{F}})+\phi (\partial_{2}^{0} 
    \overline{\mathcal{F}})
\end{displaymath}
is closed in the complex $\bigoplus
_{p}\mathcal{D}^{\ast}(E_{X},p)$. The iterated Bott-Chern form
is a differential form 
$$\phi (\overline {\mathcal{F}})\in 
\bigoplus_{p}\mathcal{D}^{2p-2}(E_{X},p)$$ 
such that  
\begin{displaymath}
  \dd_{\mathcal{D}}\phi(\overline {\mathcal{F}})=
  \phi (\partial_{1}^{-1} 
  \overline{\mathcal{F}}\oplus \partial_{1}^{1} 
  \overline{\mathcal{F}})-\phi (\partial_{1}^{0} 
  \overline{\mathcal{F}})-\phi (\partial_{2}^{-1} 
  \overline{\mathcal{F}}\oplus \partial_{2}^{1} 
  \overline{\mathcal{F}})+\phi (\partial_{2}^{0} 
  \overline{\mathcal{F}}).
\end{displaymath}

\nnpar{The second transgression bundle.} 
\begin{definition} \label{def:10}
  Let $\overline {\mathcal{F}}$
  be a hermitian exact square such that for $j=-1,0,1$, the hermitian
  metrics of the vector bundles $\mathcal{F}_{j,-1}$ and
  $\mathcal{F}_{-1,j}$ are induced by the metrics of $\mathcal{F}_{j,0}$ and
  $\mathcal{F}_{0,j}$, respectively.
  The \emph{second transgression bundle} associated to $\overline
  {\mathcal{F}}$ is the 
  hermitian vector bundle on $X\times \mathbb{P}^{1}\times
  \mathbb{P}^{1}$ given by 
  \begin{displaymath}
    \tr_{2}(\overline {\mathcal{F}})=
    \tr_{1}\left(\tr_{1}(\partial^{-1}_{2}\overline {\mathcal{F}})\longrightarrow 
    \tr_{1}(\partial^{0}_{2}\overline {\mathcal{F}})\longrightarrow
    \tr_{1}(\partial^{1}_{2}\overline {\mathcal{F}})\right).
  \end{displaymath}
\end{definition}

The second transgression bundle satisfies 
\begin{equation}
  \label{eq:29}
  \begin{aligned}
    \tr_{2}(\overline {\mathcal{F}})|_{X\times \mathbb{P}^{1}\times (0:1)}&=
    \tr_{1}(\partial^{0}_{2}\overline {\mathcal{F}}),\\
    \tr_{2}(\overline {\mathcal{F}})|_{X\times \mathbb{P}^{1}\times (1:0)}&=
    \tr_{1}(\partial^{-1}_{2}\overline {\mathcal{F}})\overset{\perp}{\oplus}
    \tr_{1}(\partial^{1}_{2}\overline {\mathcal{F}}),\\
    \tr_{2}(\overline {\mathcal{F}})|_{X\times (0:1)\times \mathbb{P}^{1}}&=
    \tr_{1}(\partial^{0}_{1}\overline {\mathcal{F}}),\\
    \tr_{2}(\overline {\mathcal{F}})|_{X\times (1:0)\times \mathbb{P}^{1}}&=
    \tr_{1}(\partial^{-1}_{1}\overline {\mathcal{F}})\overset{\perp}{\oplus}
    \tr_{1}(\partial^{1}_{1}\overline {\mathcal{F}}).
  \end{aligned} 
\end{equation}

\nnpar{The second Wang current.} On $\mathbb{P}^{1}\times
\mathbb{P}^{1}$ we put homogeneous coordinates 
$((x_{1}:y_{1}),(x_{2}:y_{2}))$; let $t_{1}=x_{1}/y_{1}$ and
$t_{2}=x_{2}/y_{2}$.
\begin{definition}
  The \emph{second Wang current} is the current on
  $\mathbb{P}^{1}\times \mathbb{P}^{1}$ given by
  $$
  W_{2}=\frac{1}{4} \left[\log( t_{1}\overline t_{1})
    \left(\frac{\dd t_{2}}{t_{2}}-
      \frac{ \dd \overline t_{2}}{\overline t_{2}} \right)-
    \log( t_{2}\overline t_{2})
    \left(\frac{\dd t_{1}}{t_{1}}-
      \frac{ \dd  \overline t_{1}}{\overline t_{1}} \right)
  \right].
  $$
\end{definition}

Observe that $W_{2}\in \mathcal{D}^{2}(D^{\ast}_{(\mathbb{P}^{1})^{2}},2)$,
where $D^{\ast}_{(\mathbb{P}^{1})^{2}}$ is the Dolbeault complex of
currents on $\mathbb{P}^{1}\times \mathbb{P}^{1}$. Moreover, we can write
\begin{equation}\label{eq:30}
  W_{2}=\left[\left(-\frac{1}{2}\log( t_{1}\overline t_{1})
\right  )
\bullet \left(-\frac{1}{2}\log( t_{2}\overline t_{2}) \right)\right],
\end{equation}
where $\bullet$ is the product in the Deligne complex (see definition \ref{def:17}).  

For $p=(x_{0},y_{0})\in \mathbb{P}^{1}$, $i=1,2$, let 
$\iota _{i,p}:\mathbb{P}^{1}\longrightarrow \mathbb{P}^{1}\times
\mathbb{P}^{1}$ be the inclusion given by 
\begin{align*}
  \iota _{1,p}(x:y)&=(x_{0}:y_{0})\times (x:y),\\
  \iota _{2,p}(x:y)&=(x:y)\times (x_{0}:y_{0}).
\end{align*}

\begin{proposition}\label{prop:27}
  We have the equality
  \begin{displaymath}
    \dd_{\mathcal{D}}W_{2}=(\iota _{1,(1:0)})_{\ast}W_{1}
    -(\iota _{1,(0:1)})_{\ast}W_{1}-
    (\iota _{2,(1:0)})_{\ast}W_{1}+
    (\iota _{2,(0:1)})_{\ast}W_{1}.
  \end{displaymath}
\end{proposition}
\begin{proof}
  This proposition follows easily from a residue computation. Formally
  we can interpret it as the Leibniz rule for the Deligne complex and
  equations \eqref{eq:31} and \eqref{eq:30}.
\end{proof}

\nnpar{The iterated Bott-Chern form.} 
\begin{definition} \label{def:6} Let $\overline {\mathcal{F}}$
be a hermitian exact square satisfying the condition of definition
\ref{def:10}.  
  The
  \emph{iterated Bott-Chern form} associated to
  $\overline{\mathcal{F}}$ is the differential form 
   given by  
  \begin{align*}
    \phi(\overline{\mathcal{F}})=
    W_{2}(\phi (\overline{\mathcal{F}}))
    =&\frac{1}{(4\pi i)^{2}}\int_{\mathbb{P}^{1}\times \mathbb{P}^{1}}
    \phi (\tr_{2}(\overline{\mathcal{F}}))\land
    \log( t_{1}\overline t_{1})
    \left(\frac{\dd  t_{2}}{t_{2}}-
      \frac{ \dd  \overline t_{2}}{\overline t_{2}} \right)-\\
    &\frac{1}{(4\pi i)^{2}}\int_{\mathbb{P}^{1}\times \mathbb{P}^{1}}
    \phi (\tr_{2}(\overline{\mathcal{F}}))\land
    \log( t_{2}\overline t_{2})
    \left(\frac{\dd  t_{1}}{t_{1}}-
      \frac{ \dd  \overline t_{1}}{\overline t_{1}} \right).
  \end{align*}
\end{definition}

When $\overline {\mathcal{F}}$ does not satisfy the condition of
definition \ref{def:10} we proceed as
follows. Let $\lambda ^{k}_{i}\overline {\mathcal{F}}$, $i=1,2$,
$k=1,2$, be the hermitian 
exact square determined by
\begin{displaymath}
  \partial_{i}^{j} (\lambda ^{k}_{i}\overline {\mathcal{F}})=
  \lambda ^{k}(\partial_{i}^{j}\overline {\mathcal{F}}) \qquad(j=-1,0,1);
\end{displaymath}
here $\lambda^k( \cdot)$ is as in definition \ref{def:3}.
 
\begin{definition} \label{def:7}
  Let $\overline {\mathcal{F}}$ be a hermitian exact square. Then, the
  \emph{iterated Bott-Chern form} associated to
  $\overline{\mathcal{F}}$ is the differential form 
   given by  
   \begin{displaymath}
     \phi (\overline {\mathcal{F}})= 
     \phi (\lambda _{1}^{1}\lambda _{2}^{1}\overline {\mathcal{F}})+
     \phi (\lambda _{1}^{1}\lambda _{2}^{2}\overline {\mathcal{F}})+
     \phi (\lambda _{1}^{2}\lambda _{2}^{1}\overline {\mathcal{F}})+
     \phi (\lambda _{1}^{2}\lambda _{2}^{2}\overline {\mathcal{F}}).
   \end{displaymath}
\end{definition}

As in the case of exact sequences, if $\overline {\mathcal{F}}$
satisfies the condition of definition \ref{def:10}, then the iterated
Bott-Chern forms obtained by means of definition \ref{def:6} and definition
\ref{def:7} agree.

\begin{theorem} \label{thm:5}
  The second iterated Bott-Chern form satisfies
  \begin{displaymath}
    \dd_{\mathcal{D}}\phi (\overline{\mathcal{F}})=
    \phi (\partial_{1}^{-1} 
    \overline{\mathcal{F}}\oplus \partial_{1}^{1} 
    \overline{\mathcal{F}})-\phi (\partial_{1}^{0} 
    \overline{\mathcal{F}})-\phi (\partial_{2}^{-1} 
    \overline{\mathcal{F}}\oplus \partial_{2}^{1} 
    \overline{\mathcal{F}})+\phi (\partial_{2}^{0} 
    \overline{\mathcal{F}}).
  \end{displaymath}
\end{theorem}
\begin{proof}
  This follows from \eqref{eq:29} and proposition \ref{prop:27}.
\end{proof}

\nnpar{The case of three different metrics.}  Let $X$ be a complex
  manifold, $E$ a holomorphic vector bundle on 
  $X$ and $h$, 
  $h'$ and $h''$ smooth hermitian metrics on $E$. We will denote by
  $\overline {\mathcal{F}}(E,h,h',h'')$ the hermitian exact square  
  \begin{displaymath}
    \begin{CD}
      0@>>> 0 @>>> 0\\
      @VVV @VVV @VVV\\
      0@>>> (E,h) @>>> (E,h'')\\
      @VVV @VVV @VVV\\
      0@>>> (E,h')@>>> (E,h'')
    \end{CD}
  \end{displaymath}
  where the faces $\partial^{j}_{1}$ are the rows and the faces
  $\partial^{j}_{2}$ are the columns. As a shorthand we will denote
  the hermitian vector bundle $\tr_{2}(\overline
  {\mathcal{F}}(E,h,h',h''))$ by $\tr_{2}(E,h,h',h'')$, or simply by
  $\tr_{2}(h,h',h'')$, if $E$ is understood.

\begin{definition}
  The
  \emph{iterated Bott-Chern form} associated to the metrics $h$, $h'$
  and $h''$ is the differential form
  given by  
  \begin{displaymath}
    \phi(E,h,h',h'')= \phi (\overline {\mathcal{F}}(E,h,h',h'')).
  \end{displaymath}
\end{definition}

\begin{proposition}
  The iterated Bott-Chern form satisfies
  \begin{displaymath}
    \dd_{\mathcal{D}}\phi (E,h,h',h'')=
    \phi(E,h,h')+\phi(E,h',h'')
    +\phi(E,h'',h).
  \end{displaymath}
\end{proposition}
\begin{proof}
  By theorem \ref{thm:5}
  \begin{displaymath}
    \dd_{\mathcal{D}} \phi (E,h,h',h'')=
        \phi(E,h',h'')-\phi(E,h,h'')
        -\phi(E,h'',h'')
    +\phi(E,h,h').
  \end{displaymath}
  but a direct computation shows that $\phi(E,h'',h'')=0$
  and that
  $\phi(E,h,h'')=-\phi(E,h'',h)$, which implies
  the result.
\end{proof}

\subsection{Chern forms for singular hermitian  metrics}
\label{sec:good-metrics}

There are various successful concepts of singular metrics in Arithmetic
and Diophantine Geometry \cite{Bost:Lfg}, \cite{Faltings:EaVZ},
\cite{Moriwaki:ipacdGc} and \cite{Mumford:Hptncc}.  For our purposes
the most important are: Faltings notion of a metric with logarithmic
singularities along a divisors with normal crossings
\cite{Faltings:EaVZ} and Mumford notion of a good metric
\cite{Mumford:Hptncc}. Both concept have in common nature that
automorphic vector bundles (equipped with the natural metric) have
the required local behavior. And indeed, the application to
automorphic bundles was the driving motivation to establish these
definitions.  For our purposes we will need a more precise description
of the kind of metrics that appear when studying automorphic vector
bundles.

\nnpar{Faltings' logarithmic singular metric.}  Let $X$ be a complex
manifold and let $D$ be a normal crossing divisor.  Put
$U=X\setminus D$, and let
$j:U\longrightarrow X$ be the inclusion. Let $L$ be a line bundle on $X$
and $L_{0}$ the restriction to $U$. A smooth metric $h$ on $L_{0}$ is
said to have logarithmic singularities along $D$ if, for any
coordinate open subset $V$ adapted to $D$ and every non vanishing local
section $s$, there exists a number $N\in \mathbb{N}$ such that
\begin{align}
\max\{ h(s ) , h^{-1}(s)  \} \prec \left| \min_{j =1,..,k}
  \{ \log|r_j| \}\right|^N. 
\end{align} 
 
Observe that  this definition does not give any  information on the
behavior of the Chern form associated to the metric.  

\nnpar{Good metrics in the sense of Mumford.} We recall the notion of good
metric in the sense of Mumford \cite{Mumford:Hptncc}.

\begin{definition}
  Let $E$ be a rank $n$ vector bundle on $X$ and $E_{0}$ the restriction to
  $U$. A smooth metric $h$ on $E_{0}$ is said to be \emph{good} on $X$ if
  for all 
  $x\in \DD$ there is a neighborhood $V$ adapted $D$ and a
  holomorphic frame
  $\xi=\{e_{1},\dots,e_{n}\}$ such that, writing 
  $h(\xi )_{ij}=h(e_{i},e_{j})$, 
  \begin{enumerate}
  \item $|h(\xi )_{ij}|,\det(h)^{-1}\prec
    \left(\prod_{i=1}^{k}\log(r_{i})\right)^{N}$ for some
    $N\in \mathbb{N}$.
  \item The $1$-forms $(\partial h(\xi )\cdot h(\xi )^{-1})_{ij}$ are good.
  \end{enumerate}
  A vector bundle provided with a good hermitian metric will be called
  a \emph{good hermitian vector bundle.}
\end{definition}

\begin{lemma} \label{lemm:good-preloglog}
  If $(E,h)$ is a good hermitian vector bundle, then the $1$-forms
  $(\partial h(\xi )\cdot h(\xi )^{-1})_{ij}$ are pre-log-log forms. 
\end{lemma}
\begin{proof}
  Since a differential form with Poincar\'e growth has log-log growth
  (see \cite{BurgosKramerKuehn:cacg} \S7.1) we have that $(\partial
  h(\xi )\cdot h(\xi )^{-1})_{ij}$ and $\dd (\partial h(\xi )\cdot h(\xi
  )^{-1})_{ij}$ have log-log growth. Since the condition of having
  log-log growth is bihomogeneous and $(\partial h(\xi )\cdot h(\xi
  )^{-1})_{ij}$ has pure bidegree $(1,0)$ we have that $\partial(\partial
  h(\xi )\cdot h(\xi )^{-1})_{ij}$ and $\bar \partial(\partial h(\xi
  )\cdot h(\xi
  )^{-1})_{ij}$ have log-log growth. Finally, since 
  \begin{displaymath}
    \partial(\partial h(\xi )\cdot h(\xi)^{-1})=\partial h(\xi
    )\cdot h(\xi)^{-1} \wedge \partial h(\xi )\cdot h(\xi)^{-1}
  \end{displaymath}
  the form $\partial \bar \partial(\partial h(\xi )\cdot h(\xi
  )^{-1})_{i,j}$ also has log-log growth.
\end{proof}

A
 fundamental property of the concept of good  metrics is the following
 result of Mumford \cite{Mumford:Hptncc}.

\begin{proposition} Let $X$, $\DD$ and $U$ be as before.
 \begin{enumerate}
\item Let $(E_0,h)$ be a vector bundle over $U$. Then it has at most
  one extension to a vector bundle $E$ to $X$ such that $h$ is
  good 
  along $\DD$.
\item If $(E,h)$ is a good hermitian vector bundle, then, for
  any power 
  series $\phi $, the Chern form $\phi(E,h)$ is good.
  Moreover, 
  its associated 
  current  $\left[\phi^{(k)}(E,h)\right]_X$  
  represents the Chern class $\phi(E)$ of $E$.\hfill $\square$  
\end{enumerate}
\end{proposition}

\nnpar{Good Metrics of infinite order.} Note that with the concept of
good metric
we have control on the local behavior of the Chern forms and of the
cohomology class represented by its associated currents. As we will see
later, we can also control the local behavior of the Bott-Chern forms.
In order to
have control on the cohomology classes
represented by the Chern forms we need a slightly stronger
definition, that is the analogue of our definition \ref{def:15} of
Poincar\'e singular forms.  

\begin{definition}\label{def:4} 
  Let $X$, $\DD$ and $U$ be as before.
  Let $E$ be a rank $n$ vector bundle on $X$ and let $E_{0}$ be the
  restriction of $E$ to $U$. A smooth metric on $E_{0}$ is said to be
  \emph{good of infinite order (along $\DD$)} if for every $x\in \DD$,
  there exist a 
  trivializing open coordinate neighborhood $V$ adapted to $\DD$ with
  holomorphic frame $\xi =\{e_{1},\dots,e_{n}\}$, such that,
  writing $h(\xi)_{ij}=h(e_{i},e_{j})$, then
  \begin{enumerate}
  \item The functions $h(\xi)_{ij},\ \det(h(\xi))^{-1}$ belong to 
    $ \Gamma (V,\mathscr{E}^{0}_{X}\left<\DD\right>)$,
  \item The $1$-forms $(\partial h(\xi)\cdot h(\xi)^{-1})_{ij}$ are
  Poincar\'e singular.
  \end{enumerate}
  A vector bundle provided with a good hermitian metric of infinite
  order will be 
  called a \emph{$\infty$-good hermitian vector bundle}.  
\end{definition}

\nnpar{Log-singular hermitian metrics.} Although the hermitian
metrics we are interested in, the automorphic hermitian metrics, are
$\infty$-good, we will consider a slightly bigger set of singular metrics, the
log-singular metrics, for which we will be able to define arithmetic
characteristic classes.

\begin{definition}\label{def:14} 
  Let $X$, $\DD$ and $U$ be as before.
  Let $E$ be a rank $n$ vector bundle on $X$ and let $E_{0}$ be the
  restriction of $E$ to $U$. A smooth metric on $E_{0}$ is said to be
  \emph{log-singular (along $\DD$)} if for every $x\in \DD$, there exist a
  trivializing open coordinate neighborhood $V$ adapted to $\DD$ with
  holomorphic frame $\xi =\{e_{1},\dots,e_{n}\}$, such that,
  writing $h(\xi)_{ij}=h(e_{i},e_{j})$, then
  \begin{enumerate}
  \item The functions $h(\xi)_{ij},\ \det(h(\xi))^{-1}$ belong to 
    $ \Gamma (V,\mathscr{E}^{0}_{X}\left<\DD\right>)$,
  \item The $1$-forms $(\partial h(\xi)\cdot h(\xi)^{-1})_{ij}$ belong to 
    $ \Gamma (V,\mathscr{E}^{1,0}_{X}\left<\left<\DD\right>\right>)$.
  \end{enumerate}
  A vector bundle provided with a log-singular hermitian metric will be
  called a \emph{log-singular hermitian vector bundle}.  
\end{definition}

Note that if a smooth metric on $E_{0}$ is log-singular then the
conditions of definition \ref{def:14} are satisfied in every
holomorphic frame in every trivializing open coordinate
neighborhood $V$ adapted to $\DD$.

\begin{remark}
By the very definition of log-singular metrics, the Chern forms $\phi (E,h)$
belong to the group $\oplus_{k}\mathcal{D}^{2k}(E_{X}\left<\left<
    \DD\right>\right>,k)$, if $(E,h)$ is a
log-singular hermitian vector bundle. Moreover,  as we will see in
proposition \ref{prop:2}, 
the form $\phi (E,h)$ represents the  Chern class $\phi (E)$ in
  $H^{\ast}_{\mathcal{D}}(X,\mathbb{R}(\ast))$.
\end{remark}

\nnpar{Basic properties of log-singular hermitian metrics.} The following
properties are easily verified.
\begin{proposition}
  Let $X$, $\DD$ and $U$ be as before. Let $E$ and $F$ be vector
  bundles on $X$ and let $E_{0}$ and $F_{0}$ be their restriction to
  $U$. Let $h_{E}$ and $h_{F}$ be smooth hermitian metrics on $E_{0}$
  and $F_{0}$. Write $\overline E=(E,h_{E})$ and $\overline
  F=(F,h_{F})$. 
  \begin{enumerate}
  \item The hermitian vector bundle $\overline E
    \overset{\perp}{\oplus} \overline F$ is log-singular along $\DD$
    if and only if $\overline E$ and $\overline F$ are log-singular
    along $\DD$.
  \item If $\overline E$ and $\overline F$ are log-singular along $\DD$,
    then the tensor product $\overline E\otimes \overline F$, the
    exterior and symmetric powers $\Lambda ^{n}\overline E$,
    $S^{n}\overline E$, the dual bundle $\overline E^{\vee}$ and the
    bundle of homomorphisms $\Hom(\overline E, \overline F)$, with
    their induced metrics, are log-singular along $\DD$.
  \end{enumerate}
  \hfill $\square$
\end{proposition}

\begin{remark} \label{rem:1}  
  Note however that the condition of being log-singular is not stable
  under taking general quotients and subbundles. That is, if $(E,h)$
  is a hermitian vector bundle, log-singular along a normal crossing
  divisor $\DD$, and $E'$ is a subbundle or a quotient bundle, then
  the induced metric on $E'$ does not need to be log-singular along
  $\DD$. For instance, let $X=\mathbb{A}^{2}$ with coordinates
  $(t,z)$. Let $E=\mathcal{O}_{X}\oplus \mathcal{O}_{X}$ be the
  trivial rank two vector bundle with hermitian metric given, in the
  frame $\{e_{1},e_{2}\}$, by the matrix
  \begin{align}
    H=\left(
      \begin{matrix}
        (\log (1/|z|))^{-1} & 0\\
        0 & 1
      \end{matrix}
    \right).
  \end{align}
  This hermitian metric is log-singular along the divisor
  $\DD=\{z=0\}$. But the 
  sub-bundle generated by the section $e_{1}+te_{2}$ with the induced
  metric does not satisfy the second condition of definition
  \ref{def:14}. Namely, let 
  $h(t,z)=\|e_{1}+te_{2}\|^{2}$. Then  
  \begin{align*}
    h(t,z)&=t\bar t+(\log (1/|z|))^{-1},\\
    \partial h/h&=\frac{\bar t \,dt}{t\bar t+(\log (1/|z|))^{-1}}+
    \frac{dz}{z(\log (1/|z|))^{2}(t\bar t+(\log (1/|z|))^{-1})}. 
  \end{align*}
  But the function $\bar t/(t\bar t+(\log (1/|z|))^{-1})$ is not log-log
  along $\DD$, as can be seen by considering the set of points
  \begin{displaymath}
    t=\sqrt{(\log(1/|z|))^{-1}}.
  \end{displaymath}
  In this concrete case, the induced metric is not far from
  being log-singular:
  If $\widetilde X$ is the blow up
  of $X$ along the point $(0,0)$ and $\widetilde \DD$ is the pre-image of
  $\DD$, then the metric $h $ is log-singular along $\widetilde
  \DD$. See proposition \ref{prop:7} for a related example.  
\end{remark}

\begin{remark}\label{rem:5}
  The condition of being log-singular is also not stable under
  extensions. That is, let 
  \begin{displaymath}
    0\longrightarrow (E',h')\longrightarrow (E,h)
    \longrightarrow (E'',h'')\longrightarrow 0
  \end{displaymath}
  be a short exact sequence with $h'$ and $h''$ the hermitian metrics
  induced by $h$. If $h'$ and $h''$ are log-singular then $h$ does not
  need to be log-singular. 
\end{remark}

\nnpar{Functoriality of log-singular metrics.} The following result is a
direct consequence of the definition and the functoriality of log
forms and log-log forms.
\begin{proposition}
  Let $X$, $X'$ be complex manifolds and let $\DD$, $\DD'$ be normal
  crossing divisors of $X$ and $X'$ respectively. If
  $f:X'\longrightarrow X$ is a holomorphic map such that
  $f^{-1}(\DD)\subseteq \DD'$ and $(E,h)$ is a log-singular hermitian
  vector bundle on $X$, then $(f^{\ast}E,f^{\ast}h)$ is a log-singular
  hermitian vector bundle on $X'$.
  \hfill $\square$
\end{proposition}

\subsection{Bott-Chern forms for singular hermitian metrics}
\label{sec:bott-chern-forms-1}

\nnpar{ Bott-Chern forms for log-singular hermitian metrics.}  In order
to define 
characteristic classes of log-singular hermitian metrics with values in the
log-log arithmetic Chow groups, we have to
show that the Bott-Chern forms associated to a change of metric between a
smooth metric and a log-singular metric is a log-log form. By the
proof of the next theorem, it is clear that, even if we restrict ourselves to
$\infty$-good hermitian metrics, the Bott-Chern forms are not necessarily
Poincar\'e singular but log-log. Therefore the  log-log forms are an essential
ingredient of the theory and not only a technical addition to have
the Poincar\'e lemma.
 
\begin{theorem} \label{thm:3} Let $X$ be a complex manifold and let
  $\DD$ be a normal crossing divisor. Put $U=X\setminus \DD$. Let $E$
  be a vector bundle on $X$.
  \begin{enumerate}
  \item \label{itemthm3:1} If $h$ is a smooth hermitian metric on $E$
    and $h'$ is a smooth hermitian metric on $E|_{U}$, which is
    log-singular 
    along $\DD$, then the Bott-Chern form $\phi (E,h,h')$ belongs to
    the group $\bigoplus_{k}\mathcal{D}^{2k-1}
    (E_{X}\left<\left<\DD\right>\right>,k)$.
  \item \label{itemthm3:2} If $h$ and $h'$ are smooth hermitian
    metrics on $E$ and $h''$ is a smooth hermitian metric on $E|_{U}$,
    which is
    log-singular along $D$, the iterated Bott-Chern form
    $\phi(E,h,h',h'')$ belongs to the group
    $\bigoplus_{k}\mathcal{D}^{2k-2}
    (E_{X}\left<\left<\DD\right>\right>,k)$.
  \end{enumerate}
\end{theorem}
\begin{proof}
  Let $V$ be a trivializing coordinate subset adapted to $\DD$, with
  coordinates $(z_{1},\dots ,z_{d})$. Thus $D$ has equation
  $z_{1}\cdot\ldots \cdot z_{k}=0$ and we put $r_{i}=|z_{i}|$. We may
  assume also 
  that $V$ is contained in a compact subset of $X$. Let $\xi=\{e
  _{i}\}$ be a local holomorphic frame for $E$.  Let $g$ be the
  hermitian metric of $\tr_{1}(E,h,h')$. Since the vector bundle
  $\tr_{1}(E,h,h')$ is isomorphic to $p_{1}^{\ast}E$, the holomorphic
  frame $\xi$ induces a holomorphic frame (also denoted $\xi$) of
  $\tr_{1}(E,h,h')$.
  
  In the rest of the proof the frame $\xi$ will be fixed. Therefore we
  drop it from the notation and we write
  \begin{displaymath}
    H=h(\xi),\qquad H_{ij}=h(\xi)_{ij}=h(e_{i},e_{j}).
  \end{displaymath}
  We use the same notation for the metrics $h'$ and $g$.
  
  Let $(x:y)$ be homogeneous coordinates of $\mathbb{P}^{1}$. Write
  $t=x/y$.  We decompose $\mathbb{P}^{1}$ into two closed sets
  \begin{displaymath}
    \mathbb{P}^{1}_{+}=\{(x:y)\in \mathbb{P}^{1}\mid |x|\ge |y|\}\ \text{
      and }\  
    \mathbb{P}^{1}_{-}=\{(x:y)\in \mathbb{P}^{1}\mid |x|\le |y|\}.
  \end{displaymath}
  Then we write
  \begin{displaymath}
    \phi (h,h')=\phi _{+}(h,h')+
    \phi _{-}(h,h'),
  \end{displaymath}
  with
  \begin{equation}\label{eq:18}
    \phi _{\pm}(h,h')=\frac{-1}{4\pi i} \int_{\mathbb{P}_{\pm}^{1}}
    \phi (\tr_{1}(h,h'))\log( t\overline t).     
  \end{equation}
  
  We first show that the form $\phi _{-}(h,h')$ is log-log along
  $\DD$.  One technical difficulty that we have to solve at this point
  is that the differential form $\phi (\tr_{1}(h,h'))$ is not, in
  general, a log-log form along $\DD\times \mathbb{P}^{1}$, because
  the vector bundle $\tr_{1}(h,h')$ does not need to be log-singular
  along $\DD\times \mathbb{P}^{1}$. This is the reason we have to
  introduce a new class of singular functions.

  \begin{definition} \label{def:18} For any pair of subsets  $I
  ,J\subset \{1,\dots ,d\} $ 
    and integers $n,K\ge 0$, we say that a smooth complex function $f$
    on $(V\setminus \DD)\times \mathbb{P}^{1}_{-}$ has
    \emph{singularities of type $(n,\alpha, \beta) $ of order $K$} if
    there is an integer $N\ge 0$ such that, for any pair of
    multi-indices $\alpha , \beta \in \mathbb{Z}_{\ge 0}^{d}$ and
    integers $a,b\ge 0$ with $|\alpha +\beta |+a+b\le K$, it holds the
    estimate
    \begin{multline*}
      \left |\frac{\partial^{|\alpha |}}{\partial z^{\alpha }}
        \frac{\partial^{|\beta |}}{\partial z^{\beta }}
        \frac{\partial^{a}}{\partial t} \frac{\partial^{b}}{\partial
          \bar t} f(z_{1},\dots ,z_{d},t)\right | \prec
      \left(\frac{1}{|t|+(\prod_{i=1}^{k} \log( 1/r_{i}))^{-N}}
      \right)^{n+a+b}\cdot\\
      \cdot\frac{ \left |\prod_{i=1}^{k} \log (\log(1/r_{i})) \right
        |^{N}}{r^{(\gamma ^{I}+\gamma ^{J}+\alpha+\beta )^{\le
            k}}(\log (1/r))^{(\gamma ^{I}+\gamma ^{J})^{\le k}}}.
    \end{multline*}     
    We say that $f$ has \emph{singularities of type $(n,\alpha, \beta)
    $ of infinite order}, if
    it has singularities of type $(n,\alpha, \beta)$ of order $K$ for
    all $K\ge 0$. 
  \end{definition}

  The singularities of the differential form $\phi (\tr_{1}(h,h'))$
  are controlled by the following result.
 
  \begin{lemma} \label{lemm:9}
    Let
    \begin{displaymath}
      \phi (\tr_{1}(h,h'))=\sum_{\substack{0\le a,b\le 1\\
          I,J}}
      f_{I ,J ,a,b}\dd z^{I}\land \dd\bar z^{J }\land
      \dd t^{a}\land \dd\bar t^{b} 
    \end{displaymath}
    be the decomposition in monomials of $\phi (\tr_{1}(h,h'))$ over
    $V\times \mathbb{P}_{-}^{1}$.  Then, the function $f_{I,J,a,b}$ has
    singularities of type $(a+b,I,J)$ of infinite order.
  \end{lemma}
  \begin{proof}
    On $V\times \mathbb{P}^{1}_{-}$, the  matrix of $g$ in the
    holomorphic frame $\xi$ is
    \begin{displaymath}
      G=\frac{1}{1+t\overline t}(H+t\overline t\, H').
    \end{displaymath}
    We write $G_{1}=H+t\overline t\, H'$.     The differential form
    $\phi (\tr_{1}(h,h'))$ is a polynomial in 
    the entries of the matrix $\bar \partial(\partial G
    G^{-1})$. Since
    \begin{displaymath}
      \partial G  G^{-1}=\frac{-\overline t\,\dd t}{1+t\overline t}\Id+ 
      \partial G_{1}  G_{1}^{-1}, 
    \end{displaymath}
    and the first summand of the right term is smooth, we are led to
    study the singularities of the matrices $\partial G_{1}
    G_{1}^{-1}$ and $\bar \partial (\partial G_{1} G_{1}^{-1})$. This
    will be done in the subsequent lemmas.
    
    We write $G_{2}=(H'{}^{-1}+t\overline t H^{-1})$. The following
    lemma is easy.
    \begin{lemma}\label{lemm:4}
      The matrices $H$, $H'$, $G_{1}$ and $G_{2}$ satisfy the rules
      \begin{alignat*}{2}
        (1)&\ H G_{1}^{-1}=G_{2}^{-1}H'{}^{-1},\qquad &
        \qquad (2)&\  G_{1}^{-1}H=H'{}^{-1}G_{2}^{-1},\\
        (3)&\ H'G_{1}^{-1}=G_{2}^{-1}H^{-1},\qquad & \qquad (4)&\ 
        G_{1}^{-1}H'=H^{-1}G_{2}^{-1}.
      \end{alignat*}
      \hfill $\square$
    \end{lemma}
    
    In order to bound the entries of $\partial G_{1}G_{1}^{-1}$ and
    the other matrices, we need the following estimates.
   
    \begin{lemma}\label{lemm:2}
      \begin{enumerate}
      \item The entries of the matrix $G_{1}^{-1}$ are bounded. In
        particular, they have singularities of type
        $(0,\emptyset,\emptyset)$ of order $0$.
      \item \label{item:bound2} The entries of the matrix $G_{2}^{-1}$
        have singularities of type $(2,\emptyset,\emptyset)$ of order
        $0$. Therefore, the entries of the matrices $t\, G_{2}^{-1}$
        and $\overline t\, G_{2}^{-1}$ have singularities of type
        $(1,\emptyset,\emptyset)$ of order $0$ and the entries of the
        matrix $t\overline t\, G_{2}^{-1}$ are bounded.
      \end{enumerate}
    \end{lemma}
    \begin{proof}
      Let $H=T^{+}T$ be the Cholesky decomposition of $H$, where
      $(\bullet)^{+}$ denotes conjugate-transpose. Since $H$ is
      smooth, the same is true for $T$. We can write
      \begin{displaymath}
        G_{1}^{-1}=T^{-1}(\Id+t\overline t
        (T^{-1})^{+}H'T^{-1})^{-1}
        (T^{-1})^{+}.
      \end{displaymath}
      But for any symmetric definite positive matrix $A$, the entries
      of $(\Id+A)^{-1}$ have absolute value less that one. Therefore
      the entries of the matrix $G_{1}^{-1}$ are bounded. This proves
      the first statement.
      
      To prove the second statement we write
      \begin{displaymath}
        G^{-1}_{2}=T^{+}( 
        TH'{}^{-1}T^{+}+t\overline t\,\Id)^{-1}T.
      \end{displaymath}
      By the first condition of log-singular metric, we can decompose
      \begin{displaymath}
        TH'{}^{-1}T^{+}=U^{+}DU,
      \end{displaymath}
      with $U$ unitary and $D$ diagonal with all the diagonal elements
      bounded from above by $(\prod_{i=1}^{k} \log( 1/r_{i}))^{N}$ and
      bounded from below by $(\prod_{i=1}^{k} \log( 1/r_{i}))^{-N}$ for
      some integer $N$. Then
      \begin{displaymath}
        G^{-1}_{2}=(UT)^{+}( 
        D+t\overline t\,\Id)^{-1}(UT).
      \end{displaymath}
      Now the lemma follows from the fact that the norm of any entry of a
      unitary matrix is less or equal than one.
    \end{proof}
   
    The remainder of the proof of lemma \ref{lemm:9} is based on lemma
    \ref{lemm:2}.

    \begin{lemma} \label{lemm:16}
      Let $\sum \psi _{I ,J ,a,b}\dd z^{I}\land \dd\bar z^{J }\land
      \dd t^{a}\land 
      \dd\bar t^{b} $ be the decomposition into monomials of an entry of
      any of the matrices $\partial G_{1}G_{1}^{-1}$, $\bar \partial (\partial
      G_{1}G_{1}^{-1})$, $\partial(\partial G_{1}G_{1}^{-1})$ and
      $\partial \bar \partial (\partial
      G_{1}G_{1}^{-1})$. Then
      $\psi _{I,J,a,b}$ has singularities of type $(a+b,I,J)$ of order
      $0$. 
    \end{lemma}
    \begin{proof}
      We start with the entries of $\partial G_{1}G_{1}^{-1}$. Using
      lemma \ref{lemm:4}, we have
      \begin{align}
        \partial G_{1}G_{1}^{-1}&=\partial HG_{1}^{-1}+\bar t\dd t
        H'G_{1}^{-1}+ t\bar t \partial H' G_{1}^{-1}\notag \\
        &=\partial HG_{1}^{-1}+(\bar t\dd t
        + t\bar t \partial H' H'{}^{-1}) G_{2}^{-1}H^{-1}.\label{eq:16}
      \end{align}
      Therefore, the bound of the entries of $\partial
      G_{1}G_{1}^{-1}$ follows from lemma \ref{lemm:2} and the fact
      that $h'$ is log-singular.

      The bound of the entries of $\partial(\partial G_{1}G_{1}^{-1})$
      follows from the previous case and the formula
      \begin{equation}\label{eq:17}
        \partial(\partial G_{1}G_{1}^{-1})=
        \partial G_{1}G_{1}^{-1}\land \partial G_{1}G_{1}^{-1}.
      \end{equation}

      Before bounding $\bar \partial (\partial G_{1}G_{1}^{-1})$, we
      compute 
      \begin{displaymath}
        \bar \partial G_{1}^{-1}=-G_{1}^{-1}\bar \partial G_{1} G_{1}^{-1}=
        -(\partial G_{1} G_{1}^{-1})^{+}G_{1}^{-1}
      \end{displaymath}
      and 
      \begin{align*}
        \bar \partial G_{2}^{-1}&=-G_{2}^{-1}\bar \partial G_{2}
        G_{2}^{-1}\\
        &=-G_{2}^{-1}(\bar \partial H'{}^{-1}+t\dd \bar t H^{-1}+t\bar
        t\bar \partial H^{-1}) G_{2}^{-1}\\
        &=G_{2}^{-1}H'{}^{-1}\bar \partial H'
        H'{}^{-1}G_{2}^{-1}-G_{2}^{-1}(t\dd
        \bar t H^{-1}+t\bar 
        t\bar \partial H^{-1}) G_{2}^{-1}\\
        &=G_{2}^{-1}(\partial H'
        H'{}^{-1})^{+}G_{1}^{-1} H-G_{2}^{-1}(t\dd
        \bar t H^{-1}+t\bar 
        t\bar \partial H^{-1}) G_{2}^{-1}.
      \end{align*}
      Therefore,

      \begin{align}\label{eq:15}
        \bar \partial (\partial G_{1}G_{1}^{-1})=&\bar
        \partial\partial H G_{1}^{-1}+\partial H\land (\partial G_{1}
        G_{1}^{-1})^{+}G_{1}^{-1}\notag\\
        &+ \bar \partial (\bar t\dd t
        + t\bar t \partial H' H'{}^{-1}) G_{2}^{-1}H^{-1}\notag\\
        &-(\bar t\dd t
        + t\bar t \partial H' H'{}^{-1}) G_{2}^{-1}\land (\partial H'
        H'{}^{-1})^{+}G_{1}^{-1}\notag\\
        &+(\bar t\dd t
        + t\bar t \partial H' H'{}^{-1}) G_{2}^{-1}\land (t\dd
        \bar t H^{-1}+t\bar 
        t\bar \partial H^{-1}) G_{2}^{-1} H^{-1}\notag\\
        &-(\bar t\dd t
        + t\bar t \partial H' H'{}^{-1}) G_{2}^{-1}\land \bar \partial H^{-1}
      \end{align}

      Thus the bound of the entries of $\bar \partial (\partial
      G_{1}G_{1}^{-1})$ follows again by lemma \ref{lemm:2} and the
      assumptions on $H$ and $H'$.

      Finally, the case of $\partial \bar \partial (\partial
      G_{1}G_{1}^{-1})$ follows from the formula 
      \begin{equation}\label{eq:20}
        \partial \bar \partial (\partial G_{1}G_{1}^{-1})=
        -\bar \partial (\partial G_{1}G_{1})\land \partial G_{1}
        G_{1}^{-1}
        +
        \partial G_{1}
        G_{1}^{-1} \land 
        \bar \partial (\partial G_{1}G_{1})
      \end{equation}
    \end{proof}

    As a direct consequence of the previous lemma, we obtain that the
    functions $f_{I,J,a,b}$ of lemma \ref{lemm:9} have singularities
    of type $(a+b,I,J)$ of order $0$. But we have to show that they
    have singularities of type $(a+b,I,J)$ of infinite order. Thus we
    have to bound all of their derivatives. As before it is enough to
    bound the derivatives of the components of the entries of the
    matrix $\bar \partial (\partial G_{1}G_{1}^{-1})$. By the formulas
    \eqref{eq:16} and \eqref{eq:15}, it is enough to bound the
    derivatives of the entries of the matrices $G_{1}^{-1}$ and
    $G_{2}^{-1}$.  The idea to accomplish this task is to use induction,
    because the derivatives of these matrices can be written in terms
    of the same matrices and the derivatives of $H$ and $H'$, that we
    control. The inductive step is provided by the next lemmas.
    
    \begin{lemma}\label{lemm:5} If the entries of the matrices
      $G_{1}^{-1}$ and $G^{-1}_{2}$ have singularities of type
      $(0,\emptyset,\emptyset)$ and $(2,\emptyset,\emptyset)$
      respectively  of order $K$, then, for every $i=1,\dots
      ,d$, the entries of $\frac{\partial}{\partial z_{i}} G_{1}
      G_{1}^{-1}$ have singularities of type $(0,\{i\},\emptyset)$ of
      order $K$ and the entries of $\frac{\partial}{\partial t} G_{1}
      G_{1}^{-1}$ have singularities of type $(1,\emptyset,\emptyset)$ of
      order $K$ . 
    \end{lemma}
    \begin{proof}
      The result is a consequence of the formulas
      \begin{align}
        \frac{\partial}{\partial z_{i}} G_{1} G_{1}^{-1}&=
        \frac{\partial}{\partial z_{i}} H G_{1}^{-1}+ t\overline t
        \left(\frac{\partial}{\partial z_{i}} H' H'{}^{-1}\right)
        G_{2}^{-1} 
        H^{-1},\\
        \frac{\partial}{\partial t} G_{1} G_{1}^{-1}&= \bar t \dd t
        G_{2}^{-1} H^{-1},
      \end{align}
      which follow from equation \eqref{eq:16}.
    \end{proof}
    
    \begin{lemma} \label{lemm:6}
      If the entries of the matrix $G_{1}^{-1}$ have singularities of
      type $(0,\emptyset,\emptyset)$ of order $K$,  for all $i=1,\dots ,d$,
      the entries of the matrix $\frac{\partial}{\partial z_{i}}G_{1}
      G_{1}^{-1} $ have singularities of type  $(0,\{i\},\emptyset)$
      of order $K$
      and the entries of the matrix $\frac{\partial}{\partial t}G_{1}
      G_{1}^{-1} $ have singularities of type
      $(1,\emptyset,\emptyset)$ of order $K$, 
      then the entries of the matrix $G_{1}^{-1}$ have singularities
      of type $(0,\emptyset,\emptyset)$ of order $K+1$.
    \end{lemma}
    \begin{proof}
      The result follows from formulas
      \begin{alignat*}{2}
        \frac{\partial}{\partial z_{i}} G_{1}^{-1}&=- G_{1}^{-1}
        \left(\frac{\partial}{\partial z_{i}} G_{1}
        G_{1}^{-1}\right),&\qquad 
        \frac{\partial}{\partial \overline z_{i}} G_{1}^{-1}&=
        -\left(\frac{\partial}{\partial z_{i}} G_{1}
        G_{1}^{-1}\right)^{+}G_{1}^{-1},\\
        \frac{\partial}{\partial t} G_{1}^{-1}&=- G_{1}^{-1}
        \left(\frac{\partial}{\partial t} G_{1} G_{1}^{-1}\right),&\qquad
        \frac{\partial}{\partial \overline t} G_{1}^{-1}&=
        -\left(\frac{\partial}{\partial t} G_{1} G_{1}^{-1}\right)
        ^{+}G_{1}^{-1}.
      \end{alignat*}
    \end{proof}
    
    \begin{lemma}\label{lemm:7}
      If the entries of the matrices $G_{1}^{-1}$ and $G^{-1}_{2}$
      have singularities of type 
      $(0,\emptyset,\emptyset)$ and $(2,\emptyset,\emptyset)$
      respectively  of order $K$,
      then the entries of the matrix $G^{-1}_{2}$ have singularities of type
      $(2,\emptyset,\emptyset)$ of order $K+1$.
    \end{lemma}
    \begin{proof}
      This result is consequence of the equations
      \begin{align*}
        \frac{\partial}{\partial z_{i}} G_{2}^{-1}&= -t\overline
        t G_{2}^{-1} \frac{\partial}{\partial z_{i}} 
        H^{-1} G_{2}^{-1} +H G_{1}^{-1}
        \left(\frac{\partial}{\partial z_{i}}H' H'{}^{-1} 
        \right)G_{2}^{-1},\\
        \frac{\partial}{\partial \overline z_{i}}
        G_{2}^{-1}
        &=-t\overline t G_{2}^{-1}
        \frac{\partial}{\partial \overline z_{i}} H^{-1} G_{2}^{-1}
        +G_{2}^{-1}
        \left(\frac{\partial}{\partial z_{i}}H' H'{}^{-1} 
        \right)^{+}G_{1}^{-1}
        H\\
        \frac{\partial}{\partial t}
        G_{2}^{-1}&=-G_{2}^{-1}(\bar t \dd t H^{-1})G_{2}^{-1},\\
        \frac{\partial}{\partial \bar t}
        G_{2}^{-1}&=-G_{2}^{-1}(t \dd \bar t H^{-1})G_{2}^{-1}.
      \end{align*}
    \end{proof}
    
    Summing up the lemmas \ref{lemm:2}, \ref{lemm:16}, \ref{lemm:5},
    \ref{lemm:6} and 
    \ref{lemm:7}, and equations \eqref{eq:16}, \eqref{eq:17},
    \eqref{eq:15} and \eqref{eq:20}, we obtain
    \begin{lemma} \label{lemm:15}
      Let $\sum \psi _{I ,J ,a,b}\dd z^{I}\land \dd\bar z^{J }\land
      \dd t^{a}\land
      \dd\bar t^{b} $ be the decomposition into monomials of an entry of
      any of the matrices $\partial G_{1}G_{1}^{-1}$, $\bar \partial (\partial
      G_{1}G_{1}^{-1})$, $\partial(\partial G_{1}G_{1}^{-1})$ and
      $\partial \bar \partial (\partial
      G_{1}G_{1}^{-1})$. Then
      $\psi _{I,J,a,b}$ has singularities of type $(a+b,I,J)$ of
      infinite order. \hfill $\square$
    \end{lemma}
    
    \noindent
    \emph{End of proof of lemma \ref{lemm:9}.} This finishes the proof
    of lemma \ref{lemm:9}. 
  \end{proof}
  
  Once we have bound the components of $\phi (\tr_{1}(h,h'))$ over
  $V\times \mathbb{P}_{-}^{1}$, to bound the components of $\phi
  _{-}(h,h')$ be have to estimate the integral \eqref{eq:18}.

  \begin{lemma}\label{lemm:10}
    Let $0\le a \le 1$ be a real number. Then
    \begin{align*}
      \int_{0}^{1}\frac{\log (1/r)}{r+a}\dd r&\le
      1+\log (1/a) +\frac{1}{2}\log^{2}(1/a),\\
      \int_{0}^{1}\frac{r \log (1/r)}{(r+a)^{2}}\dd r&\le
      1+\log (1/a) +\frac{1}{2}\log^{2}(1/a).
    \end{align*}
  \end{lemma}
  \begin{proof} We have the following estimates
    \begin{align*}
      \int_{0}^{1}\frac{r \log (1/r)}{(r+a)^{2}}\dd r&\le
      \int_{0}^{1}\frac{\log (1/r)}{r+a}\dd r\\
      &= \int_{0}^{a}\frac{\log (1/r)}{r+a}\dd r+
      \int_{a}^{1}\frac{\log (1/r)}{r+a}\dd r\\
      &\le \int_{0}^{a}\frac{\log (1/r)}{a}\dd r
      +\int_{a}^{1}\frac{\log (1/r)}{r}\dd r\\
      &= \left.\frac{r\log (1/r)+r}{a}\right|_{0}^{a}-
      \left.\frac{1}{2}\log^{2} (1/r) \right|_{a}^{1}\\
      &=\log (1/a) +1 +\frac{1}{2}\log^{2}(a).
    \end{align*}
  \end{proof}
  
  We are now in situation to bound the components of $\phi _{-}(h,h')$.
  Let
  \begin{displaymath}
    \phi _{-}(h,h')=\sum_{I,J } g_{I,J}\dd z^{I
    }\land \dd\overline z^{J} 
  \end{displaymath}
  be the decomposition of $\phi _{-}(h,h')$ in monomials. Then using
  lemma \ref{lemm:9} and lemma \ref{lemm:10} we have
  \begin{align*}
    |g_{I,J}|&=\left |\frac{1}{4\pi i}
      \int_{\mathbb{P}^{1}_{-}} f_{I ,J,1,1
      }\log( t\overline t)\, \dd t\land \dd\overline t\right |\\
    &\prec \frac{ \left |\prod_{i=1}^{k} \log (\log(1/r_{i})) \right
        |^{N}}{r^{(\gamma ^{I}+\gamma ^{J})^{\le
            k}}(\log (1/r))^{(\gamma ^{I}+\gamma ^{J})^{\le k}}}\cdot \\
    &\qquad \qquad \cdot
      \int_{\mathbb{P}^{1}_{-}}\left(\frac{1}{|t|+(\prod_{i=1}^{k} 
        \log( 1/r_{i}))^{-N}} \right)^{2} \log( t\overline t)\, \dd t\land
    \dd\overline t \\
    &\prec \frac{ \left |\prod_{i=1}^{k} \log (\log(1/r_{i})) \right
        |^{N'}}{r^{(\gamma ^{I}+\gamma ^{J})^{\le
            k}}(\log (1/r))^{(\gamma ^{I}+\gamma ^{J})^{\le k}}}.
  \end{align*}
  The derivatives of $g_{I,J }$ are bounded in the same way
  using the theorem of derivation under the integral sign. The
  components of $\partial \phi _{-}(h,h')$ and 
  $\bar \partial \phi _{-}(h,h')$ and their derivatives are bounded in
  a similar way using that 
  \begin{displaymath}
    \partial \phi _{-}(h,h')=\frac{-1}{4\pi i} \int_{\mathbb{P}_{-}^{1}}
    \phi (\tr_{1}(h,h'))\land\frac{\dd t}{t} .     
  \end{displaymath}
  and 
  \begin{displaymath}
    \bar \partial \phi _{-}(h,h')=\frac{-1}{4\pi i} \int_{\mathbb{P}_{-}^{1}}
    \phi (\tr_{1}(h,h'))\frac{\dd \bar t}{\bar t} .     
  \end{displaymath}

  To bound the components of $\phi _{+}(h,h')$, $\partial \phi
  _{+}(h,h')$ and $\bar \partial \phi _{+}(h,h')$ and their
  derivatives, we will use the same 
  technique. Let $s=1/t$ be a local coordinate in
  $\mathbb{P}^{1}_{+}$. In these coordinates we have 
  \begin{displaymath}
    G=\frac{1}{1+s\bar s}(H'+s\bar s H). 
  \end{displaymath}
  We write
  \begin{displaymath}
    G_{3}= (H'+s\bar s H),\qquad G_{4}=(H^{-1}+s\bar s H'{}^{-1}).
  \end{displaymath}
  In this case, using the adequate variant of definition \ref{def:18},
  the analogue of lemma \ref{lemm:2} is
  \begin{lemma}
    \begin{enumerate}
    \item  The entries of the matrix $G_{3}^{-1}$
        have singularities of type $(2,\emptyset,\emptyset)$ of order
        $0$. Therefore the entries of the matrices $s\, G_{3}^{-1}$
        and $\overline s\, G_{3}^{-1}$ have singularities of type
        $(1,\emptyset,\emptyset)$ of order $0$ and the entries of the
        matrix $t\overline s\, G_{3}^{-1}$ are bounded.
    \item The entries of the matrix $G_{4}^{-1}$ are bounded. In
      particular they have singularities of type
      $(0,\emptyset,\emptyset)$ of order $0$.
      \end{enumerate}
      \hfill $\square$
  \end{lemma}
  Note that the bounds of $G^{-1}_{3}$ and $G^{-1}_{4}$ are not the
  same as the bounds of $G_{1}^{-1}$ $G_{2}^{-1}$, but they are
  switched. To bound the entries of $\partial G_{3}G_{3}^{-1}$ we use
  \begin{displaymath}
    \partial G_{3}G_{3}^{-1}=
    \partial H' H'{}^{-1} G_{4}^{-1}H^{-1}+\bar s \dd s H G_{3}^{-1}+
    s\bar s \partial H G_{3}^{-1}.
  \end{displaymath}
  We left to the reader the details that remain. 

  Finally, to bound
  $\partial \bar \partial \phi(h,h')$, we use equation
  \eqref{eq:19}. This completes the proof of the first statement.

  We prove now the second statement.
  By definition
  \begin{displaymath}
    \tr_{2}(h,h',h'')=\tr_{1}(0\longrightarrow
    \tr_{1}(h,h')\longrightarrow \tr_{1}(h'',h'')). 
  \end{displaymath}
  But $\tr_{1}(h,h')$ is a smooth hermitian vector bundle on $X\times
  \mathbb{P}^{1}$ and $\tr_{1}(h'',h'')$ is isometric to
  $p_{1}^{\ast}(E,h'')$ and, in consequence, log-singular along $\DD\times
  \mathbb{P}^{1}$. Therefore  we can apply lemma \ref{lemm:9} to $\phi
  (\tr_{2}(h,h',h''))$ and, by lemma \ref{lemm:10}, 
  the form $\phi (h,h',h'')$ has log-log growth of infinite order. To
  conclude that it is a log-log form we still have to control
  $\partial \phi (h,h',h'')$, $\bar \partial \phi (h,h',h'')$, and
  $\partial \bar \partial \phi (h,h',h'')$. A residue computation
  shows that 
  \begin{align*}
    \partial \phi (h,h',h'')=&\frac{1}{2}(\phi (h,h')+\phi
    (h',h'')+\phi (h'',h))\\
    &+\frac{2}{(4\pi i)^{2}}\int_{\mathbb{P}^{1}\times \mathbb{P}^{1}}
    \phi(\tr_{2}(h,h',h'')) \land \frac{\dd t_{1}}{t_{1}}\land
    \frac{\dd t_{2}}{t_{2}}\\ 
    &- \frac{1}{(4\pi i)^{2}}\int_{\mathbb{P}^{1}\times \mathbb{P}^{1}}
    \phi (\tr_{2}(h,h',h''))\land \left(\frac{\dd \bar t_{1}}{\bar
        t_{1}}\land \frac{\dd t_{2}}{t_{2}}+ \frac{\dd
        t_{1}}{t_{1}}\land \frac{\dd \bar t_{2}}{\bar t_{2}} \right)
  \end{align*}
  and
  \begin{align*}
    \bar \partial \phi (h,h',h'')=&\frac{1}{2}(\phi (h,h')+\phi
    (h',h'')+\phi (h'',h))\\
    &-\frac{2}{(4\pi i)^{2}}\int_{\mathbb{P}^{1}\times \mathbb{P}^{1}}
    \phi(\tr_{2}(h,h',h'')) \land \frac{\dd \bar t_{1}}{\bar t_{1}}\land
    \frac{\dd \bar t_{2}}{\bar t_{2}}\\ 
    &+ \frac{1}{(4\pi i)^{2}}\int_{\mathbb{P}^{1}\times \mathbb{P}^{1}}
    \phi (\tr_{2}(h,h',h''))\land \left(\frac{\dd \bar t_{1}}{\bar
        t_{1}}\land \frac{\dd t_{2}}{t_{2}}+ \frac{\dd
        t_{1}}{t_{1}}\land \frac{\dd \bar t_{2}}{\bar t_{2}} \right)
  \end{align*}
  Hence, again by lemma \ref{lemm:10}, the forms  $\partial \phi
  (h,h',h'')$, $\bar \partial \phi (h,h',h'')$, have log-log growth of
  infinite order. Finally, since 
  \begin{displaymath}
    \partial\bar \partial \phi (h,h',h'')=(\partial- \bar \partial)
    (\phi (h,h')+\phi
    (h',h'')+\phi (h'',h)),
  \end{displaymath}
  by the first statement, the form $\partial \bar \partial \phi
  (h,h',h'')$ also has log-log growth of infinite order and therefore
  $\phi (h,h',h'')$ is a log-log form.
\vskip 2mm

  \noindent
  \emph{End of proof of theorem \ref{thm:3}.} This finishes the proof
    of theorem \ref{thm:3}.
\end{proof}

\nnpar{Bott-Chern forms for good hermitian metrics.} All the theory we
have developed so far is also valid for good hermitian vector bundles
with the obvious 
changes. For instance, if the hermitian
vector bundle is good instead of log-singular, we obtain that the
Bott-Chern forms are pre-log-log instead of log-log. 

\begin{theorem} \label{thm:secondary-good}  
  Let $X$ be a complex manifold and let $\DD$ be a normal crossing
  divisor. Put $U=X\setminus \DD$. Let $E$ be a vector bundle on $X$.
  If $h$ and $h'$ are smooth hermitian metrics on $E$ and $h''$ is a
  smooth hermitian metric on $E|_{U}$, which is good along $D$,
  then the Bott-Chern form $\phi (E,h,h'')$ and the iterated
  Bott-Chern form $\phi(E,h,h',h'')$ are pre-log-log forms.
\end{theorem}
\begin{proof}
  Observe that lemma \ref{lemm:4}, lemma \ref{lemm:2} and lemma
  \ref{lemm:16} are true in the case of good hermitian metrics by
  lemma \ref{lemm:good-preloglog}, and these results are enough to
  prove that $\phi (E,h,h')$ and $\phi(E,h,h',h'')$ are pre-log-log
  forms by the same
  arguments as before.
\end{proof}

\nnpar{The singularities of the first transgression bundle.} With the
  notation of theorem \ref{thm:3}, observe that  
  the hermitian vector bundle $\tr_{1}(E,h,h')$ does not need to be
  log-singular along the divisor $\DD\times \mathbb{P}^{1}$ (see remark
  \ref{rem:1}). Nevertheless, as we will see in the following results,
  it is close to 
  be log-singular. For instance it is log-singular along $\DD\times
  \mathbb{P}^{1}\cup X\times \{(0:1),(1:0)\}$, or it can be made a
  log-singular hermitian vector bundle after some blow-ups. This second
  statement will be useful in the axiomatic characterization of
  Bott-Chern classes.

\begin{lemma} \label{lemm:11} Let $a,b$ be real numbers with $a>0$ and
  $b>e^{1/e}$. Then
  \begin{displaymath}
    \frac{\log( a)}{\log( b)}< 1+\frac{a}{b}.  
  \end{displaymath}
\end{lemma}
\begin{proof}
  If $b\ge a$ then the statement is obvious. If $a>b$ we write $a=cb$
  with $c>1$. Then the inequality of the lemma is equivalent to
  \begin{displaymath}
    \frac{\log( c)}{c}<\log( b). 
  \end{displaymath}
  But the function $\log( c)/c$ is a bounded function that has a maximum
  at $c=e$ with value $1/e$. Therefore the result is a consequence of
  the condition on $b$.
\end{proof}

\begin{corollary}\label{cor:1} With the notation of theorem \ref{thm:3},
  the hermitian vector bundle $\tr_{1}(E,h,h')$ is log-singular along the
  divisor $$\DD\times \mathbb{P}^{1}\cup X\times \{(0:1),(1:0)\}.$$
\end{corollary}
\begin{proof} The first condition of definition \ref{def:14} is easy to
  prove. We will prove the second condition.  The lemma \ref{lemm:11}
  implies that, for $a,b\gg 0$,  it holds
  \begin{equation}\label{eq:3}
    \frac{1}{1/a+1/b}< \frac{\log( b)}{(1/a)\log( a)}.  
  \end{equation}
  Applying this equation to $a=1/|t|$ and $b=(\prod_{i=1}^{k}
  \log( 1/r_{i}))^{N}$ we obtain
  \begin{displaymath}
    \frac{1}{|t|+(\prod_{i=1}^{k} \log( 1/r_{i}))^{-N}}< 
    \frac{\log((\prod_{i=1}^{k} \log
      (1/r_{i}))^{N})}{|t|\log (1/|t|)} \prec
    \frac{\sum_{i=1}^{k} \log (\log
      (1/r_{i}))}{|t|\log (1/|t|)}.
  \end{displaymath}
  Therefore lemma \ref{lemm:15} implies that, on $V\times
  \mathbb{P}^{1}_{-}$, the entries of $\partial G G^{-1}$ are log-log
  along $\DD\times \mathbb{P}^{1}\cup X\times \{(0:1),(1:0)\}$.
  The proof for the bound on $V\times \mathbb{P}^{1}_{+}$ is
  analogous.
\end{proof}

\begin{proposition} \label{prop:7}
  With the same hypothesis of theorem \ref{thm:3}, let
  $\DD=\DD_{1}\cup \dots \cup \DD_{n}$ be the decomposition of $D$ in
  smooth irreducible components. Let $\widetilde Z$ be the variety
  obtained from $X\times \mathbb{P}^{1}$ by blowing up $\DD_{1}\times
  {(1:0)}$ and then, successively, the strict transforms of
  $\DD_{2}\times {(1:0)},\dots ,\DD_{n}\times {(1:0)},\DD_{1}\times
  {(0:1)},\dots ,\DD_{n}\times {(0:1)}$. Let $\pi :\widetilde
  Z\longrightarrow X\times \mathbb{P}^{1}$ be the morphism induced by
  the blow-ups and let $C\subset \widetilde Z$ be the pre-image by $\pi
  $ of $\DD\times \mathbb{P}^{1}$. Then
  \begin{enumerate}
  \item $C$ is a normal crossing divisor.
  \item The closed immersions $i _{0},i _{\infty}:X\longrightarrow
    X\times \mathbb{P}^{1}$, given by
    \begin{displaymath}
      i _{0}(p)=(p,(0:1)),\qquad i _{\infty}(p)=(p,(1:0)),
    \end{displaymath}
    can be lifted to closed immersions $j _{0},j
    _{\infty}:X\longrightarrow \widetilde Z$.
  \item The hermitian vector bundle $\pi ^{\ast} \tr_{1}(E,h,h')$ is
    log-singular along the divisor $C$.
  \end{enumerate}
\end{proposition}
\begin{proof}
  The first statement is obvious and the second is a direct
  consequence of the universal property of the blow-up and the fact
  that the intersection of the center of every blow-up with the
  transform of $X\times (1:0)$ or $X\times (0:1)$ is either empty or a
  divisor.
  
  To prove the third statement we will use the same notations as in the
  proof of theorem \ref{thm:3}. Let $U$ be the subset of $V\times
  \mathbb{P}_{-}^{1}$, where $|t|<1/e^{e}$. For simplicity we assume
  that the components of $\DD$ that meet $V$ are $\DD_{1},\dots
  ,\DD_{k}$ and that the component $\DD_{i}$ has equation $z_{i}=0$.
  Then $U$, with coordinates $(z_{1},\dots ,z_{d},t)$, is a coordinate
  neighborhood adapted to $\DD\times \mathbb{P}^{1}$.  The open subset
  $\pi ^{-1}(U)$ can be covered by $k+1$ coordinate neighborhoods,
  denoted by $\widetilde U_{1}, \dots , \widetilde U_{k+1}$. The
  coordinates of these subsets, the expression of $\pi $ and the
  equation of $C$ in these coordinates are given in the following
  table:

  \begin{center}
  \begin{tabular}{c|c|c|c}
    Subset & Coordinates & $\pi $  \\
    \hline
    $U_{1}$ & $(u,x_{1},\dots ,x_{n})$ 
    & $
    \begin{aligned}
      t&=u\\
      z_{1}&=ux_{1}\\
      z_{i}&=x_{i}\ (i\not=1)
    \end{aligned} $& $ux_{1}\dots x_{k}=0$
    \\
    \hline
    $
    \begin{gathered}
      U_{j}\\ (1<j<k+1)
    \end{gathered}
    $  & 
    $(u,x_{1},\dots,x_{n})$ 
    & $
    \begin{aligned}
      t&=ux_{1}\dots x_{j-1}\\ z_{j}&=ux_{j}\\ z_{i}&=x_{i}\ 
      (i\not=j)
    \end{aligned}
    $
    & 
    $ux_{1}\dots x_{k}=0$
    \\
    \hline
    $U_{k+1}$ & 
    $(u,x_{1},\dots,x_{n})$ 
    & $
    \begin{aligned}
      t&=ux_{1}\dots x_{k}\\ z_{i}&=x_{i}\ (i=1,\dots ,d)
    \end{aligned}
    $& 
    $x_{1}\dots x_{k}=0$
  \end{tabular}
  \end{center} 
  Since, for $j=1,\dots ,k$, it holds 
  $$
  \pi ^{-1}(\DD\times \mathbb{P}^{1}\cup X\times \{(0:1),(1:0)\})\cap
  U_{j}= C\cap U_{j},
  $$
  then, by corollary \ref{cor:1} and the functoriality of log-singular
  metrics, we know that $\pi ^{\ast}\tr_{1}(E,h,h')|_{U_{j}}$ is
  log-singular. Hence we have only to prove that $\pi
  ^{\ast}\tr_{1}(E,h,h')|_{U_{k+1}}$ is log-singular. 
  The first
  condition of definition \ref{def:14} follows easily from the
  definition of the metric $g$. To prove the second condition of
  definition \ref{def:14} we can proceed in two ways.
  The first method is to derive this result directly from
  lemma \ref{lemm:15} applying the chain rule. But since we have to
  bound all derivatives this is a notational nightmare. The second
  method is to bound the derivatives inductively mimicking the proof of
  lemma \ref{lemm:10}. To this end, instead of lemma \ref{lemm:2}, we
  use the following one.
  \begin{lemma}
    \begin{enumerate}
          \item The entries of the matrix $\pi ^{\ast}
        G_{1}^{-1}|_{U_{k+1}}$ are bounded in every compact subset of
        $U_{k+1}$. In
        particular, they are
        $(\emptyset,\emptyset)$-log-log growth functions of order $0$
        (see definition \ref{def:16}). 
      \item  If $\psi $ is an entry of the matrix $
        G_{2}^{-1}$ then 
        \begin{displaymath}
          |(\pi ^{\ast}\psi|_{U_{k+1}}) (x_{1},\dots ,x_{d},u)|\prec
          \left| \prod_{i=1}^{k} \log(1/|x_{i}|)\right| ^{N}
        \end{displaymath}
        for some integer $N$.
        Therefore $\pi ^{\ast} (t\psi )$
        and $\pi^{\ast}(\overline t \psi )$ are bounded in any compact
        subset of $U_{k+1}$ and, for $i=1,\dots ,k$, the function
        \begin{displaymath}
          \prod_{j\not=i}|x_{j}|\pi ^{\ast}\psi  
        \end{displaymath}
        is a $(\{i\},\emptyset)$-log-log growth function of order 0.
      \end{enumerate}
      \hfill$\square$
  \end{lemma}
   We left to the reader the care to make explicit the analogues of lemmas
   \ref{lemm:5}, \ref{lemm:6} and \ref{lemm:7} in this case.
  
  The proof that it is also log-singular in the pre-image of an open subset of
  $\mathbb{P}^{1}_{+}$ is analogous. 
\end{proof}

\nnpar{Chern forms of log-singular hermitian bundles.}
\begin{proposition} \label{prop:2}
  Let $X$ be a complex projective manifold, $\DD$ a normal crossing
  divisor of $X$, $(E,h)$ a hermitian vector bundle log-singular along
  $\DD$.  Let $\phi $ be any symmetric power series.  Then the Chern
  form $\phi(E,h)$ represents the Chern class $\phi (E)$ in
  $H^{\ast}_{\mathcal{D}}(X,\mathbb{R}(\ast))$.
\end{proposition}
\begin{proof}
  By theorem \ref{thm:fq} and theorem \ref{thm:4} the inclusion
  $$\mathcal{D}^{\ast}(E_{X},\ast)\longrightarrow
  \mathcal{D}^{\ast}(E_{X}\left<\left<\DD\right>\right>,\ast)$$
  is a
  quasi-isomorphism. Moreover, if $h'$ is a smooth hermitian metric on
  $E$, then, in the complex
  $\mathcal{D}^{\ast}(E_{X}\left<\left<\DD\right>\right>,\ast)$, we
  have
  \begin{displaymath}
    \phi (E,h)-\phi (E,h')=\dd_{\mathcal{D}}\phi(E,h',h).
  \end{displaymath}
  Therefore both forms represent the same class.
\end{proof}

\nnpar{Bott-Chern classes.}
\begin{definition} \label{def:9}
  Let $X$ be a complex manifold and $\DD$ a normal crossing divisor.
  Let
  \begin{displaymath}
    \overline {{\mathcal{E}}}: 0\longrightarrow (E',h')\longrightarrow (E,h)
    \longrightarrow (E'',h'')\longrightarrow 0
  \end{displaymath}
  be an exact sequence of hermitian vector bundles log-singular along
  $\DD$. Let $h'_{s}$, $h_{s}$ and $h''_{s}$ be smooth hermitian
  metrics on $E'$, $E$ and $E''$ respectively. We denote by $\overline
  {\mathcal{E}}_{s}$ the corresponding exact sequence of smooth vector
  bundles. Let $\phi $ be a symmetric power series. Then the
  \emph{Bott-Chern class} associated to $\overline{\mathcal{E}}$ is
  the class represented by
  \begin{displaymath}
    \phi (\overline {\mathcal{E}}_{s})+\phi (E'\oplus E'',h'_{s}\oplus
    h''_{s},h'\oplus h')-\phi (E,h_{s},h),
  \end{displaymath}
  in the group
  \begin{displaymath}
  \bigoplus_{k}\widetilde{\mathcal{D}}^{2k-1}
  (E_{X}\left<\left<\DD\right>\right>,k)=
  \bigoplus_{k}\mathcal{D}^{2k-1}
  (E_{X}\left<\left<\DD\right>\right>,k)\left/
    \dd_{\mathcal{D}}\mathcal{D}^{2k-2}
  (E_{X}\left<\left<\DD\right>\right>,k)\right. .
  \end{displaymath}
  This class is denoted by $\widetilde \phi (\overline{\mathcal{E}})$.
\end{definition}
\begin{proposition}\label{prop:8}
  The Bott-Chern classes are well defined.
\end{proposition}
\begin{proof}
  The fact that the Bott-Chern forms belong to the group
  $$\bigoplus_{k}\mathcal{D}^{2k-1}
  (E_{X}\left<\left<\DD\right>\right>)$$
  is proved in theorem
  \ref{thm:3}.
  
  Let $h'_{sa}$, $h_{sa}$ and $h''_{sa}$ be another choice of smooth
  hermitian metrics and let $\overline {\mathcal{E}}_{sa}$ be the
  corresponding exact sequence. We denote by $\overline{\mathcal{C}}$
  the exact square of smooth hermitian vector bundles
  \begin{displaymath}
    0\longrightarrow 0 \longrightarrow \overline {\mathcal{E}}_{sa}
    \longrightarrow \overline {\mathcal{E}}_{s}\longrightarrow 0.
  \end{displaymath}
  Then
  \begin{multline*}
    \phi (\overline {\mathcal{E}}_{s})+\phi (E'\oplus E'',h'_{s}\oplus
    h''_{s},h'\oplus h')-\phi (E,h_{s},h)\\
    -\phi (\overline {\mathcal{E}}_{sa})-\phi (E'\oplus
    E'',h'_{sa}\oplus
    h''_{sa},h'\oplus h')+\phi (E,h_{sa},h)=\\
    \dd_{\mathcal{D}}\phi (\overline{\mathcal{C}})-
    \dd_{\mathcal{D}}\phi (E'\oplus E'',h'_{s}\oplus
    h''_{s},h'_{sa}\oplus h''_{sa},h'\oplus h'')
    +\dd_{\mathcal{D}}\phi (E,h_{s},h_{sa},h).
  \end{multline*}
  Therefore the Bott-Chern classes do not depend on the choice of the
  smooth metrics.
\end{proof}

\nnpar{Axiomatic characterization of Bott-Chern classes.}
\begin{theorem}\label{thm:18}
  The Bott-Chern classes satisfy the following properties. If $X$ is
  a complex manifold, $\DD$  is a normal crossing divisor and 
  \begin{displaymath}
    \overline {{\mathcal{E}}}: 0\longrightarrow (E',h')\longrightarrow (E,h)
    \longrightarrow (E'',h'')\longrightarrow 0
  \end{displaymath}
  is a short exact sequence of hermitian vector bundles, log-singular along
  $\DD$, then
  \begin{enumerate}
  \item $\dd_{\mathcal{D}} \widetilde {\phi }(\overline
    {\mathcal{E}})=\phi (E'\oplus E'',h'\oplus h'')-\phi (E,h)$.
  \item If $(E,h)=(E'\oplus E'',h'\oplus h'')$ then $\widetilde {\phi}
    (\overline {\mathcal{E}})=0$.
  \item If $X'$ is another complex manifold, $\DD'$ is a normal
    crossing divisor in $X'$ and    $f:X'\longrightarrow X$ is a
    holomorphic map such that 
    $f^{-1}(\DD)\subseteq \DD'$, then $\widetilde {\phi
    }(f^{\ast}\overline {\mathcal{E}})= f^{\ast}\widetilde {\phi
    }(\overline {\mathcal{E}})$.
  \item If $\overline {\mathcal{F}}$ is a hermitian exact square of
    vector bundles on $X$, log-singular along $\DD$, then
  \begin{displaymath}
    \widetilde {\phi} (\partial_{1}^{-1} 
    \overline{\mathcal{F}}\oplus \partial_{1}^{1} 
    \overline{\mathcal{F}})-\widetilde {\phi} (\partial_{1}^{0} 
    \overline{\mathcal{F}})-\widetilde {\phi} (\partial_{2}^{-1} 
    \overline{\mathcal{F}}\oplus \partial_{2}^{1} 
    \overline{\mathcal{F}})+\widetilde {\phi} (\partial_{2}^{0} 
    \overline{\mathcal{F}})=0.
  \end{displaymath}
  \end{enumerate}
  Moreover these properties determine the Bott-Chern classes.
\end{theorem}
\begin{proof} First we prove the unicity. By \cite{GilletSoule:vbhm} 1.3.2
  (see also \cite{Soule:lag} IV.3.1) the properties (1) to (3)
  characterize the Bott-Chern classes in the case $\DD=\emptyset$. By
  functoriality, the Bott-Chern classes are determined for short exact
  sequences when the three metrics are smooth. Let $E$ be a vector
  bundle, $h$ a smooth hermitian metric on $E$ and $h'$ a hermitian
  metric log-singular along $\DD$. The vector bundle $\widetilde E =
  \tr_{1}(E,h,h')$ over $X\times \mathbb{P}^{1}$ is isomorphic (as
  vector bundle) to $p_{1}^{\ast}E$. Let $h_{1}$ be the hermitian
  metric on $\widetilde E$ induced by $h$ and this
  isomorphism. Then $h_{1}$ is a 
  smooth hermitian metric.  Let $h_{2}$ be the metric of definition
  \ref{def:8}. Let $\pi :Z\longrightarrow X\times \mathbb{P}^{1}$ and
  $C$ be as in proposition \ref{prop:7}. By this proposition $\pi
  ^{\ast}(\widetilde E,h_{2})$ is log-singular along $C$. Therefore we
  can assume the existence of the Bott-Chern class $\widetilde \phi
  (\widetilde E,h_{1},h_{2})$. Write $\pi '=p_{1}\circ \pi$. We
  consider the integral
  \begin{displaymath}
    I=-\frac{1}{2\pi i} \int_{\pi '} -2\partial\overline \partial
    \widetilde \phi 
    (\widetilde E,h_{1},h_{2})  
    \pi ^{\ast}(\frac{1}{2} \log( t\overline t)).  
  \end{displaymath}
  By property (1),
  \begin{align*}
    I&=-\frac{1}{2\pi i} \int_{\pi '} \phi (\widetilde E,h_{2}) \pi
    ^{\ast}(\frac{1}{2} \log( t\overline t)) + \frac{1}{2\pi i}
    \int_{\pi '} \phi (\widetilde E,h_{1})
    \pi ^{\ast}(\frac{1}{2} \log( t\overline t))\\
    &=-\frac{1}{2\pi i} \int_{\mathbb{P}^{1}} \phi (\tr_{1}(E,h,h'))
    \frac{1}{2} \log( t\overline t),
  \end{align*}
  because the second integral vanishes. But using Stokes theorem and
  properties (2) and (3) an in \cite{GilletSoule:vbhm} 1.3.2 or
  \cite{Soule:lag} IV.3.1,
  \begin{align*}
    I&\sim j_{\infty}^{\ast}\widetilde \phi (\widetilde E,h_{1},h_{2})
    -j_{0}^{\ast}\widetilde \phi (\widetilde E,h_{1},h_{2})\\
    &=\widetilde \phi (E,h,h')-\widetilde \phi (E,h,h)\\
    &=\widetilde \phi (E,h,h'),
  \end{align*}
  where the symbol $\sim$ means equality up to the image of
  $\dd_{\mathcal{D}}$. Therefore the class $\widetilde \phi (E,h,h')$
  is also determined by properties (1) to (3). Finally, for an
  arbitrary exact sequence, $\overline {\mathcal{E}}$, of hermitian
  vector bundles log-singular along $\DD$, property (4) implies that
  $\widetilde \phi (\overline {\mathcal{E}})$ is given by definition
  \ref{def:9}.
  
  Next we prove the existence. By proposition \ref{prop:8} it only
  remains to show that the Bott-Chern classes defined by \ref{def:9}
  satisfy properties (1) to (4). Property (1) is known for smooth
  metrics. If $E$ is a vector bundle, $h$ is a smooth hermitian metric
  and $h'$ a hermitian metric log-singular along $\DD$ then, since two
  differential forms that agree in an open dense subset are equal, 
  by the smooth case,
  \begin{displaymath}
    \dd_{\mathcal{D}}\widetilde \phi (E,h,h')=
    \phi (E,h')-\phi (E,h).
  \end{displaymath}
  The general case follows from these two cases. Property (2) follows
  directly from the case of smooth metrics and definition \ref{def:9}.
  Property (3) is obvious from the functoriality of the definition. To
  prove property (4) we consider $\overline {\mathcal{F}}'$ an exact
  square with the same vector bundles as $\overline {\mathcal{F} }$
  but smooth metrics. Then, if we use the definition of Bott-Chern
  classes in the expression
  \begin{displaymath}
    \widetilde {\phi} (\partial_{1}^{-1} 
    \overline{\mathcal{F}}\oplus \partial_{1}^{1} 
    \overline{\mathcal{F}})-\widetilde {\phi} (\partial_{1}^{0} 
    \overline{\mathcal{F}})-\widetilde {\phi} (\partial_{2}^{-1} 
    \overline{\mathcal{F}}\oplus \partial_{2}^{1} 
    \overline{\mathcal{F}})+\widetilde {\phi} (\partial_{2}^{0} 
    \overline{\mathcal{F}})=0,
  \end{displaymath}
  all the change of metric terms appear twice with opposite sign.
  Therefore, this property follows from the smooth case.
\end{proof}

\nnpar{Real varieties.} The following result follows easily.
\begin{proposition}
  Let $X_{\mathbb{R}}=(X,F_{\infty})$ be a real variety, $D$ a normal
  crossing divisor on $X_{\mathbb{R}}$, $E$ a complex vector bundle
  defined over $\mathbb{R}$, $h$, $h'$ (resp. h'') smooth hermitian
  metrics (resp. log-singular hermitian metric) on $E$ invariant under complex
  conjugation. Then the forms $\phi (E,h'')$, $\widetilde \phi
  _{1}(E,h,h'')$ and $\widetilde \phi _{2}(E,h,h',h'')$ belong to the
  group $$\bigoplus_{k}\mathcal{D}^{2k-1}
  (E_{X_{\mathbb{R}}}\left<\left<\DD\right>\right>,k)=
  \bigoplus_{k}\mathcal{D}^{2k-1}
  (E_{X_{\mathbb{R}}}\left<\left<\DD\right>\right>,k)^{\sigma },$$
  where $\sigma $ is the involution that acts as complex conjugation
  on the space and on the coefficients.\hfill $\square$
\end{proposition}

\section{Characteristic classes and $K$-theory of
  log-singular hermitian vector bundles}
\label{sec:char-class}

The arithmetic intersection theory of Gillet and Soul\'e is complemented
by an arithmetic $K$-theory and a theory of characteristic classes. In
this section we will generalize both theories to cover the kind of
singular hermitian metrics that appear naturally when considering
(fully decomposed) automorphic
vector bundles.
If $E$ be a vector bundle over a quasi-projective complex manifold
$X$, then a hermitian metric $h$ on $E$ may have arbitrary singularities 
near the boundary of $X$. Therefore, the associated Chern forms will 
also have arbitrary singularities ``at infinity''. Thus, in order
to define  
arithmetic characteristic classes for this kind of hermitian vector 
bundles, we are led to use the complex $\mathcal{D}_{\as}$.

\subsection{Arithmetic Chern classes of log-singular hermitian vector bundles}
\label{sec:char-class-good}

\nnpar{Arithmetic Chow groups with coefficients.}   Let $A$ be an
arithmetic ring. Let $\widehat X=(X,\cc)$ be a
$\mathcal{D}_{\log}$-arithmetic variety over $A$. Let $B$ be a subring
of $\mathbb{R}$. We will define the arithmetic
Chow groups of $\widehat X$ with coefficients in $B$ using the same
method than \cite{Bost:Lfg}. We follow the notations of
\cite{BurgosKramerKuehn:cacg} \S4.2. 

For an integer $p\ge 0$, let ${\rm Z}^{p}(X)_{B}={\rm Z}^{p}(X)\otimes
B$ be the group of 
algebraic cycles of $X$ with coefficients in $B$. Then, 
the \emph{group of $p$-codimensional arithmetic cycles 
of $\widehat{X}=(X,\cc)$ with coefficients in $B$} is 
\begin{displaymath}
\widehat{{\rm Z}}^{p}_{B}({X},\cc)=\left\{(y,\mathfrak{g}_{y})\in
{\rm Z}^{p}_{B}(X)\oplus\widehat{H}^{2p}_{\cc,\mathcal{Z}^{p}}(X,p)
\,\Big|\,\cl(y)=\cl(\mathfrak{g}_{y})\right\}.
\end{displaymath}
Let $\rata^{p}_{B}(X,\cc)$ be the sub-$B$-module of $\widehat{{\rm
    Z}}^{p}_{B}({X},\cc)$ generated by $\rata^{p}(X,\cc)$. We
define \emph{the $p$-th arithmetic Chow group of $\widehat{X}=(X,\cc)$
  with coefficients in $B$} by
\begin{displaymath}
\cha^{p}_{B}(X,\cc)=\widehat{{\rm Z}}_{B}^{p}(X,\cc)\big/\rata_{B}^{p}
(X,\cc).
\end{displaymath}

There is a canonical morphism 
\begin{displaymath}
  \cha^{p}_{B}(X,\cc)\longrightarrow 
  \cha^{p}(X,\cc)\otimes B.
\end{displaymath}
For instance, if $B=\mathbb{Q}$, this morphism an isomorphism,
but in general, if $B=\mathbb{R}$, it is not an isomorphism.

\nnpar{The main theorem.} Let 
 $X$ be a regular scheme flat and
  quasi-projective over $\Spec(A)$. Let $\DD$ be a 
  normal crossing divisor on $X_{\mathbb{R}}$. Write $\underline
  X=(X,\DD)$. Then,
  $(X,\mathcal{D}_{\as,\underline X})$ is a
  quasi-projective  $\mathcal{D}_{\log}$-arithmetic
  variety over $A$. A log-singular hermitian vector
  bundle over $X$ is a vector bundle over $X$, together with a
  metric on $E_{\infty}$,  smooth over $X_{\infty}\setminus
  \DD_{\infty}$ and 
  log-singular along $\DD_{\infty}$, that is invariant under complex
  conjugation.

\begin{theorem} \label{thm:11}
Let $\phi \in
  B[[T_{1},\dots ,T_{n}]]$ be a symmetric power series with
  coefficients in a subring $B$ of $\mathbb{R}$.
  Then there is a
  unique way to attach, to every log-singular hermitian vector bundle
  $\overline E=(E,h)$ of rank
  $n$ over a pair $\underline X=(X,\DD)$, a characteristic class
  \begin{displaymath}
    \widehat \phi (\overline E)\in 
    \cha_{B}^{\ast}(X,\mathcal{D}_{\as})
  \end{displaymath}
  having the following properties:
  \begin{enumerate}
  \item \emph{Functoriality.} When $f:Y\longrightarrow X$ is a
    morphism of regular schemes flat and quasi-projective over $A$, 
    $D'$ a normal crossing divisor on $Y_{R}$ with $f^{-1}(D)\subset
    D'$, and $\overline E$ a log-singular hermitian vector bundle on $X$ then  
    \begin{displaymath}
      f^{\ast}(\widehat {\phi }(\overline E))=
      \widehat {\phi }(f^{\ast}\overline E).
    \end{displaymath}
  \item \emph{Normalization.} When $\overline E=\overline L_{1}\oplus
    \dots \oplus \overline L_{n}$ is a orthogonal direct sum of hermitian
    line bundles, then 
    \begin{displaymath}
      \widehat {\phi }(\overline E)=\phi (\ca_{1}(\overline L_{1}),\dots
      ,\ca_{1}(\overline L_{n})).
    \end{displaymath}
  \item \emph{Twist by a line bundle.} Let
    \begin{displaymath}
      \phi (T_{1}+T,\dots ,T_{n}+T)=
      \sum_{i\ge 0}\phi _{i}(T_{1},\dots ,T_{n})T^{i}. 
    \end{displaymath}
    Let $\overline L$ be a log-singular hermitian vector bundle. Then 
    \begin{displaymath}
      \widehat {\phi }(\overline E\otimes \overline L)=
      \sum_{i}\widehat \phi _{i}(\overline E)\ca_{1}(\overline L). 
    \end{displaymath}
  \item \emph{Compatibility with characteristic forms.}
    \begin{displaymath}
      \omega (\widehat \phi (\overline E))=\phi (E,h).
    \end{displaymath}
  \item \emph{Compatibility with the change of metrics.} If $h'$ is
  another log-singular
  hermitian metric then 
  \begin{displaymath}
    \widehat{\phi }(E,h)=\widehat {\phi }(E,h')+
    \amap(\widetilde \phi _{1}(E,h',h)).
  \end{displaymath}
\item\label{item:1} \emph{Compatibility with the definition of Gillet
  and Soul\'e.} 
  If  $D$ is empty, let $\psi $ be the isomorphism 
  $\cha^{\ast}(X,\mathcal{D}_{\as})\longrightarrow \cha^{\ast}(X)$
  of theorem \ref{thm:8}
  and let $\widehat \phi _{\GS}(\overline E)\in \cha^{\ast}(X)$ be the
  characteristic 
  class defined in  \cite{GilletSoule:vbhm}. Then 
  \begin{displaymath}
    \psi (\widehat \phi(\overline E))=\widehat \phi _{\GS}(\overline E). 
  \end{displaymath}
  \end{enumerate}
\end{theorem}
\begin{proof}
  If $\DD$ is empty, we define $\widehat \phi
  (\overline E)=\psi^{-1} (\widehat \phi _{\GS}(\overline E))$.  If
  $\DD$ is not empty but $h$ is smooth in the whole $X_{\mathbb{R}}$
  then we define $\widehat \phi
  (\overline E)$ by functoriality, using the tautological morphism
  $(X,\DD)\longrightarrow (X,\emptyset)$.  

  If $D$ is not empty and $\overline {E}=(E,h)$ is a log-singular hermitian
  vector bundle, we choose any smooth metric $h'$, invariant under
  $F_{\infty}$. Then we define
  \begin{displaymath}
    \widehat{\phi }(\overline E)=\widehat \phi (E,h')+\amap(\widetilde
    \phi (E,h',h)).
  \end{displaymath}
  This definition is independent of the choice of the metric $h'$ because,
  if $h''$ is another smooth $F_{\infty}$-invariant metric, then 
  \begin{align*}
    \widehat \phi (E,h')+&\amap(\widetilde
    \phi (E,h',h))-\widehat \phi (E,h'')-\amap(\widetilde
    \phi (E,h'',h))\\
    &=\amap(\widetilde \phi (E,h'',h'))+
    \amap(\widetilde \phi (E,h',h))+ \amap(\widetilde \phi(E,h,h''))\\
    &=\amap(d_{\mathcal{D}}(E,h'',h',h))\\&=0.
  \end{align*}
  All the properties stated in the theorem can be checked as in
  \cite{GilletSoule:vbhm}.
\end{proof}
\begin{remark}
  If $X$ is projective, the groups
  $\cha^{\ast}(X,\mathcal{D}_{\lgi})$ and
  $\cha^{\ast}(X,\mathcal{D}_{\as})$ agree. Therefore
  the arithmetic characteristic classes also belong to the former
  group. When $X$ is quasi-projective, in order to define
  characteristic classes in the group
  $\cha^{\ast}(X,\mathcal{D}_{\lgi})$, we have to impose
  conditions on the behavior of the hermitian metrics at
  infinity. For instance one may consider 
  smooth at infinity hermitian metrics (see \cite{BurgosWang:hBC}).
\end{remark}

\begin{remark}
  If in theorem \ref{thm:11} we replace good hermitian vector bundle
  for log hermitian vector bundle and pre-log-log forms for log-log
  forms (implicit in the definition of $\mathcal{D}_{\as}$. Then the
  result remains true. 
\end{remark}

\subsection{Arithmetic $K$-theory of log-singular hermitian vector bundles }
\label{sec:arithmetic-k-theory}

\nnpar{Log-singular arithmetic $K$-theory.} We want to generalize the
definition of arithmetic $K$-theory given by Gillet and Soul\'e in
\cite{GilletSoule:vbhm} to cover log-singular hermitian metrics. 

We
write 
\begin{displaymath}
  \widetilde{\mathcal{D}}_{\as}(X)=\bigoplus _{p}\widetilde
  {\mathcal{D}}^{2p-1}_{\as}(X,p), 
\end{displaymath}
\begin{displaymath}
    {\rm Z}\mathcal{D}_{\as}(X)=\bigoplus _{p}{\rm Z}
  \mathcal{D}^{2p}_{\as}(X,p). 
\end{displaymath}
Let $\ch$ be the power series associated
with the Chern character. In particular it induces Bott-Chern forms
$\widetilde {\ch}$ and arithmetic characteristic classes
$\widehat{\ch}$.

\begin{definition} \label{def:21}
  Let $\underline X$ be as in theorem \ref{thm:11}. Then the group
  $\Ka(X,\mathcal{D}_{\as})$ is the group generated by
  pairs $(\overline E,\eta)$, where $\overline E$ is a log-singular
  hermitian metric on $X$ and $\eta\in
  \widetilde{\mathcal{D}}_{\as}(X)$, with the relations 
  \begin{displaymath}
    (\overline S,\eta')+(\overline Q,\eta'')=
    (\overline E,\eta'+\eta''+\widetilde {\ch}(\overline{\mathcal{E}})),
  \end{displaymath}
  for every
  $\eta', \eta ''\in \widetilde{\mathcal{D}}_{\as}(X)$ and every short
  exact sequence of log-singular hermitian
  vector bundles
\begin{displaymath}
  \overline{\mathcal{E}}:
  0\longrightarrow \overline S \longrightarrow \overline E
  \longrightarrow \overline Q \longrightarrow 0.
\end{displaymath}
\end{definition}

If $\DD$ is empty then this definition agrees with the definition of
Gillet and Soul\'e in \cite{GilletSoule:vbhm}. 

\nnpar{Basic properties.} The following theorem summarizes the basic
properties of the arithmetic $K$-theory groups. They are a consequence
of the corresponding results of \cite{GilletSoule:vbhm} together with
theorem \ref{thm:11}.
 
\begin{theorem} Let $\underline X=(X,\DD)$ be an arithmetic variety
  over $A$ with a fixed normal crossing divisor. Then
  \begin{enumerate}
  \item There are natural maps
    \begin{align*}
      \amap: \widetilde{\mathcal{D}}_{\as}(X) &\longrightarrow 
      \Ka(X,\mathcal{D}_{\as}),\\
      \ch: \Ka(X,\mathcal{D}_{\as})&\longrightarrow 
      Z\mathcal{D}_{\as}(X),\\
      \vmap: \Ka(X,\mathcal{D}_{\as})&\longrightarrow 
      \Ko(X),\\
      \widehat{\ch}:\Ka(X,\mathcal{D}_{\as})& \longrightarrow 
      \bigoplus \cha^{p}_{\mathbb{Q}}(X,\mathcal{D}_{\as}),
    \end{align*}
    given by
    \begin{align*}
      \amap(\eta)&=(0,\eta),\\
      \ch([\overline E,\eta])&=\ch(\overline E)+\dd_{\mathcal{D}}\eta,\\
      \vmap([\overline E,\eta])&=[E],\\
      \widehat{\ch}([\overline E,\eta])&=\widehat {\ch}(\overline
      E)+\amap(\eta), 
    \end{align*}
  \item The product
    \begin{displaymath}
      (\overline E,\eta)\otimes (\overline E', \eta')=
      \left(\overline E\otimes \overline E',(\ch(\overline E)+
      \dd_{\mathcal{D}}\eta)\bullet \eta'+ \eta \bullet \ch(\overline
      E')\right)
    \end{displaymath}
    induce a commutative and associative ring structure on
      $\Ka(X,\mathcal{D}_{\as})$. The maps $\vmap$, $\ch$ and
      $\widehat \ch$ are compatible with 
    this ring structure.
  \item If $\underline Y=(Y,\DD')$ is another arithmetic variety
  over $A$ with a fixed normal crossing divisor and
  $f:X\longrightarrow Y$ is a morphism such that $f^{-1}(\DD')\subset
  \DD$ then there is a pull-back morphism
  \begin{displaymath}
    f^{\ast}:\Ka(Y,\mathcal{D}_{\as})\longrightarrow
    \Ka(X,\mathcal{D}_{\as}), 
  \end{displaymath}
  compatible with the maps $\amap$, $\ch$, $\vmap$ and $\widehat
  {\ch}$.
\item There are exact sequences
  \begin{equation}\label{eq:5}
    \Kl(X)\overset {\rho }{\rightarrow
    }\widetilde{\mathcal{D}}_{\as}(X)
    \overset{\amap}{\rightarrow}
    \Ka(X,\mathcal{D}_{\as})
    \overset{\vmap}{\rightarrow}
    \Ko(X)\rightarrow 0,
  \end{equation}
  and
  \begin{multline}
    \label{eq:4}
    \Kl(X)\overset {\rho }{\rightarrow
    }
    \bigoplus_{p}H_{\mathcal{D}_{\as}}^{2p-1}(X,p)
    \overset{\amap}{\rightarrow}
    \Ka(X,\mathcal{D}_{\as})
    \overset{\vmap+\ch}{\longrightarrow}\\
    \Ko(X)\oplus {\rm Z}\mathcal{D}_{\as}(X)
    \rightarrow \bigoplus_{p}H_{\mathcal{D}_{\as}}^{2p}(X,p)
    \rightarrow 0.
  \end{multline}
  In these exact sequences the map
  $\rho $ is the composition 
  \begin{displaymath}
    \Kl(X)\rightarrow
    \bigoplus_{p}H_{\mathcal{D}}^{2p-1}(X_{\mathbb{R}},\mathbb{R}(p))
    \rightarrow
    \bigoplus_{p}H_{\mathcal{D}_{\as}}^{2p-1}(X,p)
    \subset \widetilde{\mathcal{D}}_{\as}(X),
  \end{displaymath}
  where the first map is Beilinson's regulator.
\item The Chern character
  \begin{displaymath}
    \widehat{\ch}:\Ka(X,\mathcal{D}_{\as})\otimes
      \mathbb{Q} \longrightarrow  
      \bigoplus \cha^{p}_{\mathbb{Q}}(X,\mathcal{D}_{\as})
  \end{displaymath}
  is a ring isomorphism. \hfill $\square$
  \end{enumerate}
  
\end{theorem}

\subsection{Variant for non regular arithmetic varieties}
\label{sec:variants-non-regular}

Since there is no general theorem of resolution of singularities it is
useful to extend the theory of arithmetic Chow groups to the case of
non regular arithmetic varieties. 

\nnpar{Arithmetic Chow groups for non regular arithmetic varieties.}
Let $(A,\Sigma
,F_{\infty})$ be an arithmetic ring with fraction field $F$.  We will
assume that $A$ is
equidimensional and Jacobson. In contrast with the rest of the paper,
in this section, an arithmetic variety 
over $A$ will be a
scheme $X$ that is
quasi-projective  and flat over $\Spec(A)$ and such that the generic
fiber $X_{F}$ is smooth, but that does not need to be regular. Since
$X_{F}$ is smooth, the analytic variety 
$X_{\Sigma }$ is a disjoint union of connected components $X_{i}$ that are
equidimensional of dimension $d_{i}$. For every cohomological complex of
sheaves $\mathcal{F}^{\ast}(\ast)$ on $X_{\Sigma }$ we write
\begin{displaymath}
  \mathcal{F}_{n}(p)(U)=\bigoplus_{i} 
  \mathcal{F}^{2d_{i}-n}(d_{i}-p)(U\cap X_{i}).
\end{displaymath}
Then the definition of Green object and of arithmetic Chow groups of
\cite{BurgosKramerKuehn:cacg} can be easily adapted to the grading by
dimension.

 In this way we can define, for $X$ regular,  homological
Chow groups with 
respect to any $\mathcal{D}_{\log}$-complex $\cc$. These homological
Chow groups will be denoted by $\cha_{\ast}(X,\cc)$. In particular we
are interested in the groups $\cha
_{\ast}(X,\mathcal{D}_{\as})$. But now we can proceed
as in \cite{GilletSoule:aRRt} and we can extend the definition to
the case of non regular arithmetic varieties.  

\nnpar{Basic properties of homological Chow groups.} Following
\cite{GilletSoule:aRRt} we can extend some of the properties of the
arithmetic Chow groups to the non regular case. The proof of the
next results are as in \cite{GilletSoule:aRRt} 2.2.7, 2.3.1 and  
2.4.2 for the algebraic cycles, but using the techniques of
\cite{BurgosKramerKuehn:cacg} for the Green objects.

 \begin{theorem} 
  Let $f:X\longrightarrow Y$ be a morphism of irreducible arithmetic
  varieties 
  over $A$ that is flat or l.c.i. Let $\DD_{Y}$ be a normal crossing 
  divisor on $Y_{\mathbb{R}}$ and $\DD_{X}$ a normal crossing divisor on
  $X_{\mathbb{R}}$ such that $f^{-1}(\DD_{Y})\subset \DD_{X}$. Write
  $\underline 
  X=(X_{\mathbb{R}},\DD_{X})$ and $\underline
  Y=(Y_{\mathbb{R}},\DD_{Y})$. Then 
  there is defined an inverse image morphism  
  \begin{displaymath}
    f^{\ast}:\cha_{p}(Y,\mathcal{D}_{\as})
    \longrightarrow 
    \cha_{p+d}(X,\mathcal{D}_{\as}),
  \end{displaymath}
  where $d$ is the relative dimension.
\hfill $\square$
\end{theorem}

\begin{theorem}
  Let $f:X\longrightarrow Y$ be a map of arithmetic varieties, with $Y$
  regular.  Let $\DD_{Y}$ be a normal crossing 
  divisor on $Y_{\mathbb{R}}$ and $\DD_{X}$ a normal crossing divisor on
  $X_{\mathbb{R}}$ such that $f^{-1}(\DD_{Y})\subset \DD_{X}$. Write
  $\underline 
  X=(X_{\mathbb{R}},\DD_{X})$ and $\underline
  Y=(Y_{\mathbb{R}},\DD_{Y})$. Then there is a cap product
  \begin{displaymath}
    \cha^{p}(Y,\mathcal{D}_{\as})\otimes
    \cha_{q}(X,\mathcal{D}_{\as})
    \longrightarrow 
    \cha_{q-p}(X,\mathcal{D}_{\as})    
  \end{displaymath}
  that turns $\cha_{\ast}(X,\mathcal{D}_{\as})$ into a
  graded $\cha^{\ast}(Y,\mathcal{D}_{\as})$
  module. Moreover 
  this cap product is compatible with inverse images (when defined).
\hfill $\square$
\end{theorem}

\nnpar{Arithmetic $K$-theory.} The definition of arithmetic $K$-theory
carries over to the case of non regular arithmetic varieties without
modification (see \cite{GilletSoule:aRRt} 2.4.2). Thus we obtain a
contravariant functor $(X,\DD)\longmapsto
\Ka(X,\mathcal{D}_{\as})$ from arithmetic varieties with
a fixed normal crossing divisor to rings.  

\begin{theorem} Let $X$ be an arithmetic variety. Let $\DD_{X}$ a
  normal crossing divisor on 
  $X_{\mathbb{R}}$. Write
  $\underline 
  X=(X_{\mathbb{R}},\DD_{X})$. Then 
 there is a biadditive pairing
\begin{align*}
\Ka(X,\mathcal{D}_{\as}) \otimes
\cha_*(X,\mathcal{D}_{\as}) \longrightarrow
\cha_*(X,\mathcal{D}_{\as})_\QQ, 
\end{align*}
which we write $\alpha \otimes x \mapsto
\widehat{\operatorname{ch}}(\alpha) \cap x$, with the following
properties
\begin{enumerate}
\item Let $f: X \to Y$ be a morphism of arithmetic varieties, with $Y$
  regular. Let $\DD_{Y}$ be a normal crossing 
  divisor on $Y_{\mathbb{R}}$ such that $f^{-1}(\DD_{Y})\subset
  \DD_{X}$. Write 
  $\underline
  Y=(Y_{\mathbb{R}},\DD_{Y})$. If $\alpha \in
  \Ka(Y,\mathcal{D}_{\as})$ and $x \in
  \cha_*(X,\mathcal{D}_{\as})$, then  
\begin{align*}
\widehat{\operatorname{ch}}(f^*\alpha) \cap x =
\widehat{\operatorname{ch}}(\alpha)._f x.
\end{align*}
\item If $(0,\eta) \in \Ka(X,\mathcal{D}_{\as})$ and $x
  \in \cha^*(X,\mathcal{D}_{\as})$, then 
\begin{align*}
\widehat{\operatorname{ch}}((0,\eta)) \cap x = \amap( \eta \omega(x)).
\end{align*}
\item If $\alpha \in \Ka(X,\mathcal{D}_{\as})$ and $x \in
  \cha_*(X,\mathcal{D}_{\as})$, then 
\begin{align*}
\omega(\widehat{\operatorname{ch}}(\alpha) \cap x) =
\operatorname{ch}(\alpha) \land \omega(x).
\end{align*}
\item The pairing makes $\cha_*(X,\mathcal{D}_{\as})_\QQ$
  into a $\Ka(X,\mathcal{D}_{\as})$ module; i.e., 
  for all $\alpha, \beta \in \Ka(X)$, and $x \in \cha^*(X)$ we have
\begin{align*}
\widehat{\operatorname{ch}}(\alpha) \cap 
(\widehat{\operatorname{ch}}(\beta) \cap  x) =
\widehat{\operatorname{ch}}(\alpha\beta) \cap  x 
\end{align*}
\item If $f: X \to Y$ is a flat or l.c.i. morphism of arithmetic
  varieties, let $\alpha\in \Ka(X,\mathcal{D}_{\as})$ and
  $x \in \cha^*(X,\mathcal{D}_{\as})$. Then 
\begin{align*}
\widehat{\operatorname{ch}}(f^* \alpha)  \cap f^* x =  
f^*( \widehat{\operatorname{ch}}(\alpha)  \cap  x)
\end{align*}   
\end{enumerate}
\end{theorem}

\begin{proof}
  Follow \cite{GilletSoule:aRRt} 2.4.2. but using theorem \ref{thm:18}
  to prove the independence of the choices.
\end{proof}

\subsection{Some remarks on the properties of
  $\cha ^{\ast}(X,\mathcal{D}_{\as})$} 
\label{sec:summary-theory}

In \cite{MaillotRoessler:cdl}, V. Maillot and D. Roessler
have announced a preliminary version of the theory developed in this
paper. The final theory has some minor differences that do not affect
the heart of \cite{MaillotRoessler:cdl}. The aim of this section is to
compare the both theories.

We fix an arithmetic ring $(A,\Sigma ,F_{\infty})$ and we consider
pairs $\underline X=(X,D)$, where $X$ is an arithmetic 
variety over $A$ and $D$ is a normal crossing divisor of $X_{\Sigma }$
defined invariant under $F_{\infty}$. 

A log-singular hermitian vector bundle $\overline E$ is a pair
$(E,h)$, where $E$ is a vector bundle over $X$ and $h$ is a hermitian
metric on $E_{\Sigma }$, invariant under $F_{\infty}$ and log-singular
along $D$. Observe that the notion of log-singular hermitian metric is
not the same as the notion of good hermitian metric. This is not
important by two 
reasons. First, as we will see in the next section, the main examples of
good hermitian vector bundles, the fully decomposed automorphic vector
bundles, are good and log-singular. Second, if one insist on using good
hermitian vector bundles, one can replace pre-log and pre-log-log
forms for log and log-log forms to obtain an analogous theory. This
alternative theory has worse cohomological properties (we have not
proved the Poincar\'e lemma for pre-log and pre-log-log forms), but
the arithmetic intersection numbers computed by both theories agree.    

To each pair $\underline X=(X,D)$, we have assigned an
$\mathbb{N}$-graded abelian group
$\cha^{\ast}(X,\mathcal{D}_{\as})$, that satisfies,
among others, the following properties:
\begin{enumerate}
\item The group $\cha^{\ast}(X,\mathcal{D}_{\as})$ is
  provided with an associative, commutative and unitary ring structure,
  compatible with the grading.
\item If $X$ is proper over $\Spec A$, there is a direct image
  group homomorphism
  $f_{\ast}:\cha^{d+1}(X,\mathcal{D}_{\as})\longrightarrow
  \cha^{1}(\Spec A)$, where $d$ is the relative dimension.
\item For every integer $r\ge 0$ and every log-singular hermitian
  vector bundle there is defined the arithmetic $r$-th Chern class
  $\ca_{r}(\overline E)\in
  \cha^{r}(X,\mathcal{D}_{\as})$.
\item   Let $g:X\longrightarrow Y$ be a morphism of arithmetic varieties
  over $A$, and  let $D$, $E$ be normal crossing 
  divisor on $X_{\mathbb{R}}$ and  $Y_{\mathbb{R}}$ respectively, such
  that $g^{-1}(E)\subset D$. Write $\underline 
  X=(X_{\mathbb{R}},D)$ and $\underline Y=(Y_{\mathbb{R}},E)$. Then
  there is defined an inverse image morphism  
  \begin{displaymath}
    g^{\ast}:\cha^{\ast}(Y,\mathcal{D}_{\as})
    \longrightarrow 
    \cha^{\ast}(X,\mathcal{D}_{\as}).
  \end{displaymath}
  Moreover, it is a morphism of rings after tensoring with
  $\mathbb{Q}$.
\item For every $r\ge 0$, it holds the equality
  $g^{\ast}(\ca_{r}(\overline E))=\ca_{r}(g^{\ast}(\overline E))$.
\item There is a forgetful morphism $\zeta:
  \cha^{\ast}(X,\mathcal{D}_{\as})\longrightarrow
  \CH^{\ast}(X)$, compatible with inverse images and Chern classes. 
\item There is a complex of groups
  \begin{displaymath}
    H^{2p-1}_{\mathcal{D}}(X_{\mathbb{R}},\mathbb{R}(p)) 
    \stackrel{\amap}{\longrightarrow} 
    \cha^p(X,\mathcal{D}_{\as})\stackrel{(\zeta,\omega)}{\longrightarrow}  
    \CH^p(X) \oplus {\rm Z}\mathcal{D}_{\as}^{2p}(X,p),
  \end{displaymath}
  that is an exact sequence when $X_{\Sigma }$ is projective. 
  Observe that the
  group ${\rm Z}\mathcal{D}_{\as}^{2p}(X,p)$ does not agree with the
  group denoted ${\rm Z}^{p,p}(X(\mathbb{C}),D)$ in
  \cite{MaillotRoessler:cdl} \S 1 (7). The former is made of forms
  that are log-log along $D$ and the later by forms that are good
  along $D$. Again this is not important by two reasons. First, the image by
  $\omega $ of the arithmetic Chern classes of fully decomposed
  automorphic vector bundles lies in the intersection of the good and
  log-log forms. Second, the complex of log-log forms shares all the
  important properties of the complex of good forms (see proposition
  \ref{prop:5}). 
\item The morphism $(\zeta, \omega)$ is a ring homomorphism; the
  image of $\amap$ is a square zero ideal. Moreover It holds the equality 
  \begin{displaymath}
    \amap(x)\cdot y=\amap(x\cdot \cl(\zeta(y))),
  \end{displaymath}
  where $x\in H^{2p-1}_{\mathcal{D}}(X_{\mathbb{R}},\mathbb{R}(p))$,
  $y\in \cha^p(X,\mathcal{D}_{\as})$, $\cl$ is the class map, the
  product in the left hand side is the product in the arithmetic Chow
  groups and the product in the right hand side is the product in
  Deligne-Beilinson cohomology.
\item When $D$ is empty, there is a
  canonical isomorphism
  $\cha^{\ast}(X,\mathcal{D}_{\as})\longrightarrow \cha^{\ast}(X)$,
  compatible with the previously discussed structures. Note that we have
  dropped the projectivity assumption in \cite{MaillotRoessler:cdl} \S
  1 (9). Observe moreover that, if we use the alternative theory with
  pre-log-log forms, then this property is not established.    
\end{enumerate}

\section{Automorphic vector bundles}
\label{sec:shim-vari}

\subsection{Automorphic bundles and log-singular hermitian metrics}

\nnpar{Fully decomposed automorphic bundles.}
Let $B$ be a bounded symmetric domain. Then, by definition 
$B=G/K$, where $G$ is a semisimple adjoint group and $K$ is a maximal
compact subgroup. Inside the complexification $G_{\mathbb{C}}$ of $G$,
there is a parabolic subgroup of the form $P_{+}\cdot K_{\mathbb{C}}$ such
that $K=G\cap P_{+}\cdot K_{\mathbb{C}}$ and $G\cdot(P_{+}\cdot
K_{\mathbb{C}})$ is open 
in $G_{\mathbb{C}}$. This induces an open $G$-equivariant immersion 
\begin{displaymath}
  \xymatrix{
    \ B\ \ar@{^{(}->}[r]^{\iota}\ar@{=}[d] & \check{B} \ar@{=}[d]\\
    \ G/K\   \ar@{^{(}->}[r] & G_{\mathbb{C}}/P_{+}\cdot K_{\mathbb{C}}.
   }
\end{displaymath}
$\check{B}=G_{\mathbb{C}}/P_{+}\cdot K_{\mathbb{C}}$ is a projective
rational variety, and this immersion is compatible with the  complex
structure of $B$.

Let $\sigma :K\longrightarrow GL(n,\mathbb{C})$ be a representation of
$K$. Then $\sigma $ defines a $G$-equivariant vector bundle $E_{0}$ on
$B$. We complexify $\sigma $ and we extend it trivially to $P_{+}\cdot
K_{\mathbb{C}}$ by letting it kill $P_{+}$. Then $\sigma $ defines a
holomorphic $G_{\mathbb{C}}$-equivariant vector bundle $\check{E}_{0}$
on $\check{B}$ with $E_{0}=\iota^{\ast}(\check{E}_{0})$. This induces
a holomorphic structure on $E_{0}$. Observe that different extensions
of $\sigma $ to $P_{+}\cdot K_{\mathbb{C}}$ will define different
holomorphic structures on $E_{0}$.
   
  Let $\Gamma $ be a 
neat arithmetic 
subgroup of $G$ acting on $B$. Then $X=\Gamma \backslash B$ is a smooth
quasi-projective 
complex variety, and $E_{0}$ defines a holomorphic vector bundle $E$
on $X$. Following \cite{HarrisZucker:BcSvIII}, the vector bundles
obtained in this way (with $\sigma 
$ extended trivially) will be called \emph{fully decomposed
  automorphic vector bundles}. Since in this paper we will not treat
more general automorphic vector bundles we will just call them
automorphic vector bundles.

Let $h_{0}$ be a $G$-equivariant hermitian metric on $E_{0}$. Such
metrics exist by the compactness of $K$. Then $h_{0}$ determines a
hermitian metric $h$ on $E$. 

\begin{definition}
A hermitian vector bundle $(E,h)$ as
before will be called an \emph{automorphic hermitian vector bundle}.   
\end{definition}

Let $\overline X$ be a smooth toroidal compactification of
$X$ with $D=\overline X\setminus X$ a normal crossing divisor.
We recall  
the following result of Mumford  \cite{Mumford:Hptncc}, Thm. 3.1:
\begin{theorem}  
The automorphic vector bundle $E$  admits a unique extension to a
 vector bundle $E_{1}$ over
 $\overline X$ such that 
 $h$ is a singular hermitian 
metric which is good along $D$. \hfill $\Box$
\end{theorem}

By abuse of language the extension $(E_{1},h)$ will be also called an
automorphic hermitian vector bundle.

Our task now is to improve slightly Mumford's theorem:
\begin{theorem} \label{thm:17}
  The automorphic hermitian vector bundle $(E_{1},h)$ is a $\infty$-good hermitian
  vector bundle, therefore it is
  log-singular along $\DD$.
\end{theorem}
\begin{proof}  
The proof of this result will take us the rest of the section. The
technique follows closely the proof of theorem \cite{Mumford:Hptncc},
3.1. Instead of repeating the whole proof of Mumford we will only
point out the results needed to bound all the derivatives of the
functions involved. 

\nnpar{Cones and Jordan algebras.}
Let $V$ be a real vector space and let $C\subset V$ be a homogeneous
self-adjoint cone. See \cite{AMRT:sclsv} for the theory of homogeneous
self-adjoint cones and their relationship with Jordan algebras.
We will only recall here some basic  facts.
 
Let $G\subset GL(V)$ be the group of linear maps that preserve
$C$. Since $C$ is homogeneous, $G$ acts transitively on $C$. We will
denote by $\mathfrak{g}$ the Lie algebra of $G$. For any point $x\in
C$, let $K_{x}=\Stab(x)$. It is a maximal compact subgroup of $G$. Let
$\mathfrak{k}_{x}$ be the Lie algebra of $K_{x}$ and let  
\begin{displaymath}
  \mathfrak{g}=\mathfrak{k}_{x}\oplus \mathfrak{p}_{x} 
\end{displaymath}
be the associated Cartan decomposition. Let $\sigma_{x} $ be the
Cartan 
involution. Let us choose a point $e\in C$. Let 
$\left< \ ,\ \right>=\left< \ ,\ \right>_{e}$ be a positive-definite
scalar product such that $\sigma _{e}(g)={}^{t}g^{-1}$ for all $g\in
G$. Then $C$ is self-adjoint with respect to this inner product. For
any point $x\in C$, let us choose $g\in G$ such that $x=ge$. We will
identify $V$ with $T_{C,x}$. For $t_{1},t_{2}\in V$ we will write    
\begin{displaymath}
  \left<t_{1},t_{2}\right>_{x}=
  \left<g^{-1}t_{1},g^{-1}t_{2}\right>_{e}.  
\end{displaymath}
The right hand side in independent of $g$ because $\left< \ ,\
\right>_{e}$ is $K_{e}$-invariant. These products define a
$G$-invariant Riemannian metric in $C$ that is denoted $ds^{2}_{C}$.  

The elements of $\mathfrak{g}$
act on $V$ by endomorphisms. This action can be seen as the
differential of the $G$ 
action at $e\in V$, or given by the inclusion $\mathfrak{g}\subset
\mathfrak{gl}(V)$. For any $x\in C$ there are isomorphisms  
\begin{displaymath}
  \mathfrak{p}_{x}\overset{\cong}{\longrightarrow
    }\mathfrak{p}_{x}.x=V \quad\text{and}\quad
  P_{x}=\exp(\mathfrak{p}_{x})\overset{\cong}{\longrightarrow
    }P_{x}.x=C.
\end{displaymath}
The elements of $\mathfrak{p}_{x}$ act on $V$ by self-adjoint
endomorphisms with respect to $\left< \ ,\ \right>_{x}$.

Every $\mathfrak{p}_{x}$ has a structure of Jordan algebra defined by 
\begin{displaymath}
  (\pi.\pi').x=\pi .(\pi '.x).
\end{displaymath}
The isomorphism $\mathfrak{p}_{x}\longrightarrow V$ define a Jordan
algebra structure on $V$ that we denote by
$t_{1}\underset{x}{.}t_{2}$. Observe that $x$ is the unit element for
this Jordan algebra structure.

We summarize the compatibility relations between the objects defined so
far and the action of the group. Let $x=g.e$
\begin{align*}
  K_{x}&=\Ad(g) K_{e}=g K_{e} g^{-1},\\
  \mathfrak{k}_{x}&=\ad(g) \mathfrak{k}_{e}=g \mathfrak{k}_{e}
  g^{-1},\\
  \mathfrak{p}_{x}&=\ad(g) \mathfrak{p}_{e}=g \mathfrak{p}_{e} g^{-1}.
\end{align*}
There is a commutative diagram
\begin{displaymath}
  \begin{CD}
    \mathfrak{p}_{e}@>\ad(g)>> \mathfrak{p}_{x}\\
    @V .e VV @V .x VV\\
    V @> g. >> V
  \end{CD}
\end{displaymath}
The horizontal arrows in the above diagram are morphisms of Jordan
algebras. In particular
\begin{displaymath}
  g.(t_{1}\underset{e}{.}t_{2})=gt_{1}\underset{x}{.}gt_{2}.
\end{displaymath}
When a unit element $e$ is chosen we will write
$t_{1}.t_{2}$ and  $\left< \ ,\ \right>$ instead of
$t_{1}\underset{e}{.}t_{2}$ and $\left< \ ,\ \right>_{e}$.

\nnpar{ Derivatives with respect to the base point.}
We now study the derivatives of the scalar product and the Jordan
algebra product when we move the base point.
\begin{lemma} \label{lemm:derandprod}
  Let $t_{1},t_{2},t_{3}\in V$. Then
  \begin{enumerate}
  \item \label{item:2}$D_{t_{1}}(\left<t_{2},x^{-1}\right>_{e})
    =-\left<t_{2},t_{1}\right>_{x}.$
  \item \label{item:3}$D_{t_{3}}\left<t_{1},t_{2}\right>_{x}=-(
    \left<t_{3}\underset{x}{.}t_{1},t_{2}\right>_{x}+
    \left<t_{1},t_{3}\underset{x}{.}t_{2}\right>_{x})=
    -2 \left<t_{1},t_{3}\underset{x}{.}t_{2}\right>_{x}.$
  \item \label{item:4}$    D_{t_{3}} (t_{1}\underset{x}{.}t_{2})=-(
    (t_{3}\underset{x}{.}t_{1})\underset{x}{.}t_{2}+
    t_{1}\underset{x}{.}(t_{3}\underset{x}{.}t_{2})).$
  \end{enumerate}
\end{lemma}
\begin{proof}
  The proof of \ref{item:2} is in \cite{Mumford:Hptncc} par. 244.
  To prove  \ref{item:3}
  write $t_{3}=M.x$ with $M\in \mathfrak{p}_{x}$. Then $\alpha (\delta
  )=\exp(\delta M).x$ is a curve  with $\alpha (0)=x$ and $\alpha
  '(0)=t_{3}$. Therefore 
  \begin{align*}
    D_{t_{3}}\left<t_{1},t_{2}\right>_{x}&=\frac{\dd}{\dd\delta }
    \left<t_{1},t_{2}\right>_{\exp(\delta M).x}|_{\delta =0}\\
    &=\frac{\dd}{\dd\delta }
    \left<\exp(\delta M)^{-1}.t_{1},
      \exp(\delta M)^{-1}.t_{2}\right>_{x}|_{\delta =0}\\
    &=-(\left<M.t_{1},t_{2}\right>_{x}+
    \left<t_{1},M.t_{2}\right>_{x})\\ 
    &=-(\left<t_{3}.t_{1},t_{2}\right>_{x}+
    \left<t_{1},t_{3}.t_{2}\right>_{x}).      
  \end{align*}
  The second equality of  \ref{item:3} follows from the fact that $M$ acts by
  an endomorphism that is self-adjoint with respect to $\left<\ ,\
  \right>_{x}$. 

  The proof of \ref{item:4} is completely analogous.
\end{proof}

We will denote by $\|\ \|_{x}$ the norm associated to the inner
product $\left<\ ,\  \right>_{x}$.
\begin{lemma}\label{lemm:boundprod}
  There is a constant $K>0$ that, for all $x\in C$ and $t_{1},t_{2}\in
  V$, 
  \begin{displaymath}
    \|t_{1}\underset{x}{.}t_{2}\|_{x}\le K
    \|t_{1}\|_{x}\|t_{2}\|_{x}.
  \end{displaymath}
\end{lemma}
\begin{proof}
  In $\mathfrak{p}_{e}$ we may define the norm
  \begin{displaymath}
    \|M\|'_{e}=\sup_{t\in V}\frac{\|M.t\|_{e}}{\|t\|_{e}}. 
  \end{displaymath}
  Via the isomorphism $\mathfrak{p}_{e}\longrightarrow V$ it induces a
  norm in $V$ given by
  \begin{displaymath}
    \|t_{1}\|'_{e}=\sup_{t\in
    V}\frac{\|t_{1}\underset{e}{.}t\|_{e}}{\|t\|_{e}}. 
  \end{displaymath}
  Since any two norms in a finite dimensional vector space are
  comparable, there is a constant $K>0$ such that, for all $t$,  
  \begin{displaymath}
    \|t\|_{e}'\le K\|t\|_{e}.
  \end{displaymath}
  Therefore 
  \begin{displaymath}
    \|t_{1}\underset{e}{.}t_{2}\|_{e}\le \|t_{1}\|_{e}'\|t_{2}\|_{e}
    \le K \|t_{1}\|_{e}\|t_{2}\|_{e}.
  \end{displaymath}
  But for any $x=ge$ we have
  $$
  \|t_{1}\underset{x}{.}t_{2}\|_{x}=
  \|g^{-1}t_{1}\underset{e}{.}g^{-1}t_{2}\|_{e}
  \le K \|g^{-1}t_{1}\|_{e}\|g^{-1}t_{2}\|_{e}=
  K \|t_{1}\|_{x}\|t_{2}\|_{x}.
  $$
\end{proof}

\nnpar{Maximal $\mathbb{R}$-split torus.} We fix a unit element $e\in
C$. This fixes also the Jordan algebra structure of $V$, and we write
$K=K_{e}$ and $\mathfrak{p}=\mathfrak{p}_{e}$. Let $A\subset
\exp(\mathfrak{p})$ be a maximal $\mathbb{R}-split$ 
torus, with $A=\exp(\mathfrak{a})$. Then
$\exp(\mathfrak{p})=K.A.K^{-1}$ and $C=K.A.e$. A useful result that
is proved in \cite{AMRT:sclsv} II \S3 is the following
\begin{proposition}
There exist a maximal set of
mutually orthogonal 
idempotents $\epsilon _{1},\dots ,\epsilon _{r}$ of $V$, with
$e=\epsilon _{1}+\dots +\epsilon _{r} $, such that  
\begin{displaymath}
  \mathfrak{a}.e = \sum_{i=1}^{r} \mathbb{R}\epsilon _{i} 
\text{ and } 
  A.e = \sum_{i=1}^{r} \mathbb{R}^{+}\epsilon _{i}. 
\end{displaymath}
Moreover $C\cap \mathfrak{a}.e =  A.e$.\hfill $\square$   
\end{proposition}

  We can introduce on $A$ the coordinates given by
  \begin{displaymath}
    A\cong A.e=\sum_{i=1}^{r}\mathbb{R}^{+}\epsilon _{i}\cong
    (\mathbb{R}^{+})^{r}.
  \end{displaymath}

As an application of the previous result we prove a bound for the norm
of $x^{-1}$. 
\begin{lemma} \label{lemm:8} Let $\sigma \in C$.
  There exists a constant $K$ such that $\|x^{-1}\|\le K$ for all
  $x\in \sigma +C$. 
\end{lemma}
\begin{proof}
  Since $\bigcup_{\lambda >0} (\lambda e+C)=C$ we may assume that
  $\sigma =\lambda e$ for some $\lambda >0$. Since $K$ is compact and 
  \begin{displaymath}
    \lambda e+C=K(\lambda e+ A.e)
  \end{displaymath}
  it is enough to bound $x^{-1}$ for $x\in \lambda e+ A.e$. 
  If $x\in  \lambda e+ A.e$, then we can write, using the above
  coordinates of $A$, $x=a.e$ with
  $a=(a_{1},\dots ,a_{r})$  and all $a_{i}\ge \lambda $. Then
  $x^{-1}=a^{-1}.e$. Since on a finite dimensional vector space any
  two norms are comparable we obtain 
  \begin{displaymath}
    \|x^{-1}\|^{2}\le K_{1}(a_{1}^{-2}+\dots +a_{k}^{-2})\le
    K_{2}/\lambda ^{2}. 
  \end{displaymath}
\end{proof}

\nnpar{equivariant symmetric representations.} Let $C_{n}$ be the cone
of positive definite $n\times n$ hermitian matrices. 
An equivariant symmetric representation of
dimension $n$ is a pair $(\rho ,H)$, where  $\rho :G\longrightarrow
GL(n,\mathbb{C})$ is a representation and $H:C\longrightarrow C_{n}$
is a map such that
\begin{enumerate}
\item (equivariant) $H(gx)=\rho (g)H(x){}^{t}\overline{\rho (g)},
  \text{ for all }x\in C,\ g\in G.$
\item (symmetric) $\rho (g^{\ast})=
  H(e).{}^{t}\overline{\rho (g)}^{-1}.H(e)^{-1}, \text{ for all }g\in G$.
\end{enumerate}

We will consider an equivariant symmetric representation $(\rho
,H_{t})$ 
with $H_{t}$ depending differentially on a parameter $t\in T$ with $T$
compact, as in 
\cite{Mumford:Hptncc} pp. 245, 246.

\nnpar{Bounds of $H$ and $\det H^{-1}$.}
The first step is to bound the entries of $H_{t}$ and $\det
H(t)^{-1}$.  This is done in \cite{Mumford:Hptncc} proposition 2.3.

\begin{proposition} \label{prop:11}
  For all $\sigma \in C$ there is a constant $K>0$ and an integer $N$
  such that  
  \begin{displaymath}
    \|H_{t}(x)\|,\ |\det H_{t}(x)|^{-1} \le K\left<x,x\right>^{N},
    \text{ for all }x\in \sigma +C. 
  \end{displaymath}
  \hfill $\square$
\end{proposition}

The following results 
of Mumford (\cite{Mumford:Hptncc} proposition 2.4 and 2.5) are the
starting point to bound the entries of $D_{v} H_{t}.H_{t}^{-1}$ but they
will also be used to bound the derivatives of $H_{t}$.  

\begin{proposition} \label{prop:12} Let $\xi \in V$.
  For all $1\le \alpha ,\beta \le n$ let $(D_{\xi }H_{t}.
  H_{t}^{-1})_{\alpha ,\beta }$ be the $(\alpha ,\beta )$-th entry
  of this matrix. There is a linear map 
  \begin{displaymath}
    C_{\alpha \beta ,t}:V\longrightarrow V
  \end{displaymath}
  depending differentially on $t$, such that 
  \begin{displaymath}
    (D_{\xi }H_{t}.H^{-1})_{\alpha \beta }(x)=
    \left< C_{\alpha \beta ,t}(\xi ),x^{-1}\right>.
  \end{displaymath}
  Moreover $C_{\alpha \beta ,t}$ has the property
  \begin{displaymath}
    \left. 
      \begin{matrix}
        \xi ,\eta \in \overline C\\
        \left< \xi ,\eta \right>=0
      \end{matrix}
      \right\} \Rightarrow \left< C_{\alpha \beta ,t}(\xi) ,\eta
      \right>=0 
  \end{displaymath}
\hfill $\square$
\end{proposition}

\begin{proposition} \label{prop:invderdelta}
  For all vector fields $\delta $ on $T$, $\delta H_{t}.
  H_{t}^{-1}(x)$ is independent of $x$. \hfill $\square$
\end{proposition}

\begin{proposition}\label{prop:16}
  Let $\sigma \in C$, let $P$ be a differential operator on $T$ and
  let $\xi _{1},\dots ,\xi _{d}\in V$. Then there is a
  constant $K>0$ and an integer $N$ such that
  \begin{displaymath}
    \|D_{\xi _{1}}\dots D_{\xi _{d}}P H_{t}(x)\|, \  
    | D_{\xi _{1}}\dots D_{\xi _{d}} P \det H_{t}(x) | \le
    K\left<x,x\right>^{N}, \text{ all }x\in \sigma +C.  
  \end{displaymath}
\end{proposition}
\begin{proof} 
  in view of proposition \ref{prop:invderdelta} and since $T$ is
  compact, it is enough to
  consider the case $P=\Id$. Now, by
  proposition \ref{prop:12} and the fact that  
  \begin{displaymath}
    D_{\xi_{i}} x^{-1}=-x^{-1}.(x^{-1}.\xi _{i})
  \end{displaymath}
  we can prove by induction that
  \begin{displaymath}
    D_{\xi _{1}}\dots D_{\xi _{d}} H_{t}(x)=
    M(C(\xi _{1},\dots ,\xi _{d},x)).H_{t}(x),
  \end{displaymath}
  where $M:V\longrightarrow M_{n}(\mathbb{C})$ is linear and
  $C:V\longrightarrow V$ is linear  on $\xi _{1},\dots ,\xi _{d}$ and
  is a polynomial in $x^{-1}$.

  Then the proposition follows from proposition \ref{prop:11} and
  lemma \ref{lemm:8}.
\end{proof}

\nnpar{Bounds of $\delta H.H^{-1}$.}
Let $e=\epsilon _{1}+\dots +\epsilon _{r}$ be a
maximal set of orthogonal idempotents and let $A$ be the corresponding
$\mathbb{R}$-split maximal torus. Let $C_{i}$ be the boundary
component containing $\epsilon _{i+1}+\dots +\epsilon _{r}$ (see
\cite{AMRT:sclsv} II \S3). Let $\widetilde C=C\cup C_{1}\cup\dots \cup
C_{r}\cup 0$ and let $P$ be the parabolic subgroup that stabilizes the
flag $\{C_{i}\}$.

In order to use proposition \ref{prop:12} to bound
$D_{v}H_{t}.H_{t}^{-1}$ and its derivatives we will need the following
result (\cite{Mumford:Hptncc} proposition 2.6).

\begin{proposition}\label{prop:13} Let $\xi _{1}, \xi _{2}\in \widetilde C$ and
  let $\xi '_{1}\in V$ satisfy
  \begin{displaymath}
    \left. 
      \begin{matrix}
        \eta \in \overline C\\
        \left< \xi _{1},\eta \right>=0
      \end{matrix}
      \right\} \Rightarrow \left< \xi'_{1} ,\eta
      \right>=0 .
  \end{displaymath}
  Then, for every compact subset $\omega \subset P$ there is a
  constant $K$ such that
  \begin{enumerate}
  \item $|\left< \xi '_{1},x^{-1}\right>|\le 
    K\| \xi _{1} \|_{x}$, for all $x\in
    \omega .A.e$. 
  \item $|\left< \xi '_{1}, \xi _{2}\right>|\le 
    K\| \xi _{1} \|_{x}
      \|\xi _{2} \|_{x}$, for all $x \in
    \omega .A.e$.
  \end{enumerate}
\end{proposition}

Now we can bound the derivatives of $D_{v}H_{t}.H^{-1}_{t}$ in terms
of the Riemannian metric $ds^{2}_{C}$. Let $N=\dim V$ and let $\xi
_{1},\dots ,\xi _{N}\in 
\widetilde C$ spawn $V$. 

\begin{proposition} \label{prop:14} Let $\delta $ be a vector field in
  $T$, let $P$ be a differential operator 
  that is a product of vector fields in $T$, let $(j _{i})_{i=1}^{n}$
  be a finite sequence of elements of $\{1,\dots 
  ,N\}$
 Let $\omega $ be a compact
  subset of $P$. Then there is a constant such that, for all $x\in
  \omega .A.e$, 
  \begin{align*}
    \|D_{\xi _{j_{1}}}\dots D_{\xi _{j_{n}}} 
    P(D_{\delta }H_{t}.H_{t}^{-1})\|&\le
    K \|\xi _{j_{1}}\|_{x}\dots \|\xi_{j_{n}}\|_{x}\\
    \|D_{\xi _{j_{1}}}\dots D_{\xi _{j_{n-1}}}P (D_{\xi
    _{j_{n}}}H_{t}.H_{t}^{-1})\|&\le 
    K \|\xi _{j_{1}}\|_{x}\dots \|\xi_{j_{n}}\|_{x}
  \end{align*}
\end{proposition}
\begin{proof}
  Since $T$ is compact and in view of proposition
  \ref{prop:invderdelta}, it is enough to prove the second inequality
  for $P=\Id$. In this case the lemma follows from propositions
  \ref{prop:13} and \ref{prop:12} and lemmas \ref{lemm:derandprod},
  and \ref{lemm:boundprod}. 
\end{proof}

Let $\sigma \subset C$ be the simplicial cone 
\begin{displaymath}
  \sigma =\sum_{i=1}^{N}\mathbb{R}^{+}\xi_{i}. 
\end{displaymath}
Let $\{l_{i}\}$ be the dual basis of $\{\xi_{i}\}$.
\begin{proposition}\label{prop:boundder}
  Let $\delta $ be a vector field in $T$, let $P$ be a differential operator
  that is a product of vector fields in $T$, let $(i _{j})_{j=1}^{n}$
  be a finite sequence of elements of $\{1,\dots 
  ,N\}$ and let $a\in \overline C$. Then there is a constant $K$ such
  that, for all integers 
  $1\le \alpha ,\beta \le n$, and $x\in \Int(\sigma +a)$
  \begin{align*}
   \left| \prod_{j=2}^{r} D_{\xi
   _{i_{j}}}P(D_{\xi_{i_{1}}}H_{t}.H^{-1}_{t}(x))_{\alpha ,\beta }
   \right |&\le \frac{K}{\prod_{j=1}^{r} l_{i_{j}}(x)-l_{i_{j}}(a)}\\ 
   \left| \prod_{j=1}^{r} D_{\xi
   _{i_{j} }}P(\delta H_{t}.H^{-1}_{t}(x))_{\alpha ,\beta }
   \right |&\le K. 
  \end{align*}
\end{proposition}
\begin{proof}
  The proof is as in \cite{Mumford:Hptncc} proposition 2.7, but using
  proposition \ref{prop:14} to estimate the higher derivatives. 
\end{proof}

\nnpar{End of the proof.}
Now the proof of theorem \ref{thm:17} goes exactly as the proof of
\cite{Mumford:Hptncc} theorem 3.1, but using proposition \ref{prop:16}
and proposition
\ref{prop:boundder} to bound the higher derivatives. 
\end{proof}

\begin{remark}
 Observe that what we really have proved is that if $\{e_{1},\dots
 ,e_{r}\}$ is a holomorphic frame of $E_{1}$,and $H=(h_{e_i,e_j})$ is the
 matrix of $h$ in this frame, then the entries of $H$ and $\det H^{-1}$
 are of polynomial growth in the local universal cover (which, by
 theorem \ref{thm:15}, is
 equivalent of being log forms) and that the entries of $\partial
 H\cdot H^{-1} $ are of logarithmic growth in the local universal
 cover  (which, by theorem \ref{thm:14}, is stronger than being
 log-log forms). 
\end{remark}

\subsection{Shimura varieties and automorphic vector bundles}
\label{sec:shim-vari-autom}

A wealth of examples where the theory developed in this paper can be
applied is provided by non compact Shimura varieties. In fact, the concrete
examples developed so far: modular curves (\cite{Kuehn:gainc}) and
Hilbert modular surfaces 
(\cite{BruinierBurgosKuehn:bpaihs}) are examples of Shimura
varieties of non compact type. 

For an algebraic group $G$, $G(\mathbb{R})^{+}$ is the identity
component of the topological group $G(\mathbb{R})$ and
$G(\mathbb{R})_{+}$ is the inverse image of $G^{\ad}(\mathbb{R})^{+}$
in $G(\mathbb{R})$; also $G(\mathbb{Q})^{+}=G(\mathbb{Q})\cap
G(\mathbb{R}))^{+}$ and $G(\mathbb{Q})_{+}=G(\mathbb{Q})\cap
G(\mathbb{R}))^{+}$.

\nnpar{Definition of Shimura varieties.}
Let $\underline S$ be the real algebraic torus
$\operatorname{Res}_{\mathbb{C}/\mathbb{R}} \mathbb{G}_m$.   
Following Deligne  \cite{Deligne:_variet_shimur}
(cf. \cite{Milne:_canon_shimur}) one
considers the data
\begin{itemize}
\item[(1)] $G$ a connected reductive group defined over $\mathbb{Q}$,
\item[(2)] $X$ a $G(\mathbb{R})$ -conjugacy class of morphisms $h_{x}:
  \underline S \longrightarrow
  G_\mathbb{R}$ of real algebraic groups, $x\in X$,
\end{itemize}
satisfying the properties
\begin{itemize}
\item[(a)] The Hodge structure on $\operatorname{Lie} G_\mathbb{R}$
  defined by $\Ad\circ h_{x}$ is of
  type $$\{(-1,1),(0,0),(1,-1)\},$$
\item[(b)] The involution $\operatorname{int} h_{x}(i)$ induces a Cartan
  involution on the adjoint group $G^{\operatorname{ad}}(\mathbb{R})$,
\item[(c)] Let $w:\mathbb{G}_{m,\mathbb{R}}\longrightarrow \underline
  S$ be the canonical conorm map. The weight map $h_{x}\circ w$ (whose
  image is
  central by a)) is defined
  over $\mathbb{Q}$.
\item [(d)] Let $Z'_{G}$ be the maximal $\mathbb{Q}$-split torus of
  $Z_{G}$, the center of $G$. Then
  $Z_{G}(\mathbb{R})/Z'_{G}(\mathbb{R})$ is compact.
\end{itemize}
Under the above assumptions $X$ is a product of hermitian symmetric
domains corresponding to the simple non-compact factors of 
$G^{\operatorname{ad}}(\mathbb{R})$. Denote by $\mathbb{A}^f$ the
finite ad\`eles of $\mathbb{Q}$ and let $K \subset G(\mathbb{A}^f)$ be
a neat (see for instance \cite{Pink:_arith_mixed_shimur_variet} for
the definition of neat) open compact subgroup. With these data the
Shimura variety  
$M_{K}(\mathbb{C})$ is defined by
\begin{align*}
M_{K}(\mathbb{C})= M_{K}(G,X)(\mathbb{C}) := G(\mathbb{Q})
 \backslash X \times 
  G(\mathbb{A}^f) / K.
\end{align*}

\nnpar{Connected components  of Shimura varieties.}
Let $X^{+}$ be a connected component of $X$, and for each $x\in
X^{+}$, let $h'_{x}$ be the composite of $h_{x}$ with
$G_{\mathbb{R}}\longrightarrow G_{\mathbb{R}}^{\ad}$. Then
$x\longmapsto h'_{x}$ identifies $X^{+}$ with a 
$G^{\ad}(\mathbb{R})^{+}$-conjugacy class of morphisms $
\underline S \longrightarrow
G_\mathbb{R}^{\ad}$ that satisfy the axioms of a connected Shimura
variety. In particular $X^{+}$ is a bounded symmetric domain and $X$
is a finite disjoint union of bounded symmetric domains (indexed by
$G(\mathbb{R})/G(\mathbb{R})_{+}$). 
    
Let $\mathcal{C}$ be a set of representatives of the finite set
$G(\mathbb{R})_{+}\backslash G(\mathbb{A}^{f})/K$ and, for each $g\in
\mathcal{C}$, let $\Gamma _{g}$ be the image in
$G^{\ad}(\mathbb{R})^{+}$ of the subgroup $\Gamma '_{g}=gKg^{-1}\cap
G(\mathbb{Q})_{+}$ of $G(\mathbb{Q})_{+}$. Then $\Gamma _{g}$ is a
torsion free arithmetic subgroup of $G^{\ad}(\mathbb{R})^{+}$ and
$M_{K}(\mathbb{C})$ 
is a finite disjoint union  
 \begin{displaymath}
 M_K(\mathbb{C}) = 
 \coprod_{g \in \mathcal{C}}
 \Gamma_g \backslash X^+.
 \end{displaymath}

The connected component $\Gamma_g \backslash X^+$ will be denoted $M_{\Gamma _{g}}$.



\nnpar{Algebraic models of Shimura varieties.}
Every Shimura variety is a quasi-projective
variety.  It has a
``minimal'' compactification, the Baily-Borel compactification that is
highly singular. The theory of toroidal compactifications provides us
with various other compactifications; among them we can choose non
singular ones whose boundary is a normal crossing divisor.
Moreover it has a model over a number field $E$ called the reflex
field and  
the toroidal compactifications are also defined over $E$
(\cite{Pink:_arith_mixed_shimur_variet}). This model 
can be extended 
to a proper regular model defined over $\mathcal{O}_{E}[N^{-1}]$,
where $\mathcal{O}_{E}$ is the ring of integers of $E$ and $N$ is
a suitable natural number.

\nnpar{Automorphic vector bundles}. Let $K_{x}$ be the
subgroup of $G(\mathbb{R})$ stabilizing a point $x\in X$ and
let 
$P_{x}$ be the parabolic subgroup of $G(\mathbb{C})$ arising from the
Cartan decomposition of $\Lie(G)$ associated to $K_{x}$. Let $\lambda
:K_{x}\longrightarrow GL_{n}$ be a finite dimensional representation
of $K_{x}$. It can be extended trivially to a representation
of $P_{x}$ and defines a $G(\mathbb{C})$-equivariant vector bundle
$\text{\sl \v V}$ on the compact dual $\text{\sl \v
  M}(\mathbb{C})=G(\mathbb{C})/P_{x}$. Let $\beta :X\longrightarrow
\text{\sl \v  M}(\mathbb{C})$  be the Borel embedding, then $V=\beta
^{\ast}(\text{\sl \v V})$ is a $G(\mathbb{R)}$ equivariant vector bundle
on $X$. For any neat open compact subgroup $K\subset
G(\mathbb{A}^{f})$ it defines an algebraic vector bundle
\begin{displaymath}
  V_{K}= G(\mathbb{Q})
 \backslash V \times 
  G(\mathbb{A}^f) / K
\end{displaymath}
on the Shimura variety $M_{K}$. This vector bundle is algebraic and is
defined over the reflex field $E$. Following
\cite{Harris:_bound_shimur_i}, the vector bundles obtained in this
way will be called \emph{fully decomposed automorphic vector bundles.}  

The
restriction to any component 
$M_{\Gamma _{g}}$ will be denoted by $V_{\Gamma _{g}}$. It is a fully
decomposed automorphic vector bundle in the sense of the previous section.

\nnpar{Canonical extensions.} Let $M_{K,\Sigma }$ be a smooth toroidal
compactification of $M_{K}$ and let $V_{K}$ be an
automorphic vector bundle on $M_{K}$. Then there exists a canonical
extension of $V_{K}$ to a vector bundle
$V_{K,\Sigma }$ over  $M_{K,\Sigma }$
(\cite{Mumford:Hptncc}, cf. \cite{Milne:_canon_shimur},
\cite{Harris:_funct}). This canonical 
extension can be characterized in terms of an invariant hermitian
metric on $V$.

\nnpar{} Let $M_{K}$ be a Shimura variety defined over the reflex field
$E$. Let $M_{K,\Sigma }$ be a smooth toroidal compactification of
$M_{K}$ defined over $E$ such that $D_{E}=M_{K,\Sigma }\setminus
M_{K}$ is a normal crossing divisor. Let $V_{K}$ be an
automorphic vector bundle defined over $E$ with canonical extension
$V_{K,\Sigma }$. Let $h$ be an
$G^{\operatorname{der}}(\mathbb{R})$-invariant 
hermitian metric on $V$. It induces a hermitian metric,
denoted also by $h$ on $V_{K}$. We denote also by $h$ the singular
hermitian metric induced on $V_{K,\Sigma }$.
Let $\mathcal{M}_{K,\Sigma }$ be a
regular model of $M_{K,\Sigma }$ over
$\mathcal{O}_{E}[N^{-1}]$. Assume that $V_{K,\Sigma }$ can
be extended  to a vector bundle $\mathcal{V}_{K,\Sigma
}$ over $\mathcal{M}_{K,\Sigma }$. Then theorem \ref{thm:17} implies

\begin{theorem}\label{thm:13}
  The pair $(\mathcal{V}_{K,\Sigma
}, h)$ is a log-singular hermitian vector bundle on
$\mathcal{M}_{K,\Sigma }$.\hfill $\square$
\end{theorem}

\bibliographystyle{amsplain}

\newcommand{\noopsort}[1]{} \newcommand{\printfirst}[2]{#1}
  \newcommand{\singleletter}[1]{#1} \newcommand{\switchargs}[2]{#2#1}
\providecommand{\bysame}{\leavevmode\hbox to3em{\hrulefill}\thinspace}
\providecommand{\MR}{\relax\ifhmode\unskip\space\fi MR }
\providecommand{\MRhref}[2]{%
  \href{http://www.ams.org/mathscinet-getitem?mr=#1}{#2}
}
\providecommand{\href}[2]{#2}

\end{document}